\documentclass[final]{siamltex}
\usepackage{latexsym,url,amscd,amsmath,amssymb}
\usepackage{graphics}
\usepackage{graphicx}
\usepackage{epstopdf}
\usepackage[usenames]{color}
\usepackage[hyperref]{hyperref}
\usepackage{comment}

\newtheorem{remark}{Remark}

\newcommand{\mat}[1]{\left[ \begin{array}{#1} }
\newcommand{\rix}{\end{array} \right]}

\newcommand{\R}{\mathbb{R}}
\newcommand{\B}{\mathbf{B}}
\newcommand{\mP}{\mathcal{P}}
\newcommand{\C}{\mathbb{C}}

\newcommand{\real}{\mbox{\bf Re}}

\usepackage[usenames]{color}

\begin{document}

\title{Alternating Direction Algorithms for $\ell_1$-Problems \\in Compressive Sensing}

\author{Junfeng Yang\footnotemark[1] \and Yin Zhang \footnotemark[2]}

\renewcommand{\thefootnote}{\fnsymbol{footnote}}
\footnotetext[1]{Department of Mathematics, Nanjing University, 22
Hankou Road, Nanjing, 210093, P.R. China (jfyang@nju.edu.cn)}
\footnotetext[2]{Department of Computational and Applied
Mathematics, Rice University, 6100 Main Street, MS-134, Houston,
Texas, 77005, U.S.A. (yzhang@rice.edu)}

\renewcommand{\thefootnote}{\arabic{footnote}}

\date{\today}
\maketitle

\begin{abstract}
In this paper, we propose and study the use of alternating direction
algorithms for several $\ell_1$-norm minimization problems arising
from sparse solution recovery in compressive sensing, including the
basis pursuit problem, the basis-pursuit denoising problems of both
unconstrained and constrained forms, as well as others. We present
and investigate two classes of algorithms derived from either the
primal or the dual forms of the $\ell_1$-problems. The construction
of the algorithms consists of two main steps: (1) to reformulate an
$\ell_1$-problem into one having partially separable objective
functions by adding new variables and constraints; and (2) to apply
an exact or inexact alternating direction method to the resulting
problem. The derived alternating direction algorithms can be
regarded as first-order primal-dual algorithms because both primal
and dual variables are updated at each and every iteration.
Convergence properties of these algorithms are established or
restated when they already exist. Extensive numerical results in
comparison with several state-of-the-art algorithms are given to
demonstrate that the proposed algorithms are efficient, stable and
robust.
Moreover, we present numerical results to emphasize two practically
important but perhaps overlooked points.  One point is that algorithm
speed should always be evaluated relative to appropriate solution
accuracy; another is that whenever erroneous measurements possibly
exist, the $\ell_1$-norm fidelity should be the fidelity of choice
in compressive sensing.

\end{abstract}



\pagestyle{myheadings} \thispagestyle{plain} \markboth{J.-F. Yang
and Y. Zhang}{Alternating Direction Algorithms for $\ell_1$-Problems
in Compressive Sensing}


\section{Introduction}

In the last few years, algorithms for finding sparse solutions of
underdetermined linear systems have been intensively studied,
largely because solving such problems constitutes a critical step in
an emerging methodology in digital signal processing --- compressive
sensing or sampling (CS). In CS, a digital signal is encoded as
inner products between the signal and a set of random (or
random-like) vectors where the number of such inner products, or
linear measurements, can be significantly fewer than the length of
the signal.  On the other hand, the decoding process requires
finding a sparse solution of an underdetermined linear system. What
makes such a scheme work is sparsity; i.e.,  the original signal
must have a sparse or nearly sparse representation under some known
basis. Throughout this paper we will allow all involved quantities
(signals, acquired data and encoding matrices) to be complex. Let
$\bar x\in\C^n$ be an original signal that we wish to capture.
Without loss of generality, we assume that $\bar x$ is sparse under
the canonical basis, i.e., the number of nonzero components in $\bar
x$, denoted by $\|\bar x\|_0$, is far fewer than its length. Instead
of sampling $\bar x$ directly, in CS one first obtains a set of
linear measurements
\begin{eqnarray}\label{encoding}
b = A\bar x \in \C^m,
\end{eqnarray}
where $A \in \C^{m \times n}$ ($m<n$) is an encoding matrix. The
original signal $\bar x$ is then reconstructed from the
underdetermined linear system $Ax=b$ via certain reconstruction
technique.  Basic CS theory presented in
\cite{Candes-Romberg-Tao-05, Candes-Romberg-Tao-06, Donoho06} states
that it is extremely probable to reconstruct $\bar x$ accurately or
even exactly from $b$ provided that $\bar x$ is sufficiently sparse
(or nearly sparse) relative to the number of measurements, and the
encoding matrix $A$ possesses certain desirable attributes.

In the rest of this section, we briefly review the
essential ingredients of CS decoding process
and some existing methods for the relevant
optimization problems, summarize our main contributions in this
paper, and describe the notation and organization of the paper.

\subsection{Signal Decoding in CS}

To make CS successful, two ingredients must be addressed carefully.
First, a sensing matrix $A$ must be designed so that the compressed
measurement $b=A\bar x$ contains enough information for a successful
recovery of $\bar x$. Second, an efficient, stable and robust
reconstruction algorithm must be available for recovering $\bar x$
from $A$ and $b$. In the present paper, we will only concentrate on the
second aspect.

In order to recover the sparse signal $\bar x$ from the
underdetermined system \eqref{encoding}, one could naturally
consider seeking among all solutions of \eqref{encoding} the
sparsest one, i.e., solving
\begin{eqnarray}\label{Decoder-L0}
\min_{x\in \C^n}\{\|x\|_0: Ax=b\}
\end{eqnarray}
where $\|x\|_0$ is the number of nonzeros in $x$. Indeed, with
overwhelming probability decoder \eqref{Decoder-L0} can recover
sparse signals exactly from a very limited number of random
measurements (see e.g., \cite{decoder-l0}). Unfortunately, this
$\ell_0$-problem is combinatorial and generally computationally
intractable.   A fundamental decoding model in CS is the so-called
basis pursuit (BP) problem \cite{BP}:
\begin{equation}\label{decoder-L1}
\min_{x \in \C^n} \{\|x\|_1 : A x = b\}.
\end{equation}
Minimizing the $\ell_1$-norm in \eqref{decoder-L1} plays a central
role in promoting solution sparsity. In fact, problem
\eqref{decoder-L1} shares common solutions with \eqref{Decoder-L0}
under some favorable conditions (see, for example, \cite{l1eql0}).
When $b$ contains noise, or when $\bar x$ is not exactly sparse but
only compressible, as are the cases in most practical applications,
certain relaxation to the equality constraint in \eqref{decoder-L1}
is desirable. In such situations, common relaxations to
\eqref{decoder-L1} include the constrained basis pursuit denoising
(BP$_\delta$) problem \cite{BP}:
\begin{equation}\label{decoder-BPDN}
\min_{x \in \C^n} \{\|x\|_1 : \|A x - b\|_2 \leq \delta\},
\end{equation}
and its variants including the unconstrained basis pursuit denoising
(QP$_\mu$) problem
\begin{eqnarray}\label{decoder-L1L2}
\min_{x  \in \C^n} \|x\|_1 + \frac{1}{2\mu}\|Ax-b\|_2^2,
\end{eqnarray}
where $\delta, \mu>0$ are parameters. From optimization theory,
it ie well known that
problems \eqref{decoder-BPDN} and \eqref{decoder-L1L2} are
equivalent in the sense that solving one will determine a parameter
value in the other so that the two share the same solution. As
$\delta$ and $\mu$ approach zero, both BP$_\delta$ and QP$_\mu$
converge to \eqref{decoder-L1}. In this paper, we also consider the
use of an $\ell_1/\ell_1$ model of the form
\begin{eqnarray}\label{decoder-L1L1}
\min_{x \in \C^n} \|x\|_1 + \frac{1}{\nu}\|Ax-b\|_1,
\end{eqnarray}
whenever $b$ might contain erroneous measurements. It is well-known
that unlike \eqref{decoder-L1L2} where squared $\ell_2$-norm
fidelity is used, the $\ell_1$-norm fidelity term makes
\eqref{decoder-L1L1} an exact penalty method in the sense that it
reduces to  \eqref{decoder-L1} when $\nu>0$ is less than some
threshold.

It is worth noting that problems \eqref{decoder-L1},  \eqref{decoder-BPDN},
\eqref{decoder-L1L2} and \eqref{decoder-L1L1} all have their
``nonnegative counterparts" where the signal $x$ is real and
nonnegative.  These nonnegative counterparts will be briefly
considered later. Finally, we mention that aside from
$\ell_1$-related decoders, there exist alternative decoding
techniques such as greedy algorithms (e.g., \cite{Tropp07}) which,
however, are not a subject of concern in this paper.

\subsection{Some existing methods}

In the last few years, quite a number of algorithms have been proposed
and studied for solving the aforementioned $\ell_1$-problems arising in
CS.   Although these problems are convex programs with relatively
simple structures (e.g., the basis pursuit problem is a linear program
when $x$ is real), they do demand dedicated algorithms because
standard methods, such as interior-point algorithms for linear and
quadratic programming, are simply too inefficient on them.
This is the consequence of several factors, most prominently the
fact that the data matrix $A$ is totally dense while the solution is
sparse.   Clearly, the existing standard algorithms were not designed to
handle such a feature.
Another noteworthy structure is that encoding matrices in CS are
often formed by randomly taking a subset of rows from orthonormal
transform matrices, such as DCT (discrete cosine transform), DFT
(discrete Fourier transform) or DWHT (discrete Walsh-Hadamard
transform) matrices.   Such encoding matrices do not require storage
and enable fast matrix-vector multiplications. As a result,
first-order algorithms that are able to take advantage of such a
special feature lead to better performance and are highly desirable.

One of the earliest first-order methods for solving
\eqref{decoder-L1L2} is the gradient projection method suggested in
\cite{GPSR}, where the authors reformulated \eqref{decoder-L1L2} as
a box-constrained quadratic program and implemented a gradient
projection method with line search. To date, the most widely studied
first-order method for solving \eqref{decoder-L1L2} is the iterative
shrinkage/thresholding (IST) method, which is first proposed in
\cite{Fig-Nowak-03, Nowak-Fig-01, DeMol-Defrise02} for wavelet-based
image deconvolution and then independently discovered and analyzed
by many others \cite{Elad06, Starck-Candes-03, Starck-etal-03,
Daubechies-04}. In \cite{FPC, FPC2}, Hale, Yin and Zhang derived the
IST algorithm from an operator splitting framework and combined it with
a continuation strategy. The resulting
algorithm, which is named fixed-point continuation (FPC), is also
accelerated via a non-monotone line search
with Barzilai-Borwein steplength \cite{BB88}. A similar sparse
reconstruction algorithm called SpaRSA was  also studied by Wright,
Nowak and Figueiredo in \cite{SpaRSA}. Recently, Beck and Teboulle
proposed a fast IST algorithm (FISTA) in \cite{FISTA}, which attains
the same optimal convergence in function values as Nesterov's
multi-step gradient method \cite{Nesterov-07} for minimizing
composite functions. Lately, Yun and Toh also studied a block
coordinate gradient descent (CGD) method in \cite{CGD} for solving
\eqref{decoder-L1L2}.

There exist also algorithms for solving constrained
$\ell_1$-problems \eqref{decoder-L1} and \eqref{decoder-BPDN}. The
Bregman iteration \cite{BregmanTV} was applied to the basis pursuit
problem in \cite{BregmanL1}. In the same paper, a   linearized
Bregman method was also suggested and analyzed subsequently in
\cite{LBregmanConver1, LBregmanConver2, LBregmanDual}. In
\cite{SPGL1}, Friedlander and Van den Berg proposed a spectral
projection gradient method (SPGL1), where \eqref{decoder-BPDN} is
solved by a root-finding framework applied to a sequence of LASSO
problems \cite{LASSO}. Moreover, based on a smoothing technique
studied in \cite{Nesterov-smoothing}, a fast and accurate
first-order algorithm called NESTA  was proposed in \cite{NESTA} for
solving \eqref{decoder-BPDN}.

In Section \ref{sc:numerical}, we present extensive comparison
results with several state-of-the-art algorithms including FPC,
SpaRSA, FISTA and CGD for solving \eqref{decoder-L1L2}, and SPGL1
and NESTA for solving \eqref{decoder-L1} and \eqref{decoder-BPDN}.

\subsection{Contributions}
After years of intensive research on $\ell_1$-problem solving,
it would appear that most relevant algorithmic ideas have been either tried or,
in many cases, re-discovered.   Yet interestingly, until the writing of the present
paper,  the classic idea of alternating direction method (ADM) had not,
to the best of our knowledge, been seriously investigated.

The {\em main contributions of this paper} are to introduce ADM
technique to the area of solving $\ell_1$-problems in CS and other
applications, and to demonstrate its usefulness as a versatile and
powerful algorithmic tool. Indeed, based on ADM technique we have
derived first-order primal-dual algorithms for  \eqref{decoder-L1},
\eqref{decoder-BPDN}, \eqref{decoder-L1L2} and \eqref{decoder-L1L1}.
Furthermore, the versatility of ADM approach has allowed us to
develop a Matlab package called {\tt YALL1} \cite{YALL1} that can
solve eight different $\ell_1$-minimization models, i.e.,
\eqref{decoder-L1}, \eqref{decoder-BPDN}, \eqref{decoder-L1L2},
\eqref{decoder-L1L1} and their nonnegative counterparts, where local
weights are also permitted in the $\ell_1$-norm.

In this paper, we present extensive computational results to
document the numerical performance of the proposed algorithms in
comparison to several state-of-the-art algorithms for solving
$\ell_1$-problems under various situations. As by-products, we also
address a couple of issues of practical importance; i.e., choices of
optimization models and proper evaluation of algorithm speed.

\subsection{Notation}
We let $\|\cdot\|$ be the $\ell_2$-norm and
$\mathcal{P}_{\Omega}(\cdot)$ be a projection operator onto a convex
set $\Omega$ under the $\ell_2$-norm. Superscripts ``$\top$'' and
``$*$'' denote, respectively, the transpose and the conjugate
transpose operators for real and complex quantities. We let
$\real(\cdot)$ and $|\cdot|$ be, respectively, the real part and the
magnitude of a complex quantity, which are applied component-wise to
complex vectors. Further notation will be introduced  wherever  it
occurs.

\subsection{Organization}
This paper is organized as follows. In Section~2, we first review
the basic idea of the classic ADM technique and then derive
alternating direction algorithms for solving \eqref{decoder-L1},
\eqref{decoder-BPDN} and \eqref{decoder-L1L2}. We also establish
convergence of the primal-based algorithms, while that of the
dual-based algorithms follows from classic results in the
literature.  In Section~3, we present numerical results to compare
the behavior of model \eqref{decoder-L1L1} to that of models
\eqref{decoder-BPDN} and \eqref{decoder-L1L2} under various
scenarios of noise. In Section~4, we first re-emphasize the
sometimes overlooked common sense on appropriate evaluations of
algorithm speed, and then present extensive comparison results with
several state-of-the-art algorithms to demonstrate the performance
of the proposed algorithms.    Finally, we conclude the paper in
Section 5 and discuss several extensions of ADM approach to other
$\ell_1$-like problems.

\section{ADM-based first-order primal-dual algorithms}

In this section, based on the classic ADM technique, we propose
first-order primal-dual algorithms that update both primal and dual
variables at each iteration for the solution of $\ell_1$-problems.
We start with a brief review on a general framework of ADM.

\subsection{General framework of ADM}\label{sc:generalADM}
Let $f(x): \R^m\rightarrow \R$ and $g(y): \R^n\rightarrow \R$ be
convex functions, $A\in\R^{p\times m}$, $B\in\R^{p\times n}$ and
$b\in\R^p$. We consider the structured optimization problem
\begin{eqnarray}\label{ADM-P}
\min_{x,y}\left\{ f(x) + g(y): Ax + By = b\right\},
\end{eqnarray}
where variables $x$ and $y$ are separate in the objective, and
coupled only in the constraint.   The augmented Lagrangian function
of this problem is given by
\begin{eqnarray}\label{ADM-P-AL}
\mathcal{L}_{\mathcal{A}}(x,y,\lambda) = f(x) + g(y) - \lambda^\top
(Ax + By - b) + \frac{\beta}{2}\|Ax+By-b\|^2,
\end{eqnarray}
where  $\lambda\in\R^p$ is the Lagrangian multiplier and $\beta>0$
is a penalty parameter. The classic augmented Lagrangian method
\cite{Hestenes69, Powell69} iterates as follows: given
$\lambda^k\in\R^p$,
\begin{eqnarray}\label{ALM}
\left\{
  \begin{array}{l}
  (x^{k+1},y^{k+1})\gets \arg\min_{x,y} \mathcal{L}_{\mathcal{A}}(x,y,\lambda^k),\\
  \lambda^{k+1} \gets \lambda^k - \gamma\beta(Ax^{k+1} +
  By^{k+1}-b),
  \end{array}
\right.
\end{eqnarray}
where $\gamma\in(0,2)$ guarantees convergence, as long as
the subproblem is solved to an increasingly high accuracy at every
iteration \cite{Rockafellar73}. However, an accurate, joint
minimization with respect to $(x,y)$ can become costly without
taking advantage of the separable form of the objective function
$f(x)+g(y)$. In contrast, ADM utilizes the separability structure in
\eqref{ADM-P} and replaces the joint minimization by two simpler
subproblems.  Specifically, ADM minimizes
$\mathcal{L}_{\mathcal{A}}(x,y,\lambda)$ with respect to $x$ and $y$
separately via a Gauss-Seidel type iteration. After just one sweep
of alternating minimization with respect to $x$ and $y$, the
multiplier $\lambda$ is updated immediately.   In short, given
$(y^k,\lambda^k)$,  ADM iterates as follows
\begin{eqnarray}\label{GADM}
\left\{
  \begin{array}{ll}
  x^{k+1} \gets\arg\min_{x} \mathcal{L}_{\mathcal{A}}(x,y^k,\lambda^k),\\
   y^{k+1} \gets\arg\min_{y} \mathcal{L}_{\mathcal{A}}(x^{k+1},y,\lambda^k),\\
  \lambda^{k+1} \gets \lambda^k - \gamma\beta(Ax^{k+1} +
  By^{k+1}-b).
  \end{array}
\right.
\end{eqnarray}

The basic idea of ADM goes back to the work of Glowinski and
Marocco \cite{Glow75} and Gabay and Mercier \cite{Gabay76}. Let
$\theta_1(\cdot)$ and $\theta_2(\cdot)$ be convex functionals, and
$A$ be a continuous linear operator. The authors of \cite{Gabay76}
considered minimizing an energy function of the form
\begin{eqnarray*}\label{Gabay-Prob}
\min_u \theta_1(u) + \theta_2(Au).
\end{eqnarray*}
By introducing an auxiliary variable $v$, the above problem 
was equivalently transformed to
\begin{eqnarray*}\label{Gabay-Prob2}
\min_{u,v} \left\{\theta_1(u) + \theta_2(v): Au - v = 0\right\},
\end{eqnarray*}
which has the form of \eqref{ADM-P} and to which the ADM approach was applied.
Subsequently, ADM was studied extensively in optimization and
variational analysis. In \cite{Glow89}, ADM is interpreted as the
Douglas-Rachford splitting method \cite{Douglas56} applied to a dual
problem. The equivalence between ADM and a proximal point method is
shown in \cite{Eckstein92}. ADM is also studied in convex
programming \cite{Fukushima92} and variational inequalities
\cite{Tseng91, He-Yang-98}. Moreover, ADM has been extended to
allowing inexact subproblem minimization \cite{Eckstein92, He-MP02}.

In \eqref{GADM}, a steplength $\gamma>0$ is attached to the update
of $\lambda$. Under certain technical assumptions, convergence of
ADM with a steplength $\gamma\in(0,(\sqrt{5}+1)/2)$ was established
in \cite{Glow84, Glow89} in the context of variational inequality.
The shrinkage in the permitted range from $(0,2)$ in the augmented
Lagrangian method to $(0,(\sqrt{5}+1)/2)$ in ADM is related to
relaxing the exact minimization of
$\mathcal{L}_{\mathcal{A}}(x,y,\lambda^k)$ with respect to $(x,y)$
to merely one round of alternating minimization. Recently, ADM has
been applied to total variation based  image restoration and
reconstruction in \cite{ADM09, RecPF}. In the following,
we apply ADM technique to \eqref{decoder-BPDN} and
\eqref{decoder-L1L2}, while the application to \eqref{decoder-L1}
and \eqref{decoder-L1L1} will be a by-product.

\subsection{Applying ADM to primal problems}\label{sec:p-alg}
In this subsection, we apply ADM to primal $\ell_1$-problems
\eqref{decoder-BPDN} and \eqref{decoder-L1L2}.  First, we introduce
auxiliary variables to reformulate these problems into the form of
\eqref{ADM-P}.   Then, we apply alternating minimization to the
corresponding augmented Lagrangian functions, either exactly or
approximately,  to obtain ADM-like algorithms.

With an auxiliary variable $r\in\C^m$, problem
\eqref{decoder-L1L2} is clearly equivalent to
\begin{eqnarray}\label{decoder-L1L2-2}
\min_{x\in\C^n, \; r\in\C^m} \left\{f_p(r,x)\triangleq \|x\|_1 +
\frac{1}{2\mu}\|r\|^2: Ax + r = b\right\},
\end{eqnarray}
that has an augmented Lagrangian subproblem of the form
\begin{eqnarray}\label{aug-decoder-L1L2-2}
\min_{x\in\C^n, \; r\in\C^m} \left\{\|x\|_1 + \frac{1}{2\mu}\|r\|^2
-\real( y^* (Ax + r - b)) + \frac{\beta}{2}\|Ax + r - b\|^2\right\},
\end{eqnarray}
where   $y\in\C^m$ is a multiplier and $\beta>0$ is a penalty
parameter. Given $(x^k,y^k)$, we obtain $(r^{k+1},x^{k+1},y^{k+1})$
by applying alternating minimization to \eqref{aug-decoder-L1L2-2}.
First, it is easy to show that, for $x=x^k$ and $y=y^k$ fixed, the
minimizer of \eqref{aug-decoder-L1L2-2} with respect to $r$ is given
by
\begin{eqnarray}\label{iter-s}
r^{k+1} = \frac{\mu\beta}{1+ \mu\beta} \left(y^k/\beta-(Ax^k -
b)\right).
\end{eqnarray}
Second, for $r=r^{k+1}$ and $y=y^k$ fixed, simple manipulation shows
that the minimization of \eqref{aug-decoder-L1L2-2} with respect to
$x$ is equivalent to
\begin{eqnarray}\label{aug-decoder-L1L2-fixeds}
\min_{x\in\C^n} \|x\|_1 + \frac{\beta}{2}\|Ax + r^{k+1} - b -
y^k/\beta\|^2,
\end{eqnarray}
which itself is  in the form of \eqref{decoder-L1L2}. However,
instead of solving \eqref{aug-decoder-L1L2-fixeds} exactly, we
approximate it by
\begin{eqnarray}\label{aug-decoder-L1L2-fixeds-Lzation}
\min_{x\in\C^n} \|x\|_1 + \beta \left((g^k)^*(x-x^k) +
\frac{1}{2\tau}\|x-x^k\|^2\right),
\end{eqnarray}
where $\tau>0$ is a proximal parameter and
\begin{eqnarray}\label{def:gt}
g^k \triangleq A^*(Ax^k + r^{k+1} - b - y^k/\beta)
\end{eqnarray}
is the gradient of the quadratic term in \eqref{aug-decoder-L1L2-fixeds},
$\beta$ not included, at $x=x^k$.
The solution of \eqref{aug-decoder-L1L2-fixeds-Lzation} is given
explicitly by
\begin{eqnarray}\label{iter-x}
x^{k+1}  = \text{Shrink}\left(x^k - \tau g^k,
\frac{\tau}{\beta}\right) \triangleq \max\left\{|x^k - \tau g^k| -
\frac{\tau}{\beta},0\right\}\circ\text{sign}(x^k - \tau  g^k),
\end{eqnarray}
where ``$\circ$" represents component-wise multiplication.
The operation defined in \eqref{iter-x} is well-known as the
one-dimensional shrinkage (or soft thresholding).
Finally,  we update the multiplier $y$ by
\begin{eqnarray}\label{iter-y}
y^{k+1} = y^k - \gamma\beta (A x^{k+1} + r^{k+1} - b),
\end{eqnarray}
where $\gamma>0$ is a constant. In short, ADM applied to
\eqref{decoder-L1L2} produces the iteration:
\begin{eqnarray}\label{alg-primal}
\left\{
  \begin{array}{ll}
    r^{k+1}  = \frac{\mu\beta}{1+ \mu\beta} \left(y^k/\beta - (Ax^k - b)\right),  \\
    x^{k+1} = \text{Shrink}\left(x^k - \tau g^k, \frac{\tau}{\beta}\right), \\
    y^{k+1} = y^k - \gamma\beta (A x^{k+1} + r^{k+1} - b).
  \end{array}
\right.
\end{eqnarray}
We note that \eqref{alg-primal} is an inexact ADM because the
$x$-subproblem is solved approximately. The convergence of
\eqref{alg-primal} is not covered by the analysis given in
\cite{Eckstein92}  where each ADM subproblem is required to be
solved more and more accurately as the algorithm proceeds. On the
other hand, the analysis in \cite{He-MP02} does cover the
convergence of \eqref{alg-primal} but only for the case $\gamma=1$.
A more general convergence result for \eqref{alg-primal} is
established below that allows $\gamma > 1$.

\begin{theorem}\label{conv-primal}
Let $\tau, \gamma>0$ satisfy $\tau\lambda_{\max} + \gamma < 2$,
where $\lambda_{\max}$ denotes the maximum eigenvalue of $A^*A$. For
any fixed $\beta>0$, the sequence $\{(r^k,x^k,y^k)\}$ generated by
\eqref{alg-primal} from any starting point $(x^0,y^0)$ converges
to $(\tilde{r},\tilde{x},\tilde{y})$, where $(\tilde{r},\tilde{x})$
is a solution of \eqref{decoder-L1L2-2}.
\end{theorem}
\begin{proof}
The proof is given in the Appendix.
\end{proof}

A similar alternating minimization idea can also be applied to problem
\eqref{decoder-BPDN}, which is equivalent to
\begin{eqnarray}\label{ConL1L2-aug-r}
\min_{x\in\C^n,\;r\in\C^m}\left\{\|x\|_1: Ax+r=b,
\|r\|\leq\delta\right\},
\end{eqnarray}
and has an augmented Lagrangian subproblem of the form
\begin{eqnarray}\label{ALF-ConL1L2}
\min_{x\in\C^n,\;r\in\C^m}\left\{\|x\|_1-y^*(Ax+r-b) +
\frac{\beta}{2}\|Ax+r-b\|^2: \|r\|\leq\delta\right\}.
\end{eqnarray}
Similar to the derivation of \eqref{alg-primal}, applying inexact
alternating minimization to \eqref{ALF-ConL1L2} yields the following
iteration scheme:
\begin{eqnarray}\label{alg-primal-Con}
\left\{
  \begin{array}{ll}
    r^{k+1} = \mathcal{P}_{B_\delta}\left(y^k/\beta-(Ax^k-b)\right), \\
    x^{k+1} = \text{Shrink}(x^k-\tau g^k, \tau/\beta), \\
    y^{k+1} = y^k - \gamma\beta (Ax^{k+1} + r^{k+1} - b),
  \end{array}
\right.
\end{eqnarray}
where $g^k$ is as defined in \eqref{def:gt}, and $\mathcal{P}_{B_\delta}$ is
the projection onto the set $B_\delta \triangleq \{\xi\in\C^m: \|\xi\|\leq \delta\}$
(in Euclidean norm).  This algorithm also has a similar convergence result as
\eqref{alg-primal}.

\begin{theorem}\label{conv-primal-Con}
Let $\tau, \gamma>0$ satisfy $\tau\lambda_{\max} + \gamma < 2$,
where $\lambda_{\max}$ denotes the maximum eigenvalue of $A^*A$. For
any fixed $\beta>0$, the sequence $\{(r^k,x^k,y^k)\}$ generated by
\eqref{alg-primal-Con} from any starting point $(x^0,y^0)$
converges to $(\tilde{r},\tilde{x},\tilde{y})$, where
$(\tilde{r},\tilde{x})$ solves \eqref{ConL1L2-aug-r}.
\end{theorem}
\begin{proof}
The proof of Theorem \ref{conv-primal-Con} is similar to that of
Theorem \ref{conv-primal}, and thus is omitted.
\end{proof}

\begin{remark}
We point out that, when $\mu=\delta=0$, both
\eqref{alg-primal} and \eqref{alg-primal-Con} reduce to
\begin{eqnarray}\label{alg-primal2BP}
\left\{
  \begin{array}{ll}
     x^{k+1} = \text{Shrink}\left(x^k - \tau A^*(Ax^k - b - y^k/\beta), \tau/\beta\right), \\
     y^{k+1} = y^k - \gamma\beta(Ax^{k+1}  - b),
  \end{array}
\right.
\end{eqnarray}
which is an algorithm for solving the basis pursuit problem \eqref{decoder-L1}.
Using similar techniques as in the proof for Theorem \ref{conv-primal}, we can
establish the convergence of \eqref{alg-primal2BP} to a solution of
\eqref{decoder-L1}. The only difference between
\eqref{alg-primal2BP} and the linearized Bregman method proposed in
\cite{BregmanL1} lies in the updating of multiplier. The advantage
of \eqref{alg-primal2BP} is that it solves \eqref{decoder-L1}, while
the linearized Bregman method solves a penalty approximation of
\eqref{decoder-L1} (see e.g., \cite{LBregmanDual}).
\end{remark}

Since we applied the ADM idea to the primal problems \eqref{decoder-L1},
\eqref{decoder-BPDN} and \eqref{decoder-L1L2}, we name the resulting
algorithms \eqref{alg-primal},  \eqref{alg-primal-Con} and \eqref{alg-primal2BP}
primal-based ADMs or PADMs in short.   In fact, these algorithms are really
of primal-dual nature because both the primal and the dual variables are updated
at each and every iteration.  In addition, these are obviously first-order algorithms.

\subsection{Applying  ADM to dual problems}\label{sc:DADM}
Similarly,  we can apply the ADM idea to the dual problems of
\eqref{decoder-BPDN} and \eqref{decoder-L1L2}, which results to
equally simple yet more efficient algorithms. First, we assume that
the rows of $A$ are orthonormal, i.e., $A$ satisfies $AA^*=I$. We
will remark on how to extend the treatment to a non-orthonormal matrix
$A$ at the end of this section.

Simple computation shows that the dual of \eqref{decoder-L1L2} or
equivalently \eqref{decoder-L1L2-2} is given by
\begin{eqnarray}\label{dualL1L2}
\nonumber &&\max_{y\in\C^m}\min_{x\in\C^n,r\in\C^m}\left\{ \|x\|_1 +
\frac{1}{2\mu}\|r\|^2 -
\real(y^*(Ax+r-b))\right\}\\
\nonumber &=&\max_{y\in\C^m} \left\{\real(b^*y) -\frac{\mu}{2}
\|y\|^2+ \min_{x\in\C^n}
\left(\|x\|_1 - \real(y^*Ax)\right) + \frac{1}{2\mu} \min_{r\in\C^m}\|r-\mu y\|^2\right\}\\
&=&\max_{y\in\C^m}\left\{\real(b^*y) - \frac{\mu}{2}\|y\|^2: A^*y
\in \B_1^\infty\right\},
\end{eqnarray}
where $\B_1^\infty\triangleq \{\xi \in \C^n: \|\xi\|_\infty \le
1\}$. By introducing $z\in\C^n$, \eqref{dualL1L2} is equivalently
transformed to
\begin{eqnarray}\label{dualL1L2-2}
\max_{y\in\C^m}\left\{f_d(y)\triangleq\real(b^*y) -
\frac{\mu}{2}\|y\|^2: z-A^*y = 0, z\in \B_1^\infty\right\},
\end{eqnarray}
which has an augmented Lagrangian subproblem of the form
\begin{eqnarray}\label{aug-dualL1L2}
\min_{y\in\C^m, z\in\C^n}\left\{-\real(b^*y) + \frac{\mu}{2}\|y\|^2
-\real(x^* (z-A^*y)) + \frac{\beta}{2}\|z-A^*y\|^2, z\in
\B_1^\infty\right\},
\end{eqnarray}
where $x\in\C^n$ is a multiplier (in fact, the primal variable)
and $\beta>0$ is a penalty
parameter.  Now we apply the ADM scheme to \eqref{dualL1L2-2},
i.e., alternatingly update the multiplier (or the primal variable) $x
\in \C^n$ and the dual variables $y\in \C^m$ and $z\in\C^n$. First,
it is easy to show that, for  $x=x^k$ and $y=y^k$, the minimizer
$z^{k+1}$ of \eqref{aug-dualL1L2} with respect to $z$ is given
explicitly by
\begin{eqnarray}\label{update-z}
z^{k+1} = \mP_{\B_1^\infty} (A^*y^k + x^k/\beta),
\end{eqnarray}
where, as in the rest of the paper, $\mP$ represent a projection
(in Euclidean norm) onto a convex set denoted as a subscript.
Second, for $x=x^k$ and $z=z^{k+1}$, the minimization of
\eqref{aug-dualL1L2} with respect to $y$ is a least squares problem
 and the corresponding normal equations are
\begin{eqnarray}\label{y-opt}
(\mu I + \beta AA^*)y = \beta Az^{k+1} - (Ax^k-b).
\end{eqnarray}
Under the assumption  $AA^*=I$, the solution $y^{k+1}$ of
\eqref{y-opt} is given by
\begin{eqnarray}\label{update-y}
y^{k+1} = \frac{\beta}{\mu + \beta}\left( Az^{k+1}
-(Ax^k-b)/\beta\right).
\end{eqnarray}
Finally, we update $x$ as follows
\begin{eqnarray}\label{update-x}
x^{k+1} = x^k - \gamma \beta (z^{k+1} - A^*y^{k+1}),
\end{eqnarray}
where $\gamma\in(0,(\sqrt{5}+1)/2)$. Thus, the ADM scheme
for \eqref{dualL1L2-2} is as follows:
\begin{eqnarray}\label{alg-pd}
\left\{
  \begin{array}{ll}
    z^{k+1} = \mP_{\B_1^\infty} (A^*y^{k} + x^k/\beta),  \\
    y^{k+1} = \frac{\beta}{\mu + \beta}\left( Az^{k+1}-(Ax^k-b)/\beta\right), \\
    x^{k+1} = x^k - \gamma \beta (z^{k+1} - A^*y^{k+1}).
  \end{array}
\right.
\end{eqnarray}
Similarly, the ADM technique can also be applied to the dual of
\eqref{decoder-BPDN} given by
\begin{eqnarray}\label{dualConL1L2}
\max_{y\in\C^m}\left\{b^* y-\delta\|y\|: A^* y \in
\B_1^\infty\right\},
\end{eqnarray}
and produces the iteration scheme
\begin{eqnarray}\label{alg-pd-Con}
\left\{
  \begin{array}{ll}
    z^{k+1} = \mathcal{P}_{\B^\infty_1}(A^* y^{k}+x^k/\beta), \\
    y^{k+1} = \mathcal{S}\left(Az^{k+1}-(Ax^k-b)/\beta,\delta/\beta\right), \\
    x^{k+1} = x^k - \gamma\beta (z^{k+1} - A^* y^{k+1}),
  \end{array}
\right.
\end{eqnarray}
where  $\mathcal{S}(v,\delta/\beta)\triangleq v -
\mP_{\B_{\delta/\beta}}(v)$ with $\B_{\delta/\beta}$ being the
Euclidian ball in $\C^m$ with radius $\delta/\beta$.

\begin{remark}
Under the assumption $AA^*=I$,  \eqref{alg-pd} is an
exact ADM in the sense that each subproblem is solved exactly. From
convergence results in \cite{Glow84, Glow89}, for any $\beta>0$ and
$\gamma\in(0,(\sqrt{5}+1)/2)$, the sequence $\{(x^k,y^k,z^k)\}$
generated by \eqref{alg-pd} from any starting point $(x^0,y^0)$
converges to $(\tilde{x},\tilde{y},\tilde{z})$, which solves the
primal-dual pair \eqref{decoder-L1L2} and \eqref{dualL1L2-2}.
Similar arguments apply to \eqref{alg-pd-Con} and the primal-dual
pair \eqref{decoder-BPDN} and \eqref{dualConL1L2}.
\end{remark}

Derived from the dual problems, we name the algorithms \eqref{alg-pd}
and \eqref{alg-pd-Con} dual-based ADMs or simply DADMs.  Again we note
these are in fact first-order primal-dual algorithms.

It is easy to show that the dual of \eqref{decoder-L1} is given by
\begin{eqnarray}\label{dualBP}
\max_{y\in\C^m}\left\{\real(b^*y): A^*y \in \B_1^\infty\right\},
\end{eqnarray}
which is a special case of \eqref{dualL1L2} and \eqref{dualConL1L2}
with $\mu=\delta=0$. Therefore, both \eqref{alg-pd} and
\eqref{alg-pd-Con} can be applied to solve \eqref{decoder-L1}.
Specifically, when $\mu=\delta=0$, both \eqref{alg-pd} and
\eqref{alg-pd-Con} reduce to
\begin{eqnarray}\label{DADM-BP}
\left\{
  \begin{array}{ll}
    z^{k+1} = \mathcal{P}_{B^\infty_1}(A^* y^{k}+x^k/\beta), \\
    y^{k+1} = Az^{k+1}-(Ax^k-b)/\beta, \\
    x^{k+1} = x^k - \gamma\beta (z^{k+1} - A^* y^{k+1}).
  \end{array}
\right.
\end{eqnarray}

We note that the last equality in \eqref{dualL1L2} holds if and only
if $r=\mu y$. Therefore, the primal-dual residues and the duality
gap between \eqref{decoder-L1L2-2} and \eqref{dualL1L2-2} can be
defined by
\begin{eqnarray}
\left\{
  \begin{array}{ll}
    r_p \triangleq Ax + r - b \equiv Ax + \mu y -b,  \\
    r_d \triangleq A^*y - z,  \\
    \Delta \triangleq f_d(y) - f_p(x,r)\equiv \real(b^*y) -
\mu\|y\|^2 - \|x\|_1.
  \end{array}
\right.
\end{eqnarray}
In computation,  algorithm \eqref{alg-pd} can be terminated by
\begin{eqnarray}\label{def-OptCond}
    \text{Res}\triangleq \left\{\|r_p\|/\|b\|,\; \|r_d\|/\sqrt{m},\; \Delta/f_p(x,r)\right\}\leq\epsilon
\end{eqnarray}
where $\epsilon>0$ measures residue in optimality.

When $AA^*\neq I$, the solution of \eqref{y-opt}  is costly.  In
this case, we take a steepest descent step in the $y$ direction and
obtain the following iteration scheme:
\begin{eqnarray}\label{alg-pd-nonorth}
\left\{
  \begin{array}{ll}
    z^{k+1} =  \mP_{\B_1^\infty} (A^* y^k + x^k/\beta),  \\
    y^{k+1} = y^k - \alpha_k^*g^k, \\
    x^{k+1} = x^k - \gamma \beta (z^{k+1}- A^*y^{k+1}),
  \end{array}
\right.
\end{eqnarray}
where, at the current point $(x^k,y^k,z^k)$, $g^k$ and $\alpha_k^*$
are given by
\[
g^k = \mu y^k + Ax^k - b + \beta
A(A^*y^k-z^{k+1})\quad\text{and}\quad
\alpha_k^*=\frac{(g^k)^*g^k}{(g^k)^*\left(\mu I + \beta
AA^*\right)g^k}.
\]
In our experiments, algorithm \eqref{alg-pd-nonorth} converges very
well for random  matrices where $AA^*\neq I$, although its
convergence remains an issue of further research. Similar arguments
apply to \eqref{dualConL1L2}.

The ADM idea can also be easily applied to $\ell_1$-problems for
recovering real and nonnegative signals.   As an example, we consider
model \eqref{decoder-L1L2} plus nonnegativity constraints:
\begin{eqnarray}\label{decoder-L1L2-pos}
\min_{x  \in \R^n} \left\{\|x\|_1 + \frac{1}{2\mu}\|Ax-b\| ^2: x\geq
0\right\},
\end{eqnarray}
where $(A,b)$ can remain complex, e.g., $A$ being a partial  Fourier matrix.
A similar derivation as for \eqref{dualL1L2} shows that a dual problem of
\eqref{decoder-L1L2-pos} 
is equivalent to
\begin{eqnarray}\label{dualL1L2-pos}
\max_{y\in\C^m}\left\{\real(b^*y) - \frac{\mu}{2}\|y\|^2: z - A^*y =
0, z\in\mathcal{F}\right\},
\end{eqnarray}
where $\mathcal{F} \triangleq \{z\in\C^n: \real(z) \leq 1\}$. The
only difference between \eqref{dualL1L2-pos} and \eqref{dualL1L2-2}
lies in the changing of constraints on $z$ from $z\in\B^\infty_1$ to
$z\in\mathcal{F}$. Applying the ADM idea to \eqref{dualL1L2-pos},
i.e., alternately   update the primal and dual variables, yields an
iterative algorithm with the  same updating formulae as
\eqref{alg-pd} except the computation for $z^{k+1}$ is replaced by
\begin{eqnarray}\label{update-z+}
z^{k+1}=\mP_{\mathcal{F}} (A^*y^{k} + x^k/\beta).
\end{eqnarray}
It is clear that the projection onto $\mathcal{F}$ is trivial.
The same procedure applies to the dual problems of other $\ell_1$-problems
with nonnegativity constraints as well.
Currently, with simple optional parameter settings, our Matlab package {\tt YALL1} can
be applied to solve \eqref{decoder-L1}, \eqref{decoder-BPDN}, \eqref{decoder-L1L2},
\eqref{decoder-L1L1} and their nonnegative counterparts.

\section{Choice of denoising models}

In most practical applications, measurements are invariably
contaminated by some noise. To date, the most widely used denoising
models in CS include \eqref{decoder-L1L2} and its variants, which
are  based on the $\ell_2$-norm fidelity. An alternative approach is
to penalize $Ax-b$ by the $\ell_1$-norm, resulting the
$\ell_1/\ell_1$ denoising model \eqref{decoder-L1L1}.
Between the two types of models, which type should be preferred in
general?  So far, model \eqref{decoder-L1L2}, along with its
variants, is widely used and appears to be the {\it de facto} model
of choice. However, we provide supporting evidence to argue that
perhaps in most practical situations the $\ell_1/\ell_1$ model
\eqref{decoder-L1L1} should be preferred unless there is an
excessive amount of white noise in observed data, in which case one
can only expect to obtain low-quality solutions. On the other hand,
if there is any possibility that observed data may contain large
measurement errors or impulsive noise, model \eqref{decoder-L1L1}
can potentially be dramatically better than \eqref{decoder-L1L2}
(see also, e.g., \cite{Wright-Ma08}).

We conducted  a set of experiments comparing models
\eqref{decoder-BPDN}, \eqref{decoder-L1L2} with \eqref{decoder-L1L1}
on random problems with $n=1000$, $m=300$ and $k=60$, using the
solver YALL1~\cite{YALL1}, which implements the dual ADMs described
in Subsection \ref{sc:DADM}  for eight different models, including
model \eqref{decoder-L1L1} that is reformulated and solved as a
basis pursuit problem in the form of \eqref{decoder-L1}. The
reformulation is as follows. Clearly, problem \eqref{decoder-L1L1}
is equivalent to
\begin{equation*}
\min_{x \in \C^n,r \in \C^m} \left\{\nu\|x\|_1  + \|r\|_1: Ax + r =
b\right\},
\end{equation*}
which, after a change of variables, can be rewritten as $\min_{\hat
x \in \C^{n+m}} \{ \|\hat x\|_1: \hat A \hat x = \hat b \}$, where
\begin{equation}
\hat{A} = \frac{\left(A, \nu I\right)} {\sqrt{1+\nu^2}}, \;\;\;
\hat{b}  = \frac{\nu }{\sqrt{1+\nu^2}}b \;\mbox{ ~and~ } \hat{x}  =
\left(\begin{array}{c}\nu x \\ r
\end{array}\right).
\end{equation}
We note that $\hat{A}\hat{A}^*=I$ provided that $AA^* = I$. In our
experiments, each model is solved for a sequence of parameter values
varying from 0 to 1. The simulation of data acquisition is as
follows
\begin{equation}\label{b-model}
b = A \bar x + p  \equiv A \bar x + p_W + p_I \equiv b_W + p_I,
\end{equation}
where $p, p_W, p_I$ represents, respectively, the total, the white and the
impulsive noise, and $b_W$ is the data containing white noise only
(or noiseless whenever $p_W=0$). White noise is generated by the
{\tt MATLAB} command {\tt randn(n,1)} multiplied by an appropriately
chosen constant to attain a desired signal-to-noise ratio, while
impulsive noise values are set to be $\pm 1$ at random positions of
$b$ that is always scaled so that $\|b\|_\infty=1$. As is usually
done, we measure the signal-to-noise ratio (SNR) of both $b_W$ and
$b$ in terms of decibel (dB). In this paper, we use the following
definition of SNR of $b$:
\begin{equation*}
\mbox{SNR}(b) = 20\log_{10}\left( \frac{\|b-\mbox{\bf
E}(b)\|}{\|p\|}\right),
\end{equation*}
where $\mbox{\bf E}(b)$ represents the mean value of $b$. The SNR of
$b_W$ is defined similarly. Impulsive noise is measured by
percentage, e.g., 5\% impulsive noise means 5\% of elements in $b$
are erroneous. For a computed solution $x$ from a data vector $b$
defined in (\ref{b-model}), the relative error in $x$ is defined in
terms of percentage:
\begin{eqnarray}\label{def-RelErr}
\mbox{RelErr}(x) = \frac{\|x- \bar x\|}{\|\bar x\|}\times 100\%.
\end{eqnarray}

\begin{figure}[htbp]
\vspace{-0cm}
\centering{
\includegraphics[scale=.38]{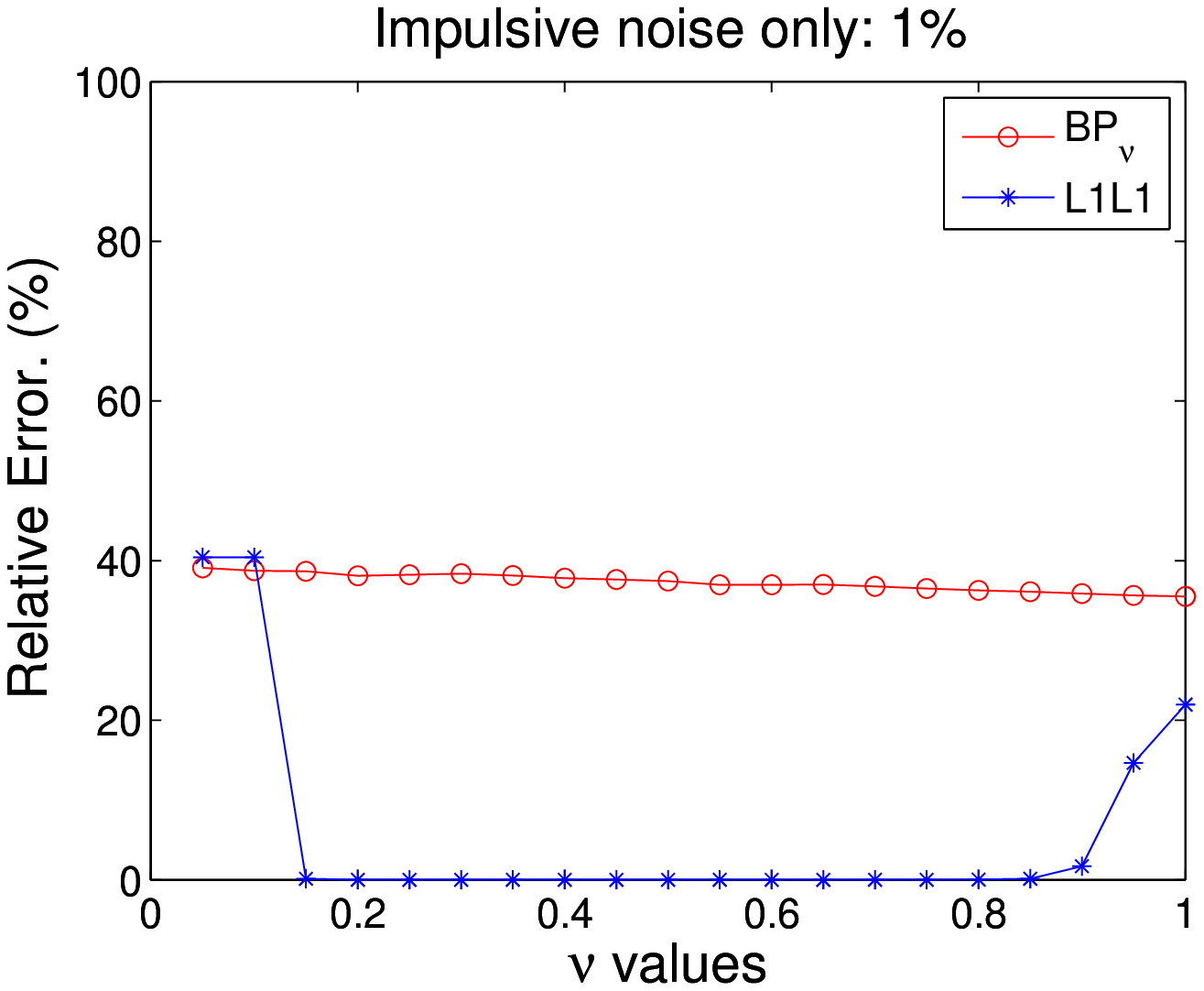}\hspace{-.5cm}
\includegraphics[scale=.38]{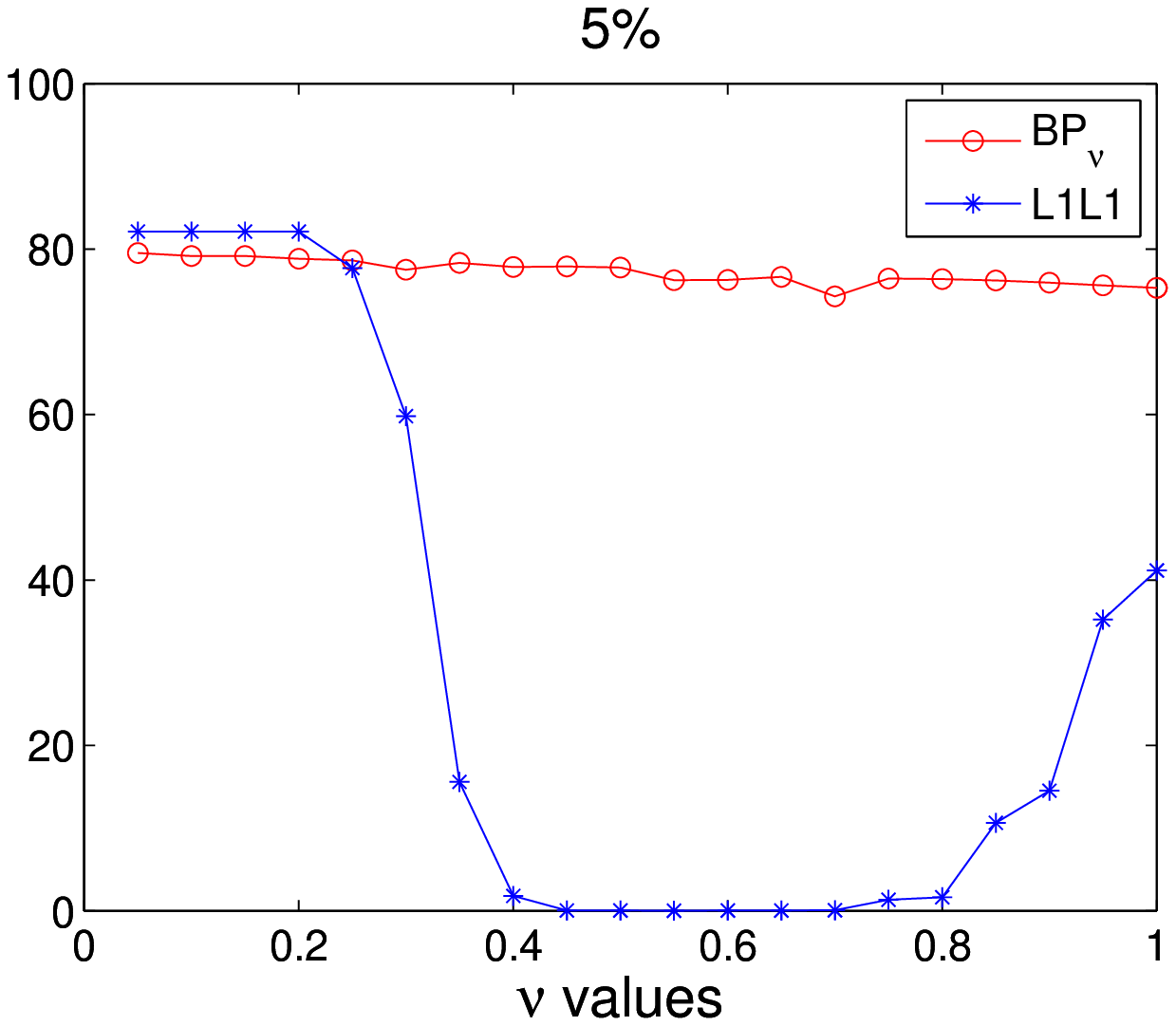}\hspace{-.5cm}
\includegraphics[scale=.38]{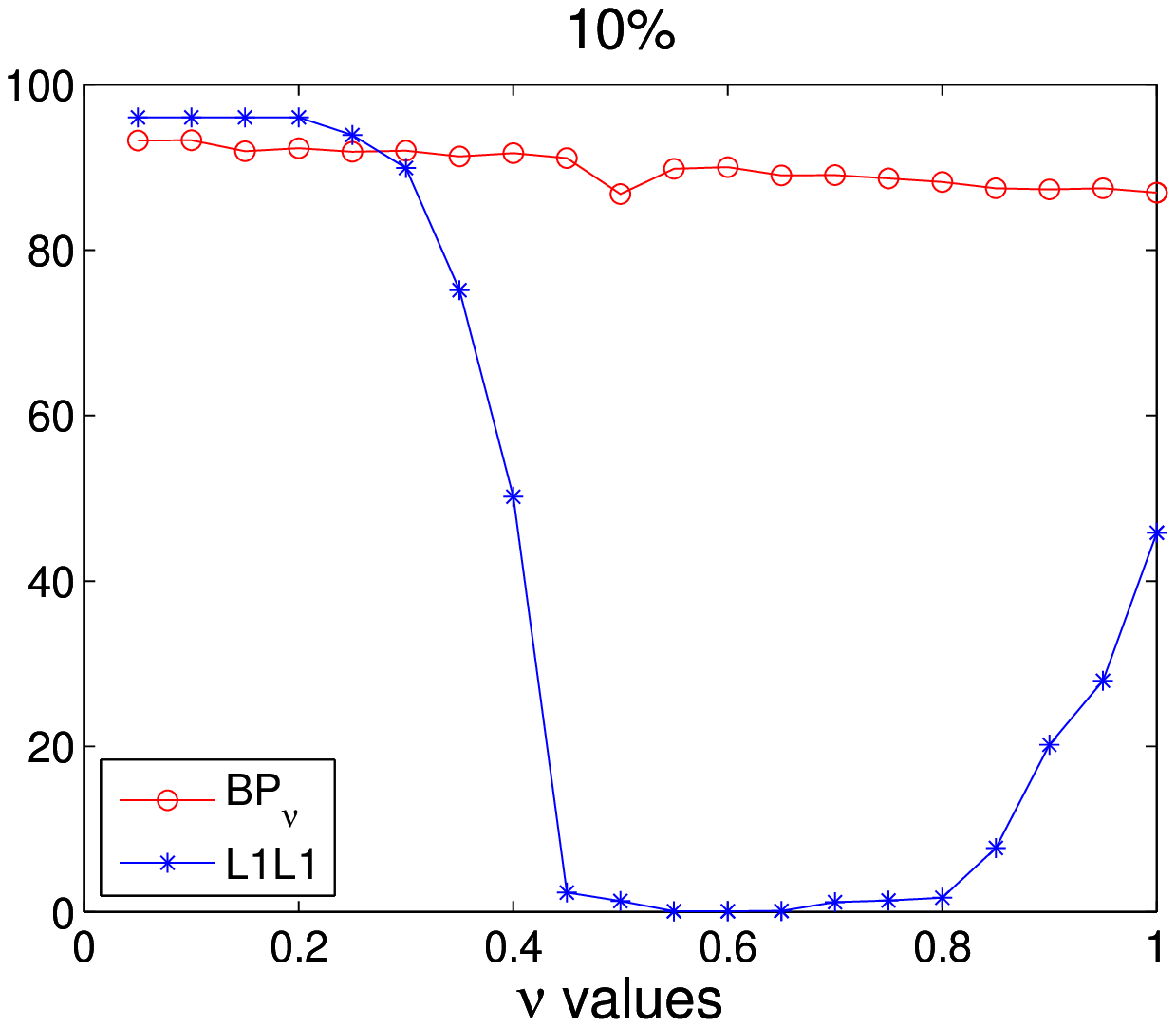}\\
\includegraphics[scale=.38]{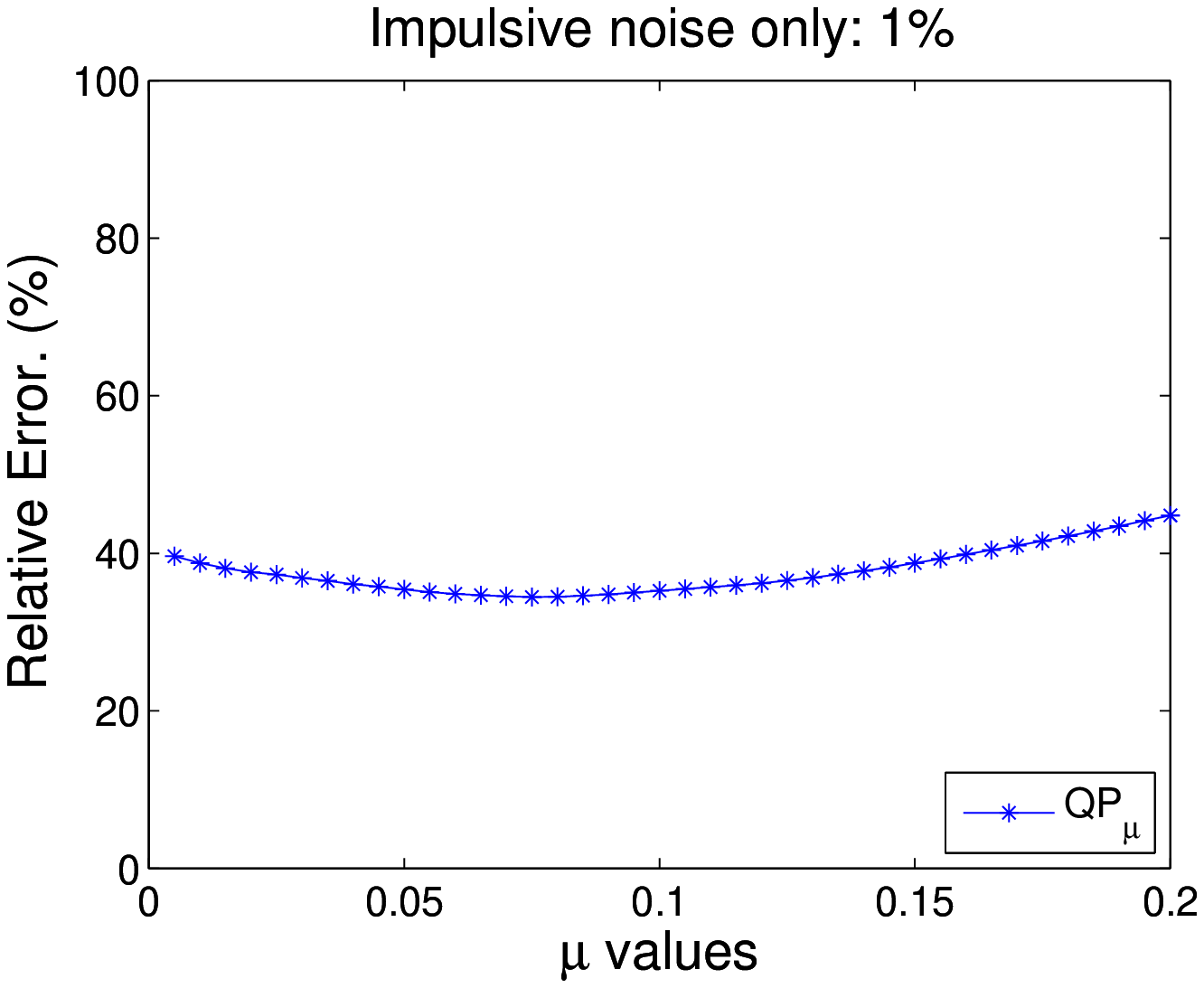}\hspace{-.5cm}
\includegraphics[scale=.38]{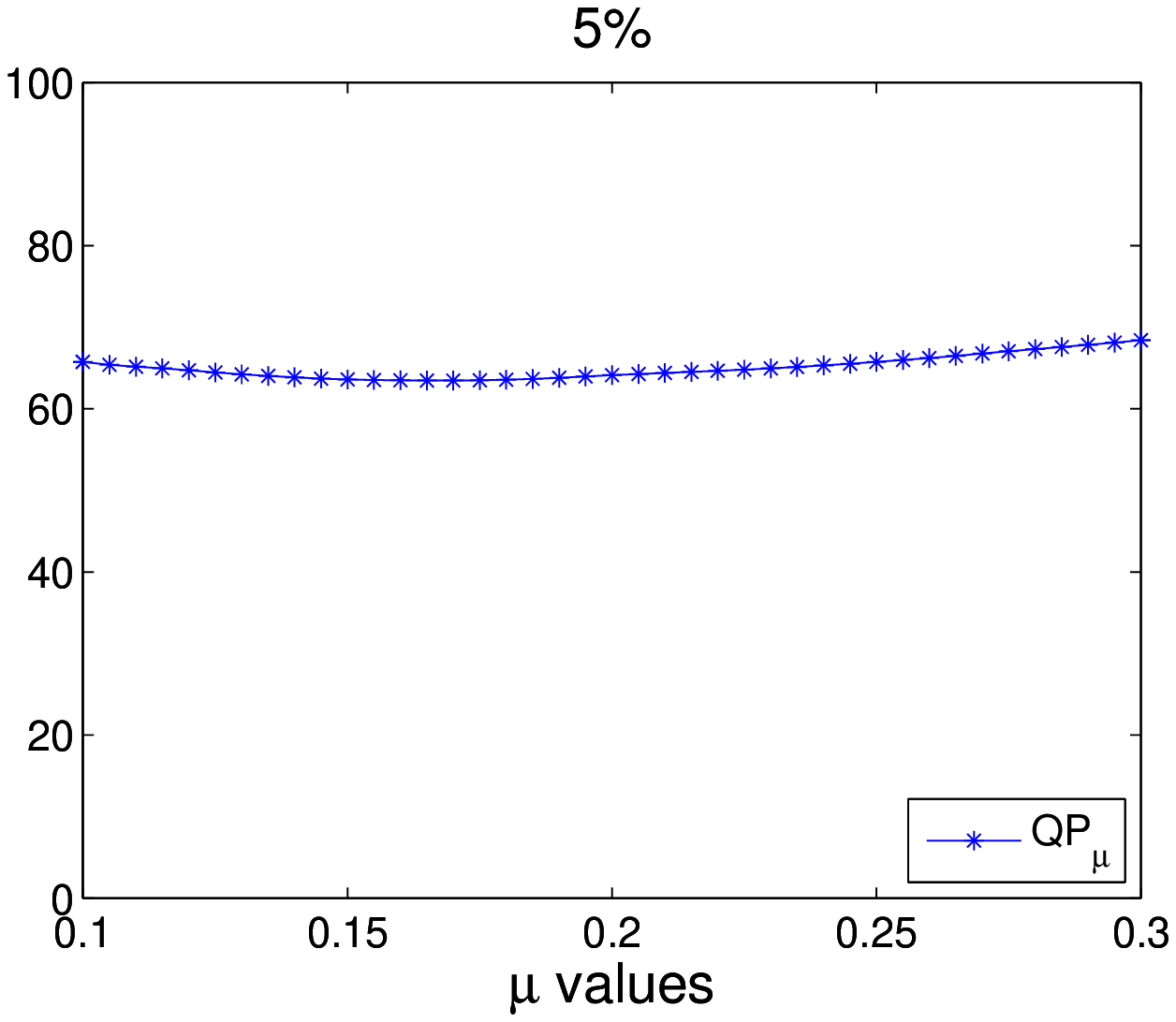}\hspace{-.5cm}
\includegraphics[scale=.38]{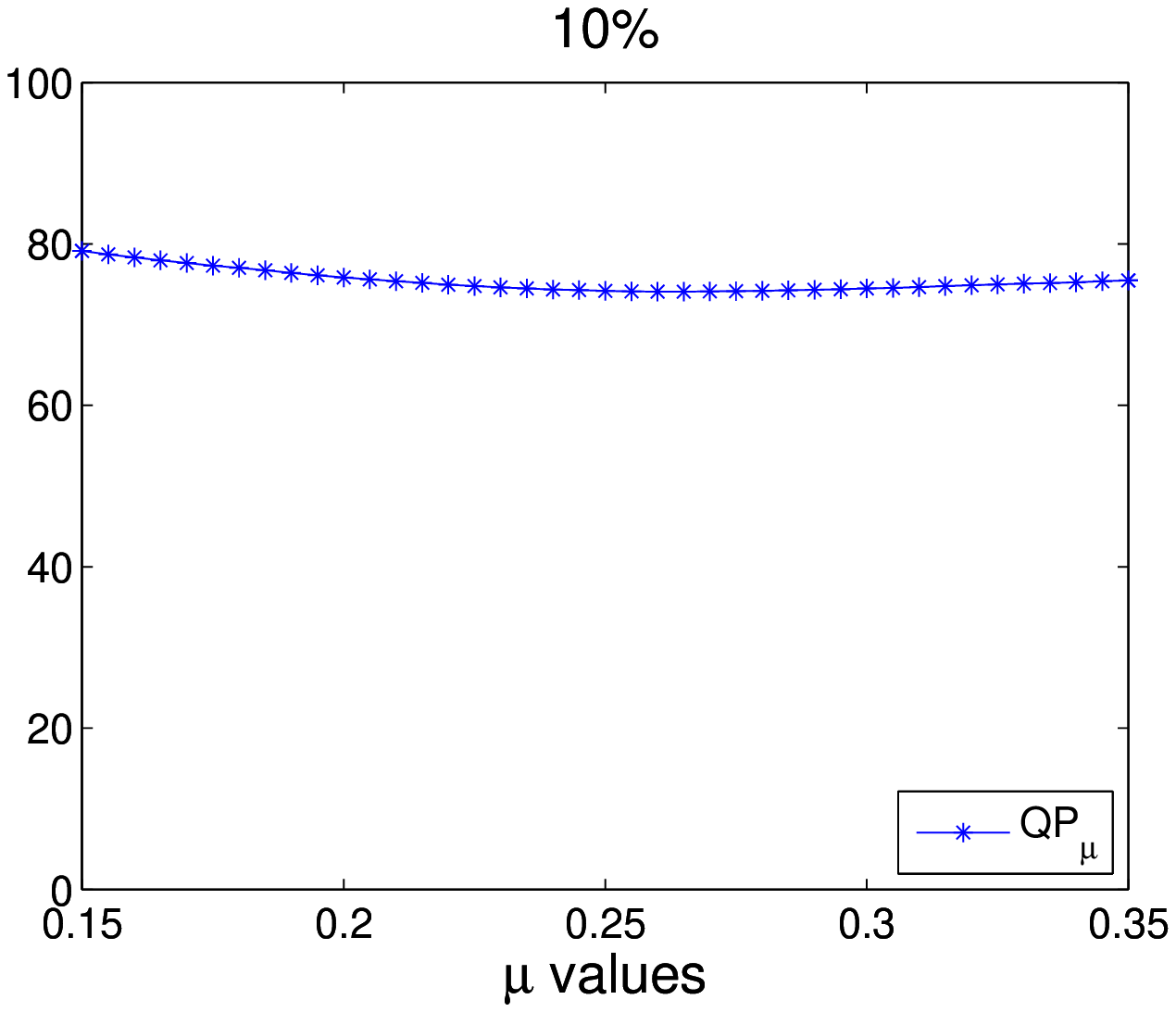}
}\caption{Recovered results from data corrupted by impulsive noise
(from left to right in both rows: 1\%, 5\% and 10\%). Top row:
recovered solutions from BP$_\nu$ ($\delta\leftarrow\nu$ in
\eqref{decoder-BPDN}) and $\ell_1/\ell_1$ problem
\eqref{decoder-L1L1} (the $x$-axes represent $\nu$ values in both
models); Bottom row: recovered solutions from QP$_\mu$ (the $x$-axes
represent $\mu$ values in the model). In all plots, the $y$-axes
represent relative errors of recovered solutions to the true sparse
signal. \label{fig-ip}}
\end{figure}

Figure~\ref{fig-ip} presents results for $p_W = 0$  and $p_I$
changes from 1\% to  10\%. From the results on the top row of Figure
\ref{fig-ip}, it is quite clear that model \eqref{decoder-L1L1} is
able to recover the exact solution $\bar x$ to a high accuracy for a
range of $\nu$ values (although the range for high quality recovery
shrinks when the corruption becomes increasingly severe), while model
\eqref{decoder-BPDN} is not, even though in all cases it reduces
relative errors by about 5\% when $\nu$ is close to 1 (we tried even
larger $\nu$ values but the achievable improvement soon saturates at
that level). From the results on the bottom row of Figure
\ref{fig-ip}, model \eqref{decoder-L1L2} behaves similarly as
\eqref{decoder-BPDN}, i.e., the quality of recovered signals is much
worse than those from \eqref{decoder-L1L1}.  Therefore, there should
be no doubt that model \eqref{decoder-L1L1} is superior to
\eqref{decoder-BPDN} and \eqref{decoder-L1L2} provided that $\nu$ is
chosen properly and should be the model of choice whenever there is
a possibility of erroneous measurements in data no matter how small
the percentage might be.

\begin{figure}[htbp]
\vspace{-0cm} \centering{
\includegraphics[scale=.38]{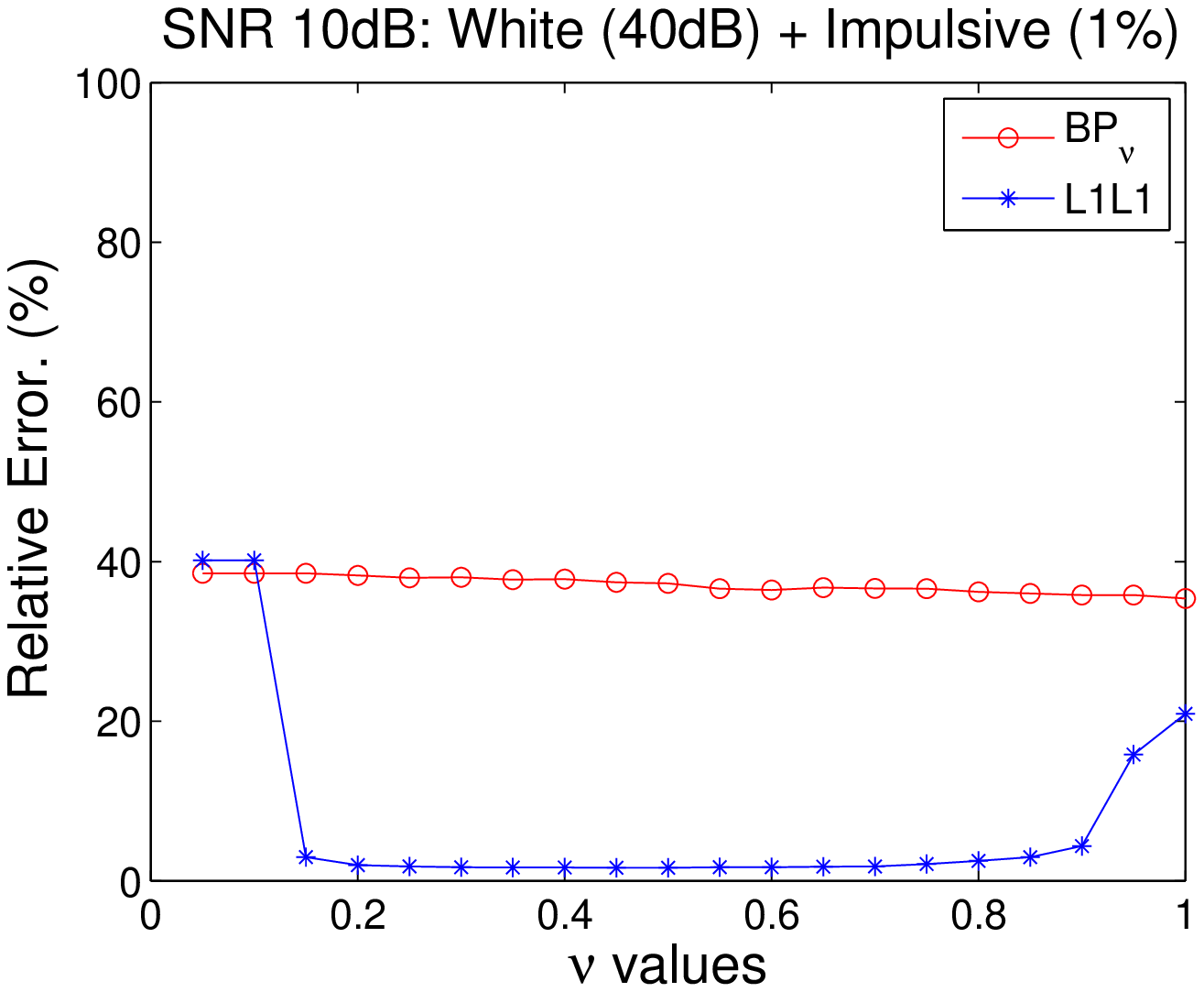}\hspace{-.5cm}
\includegraphics[scale=.38]{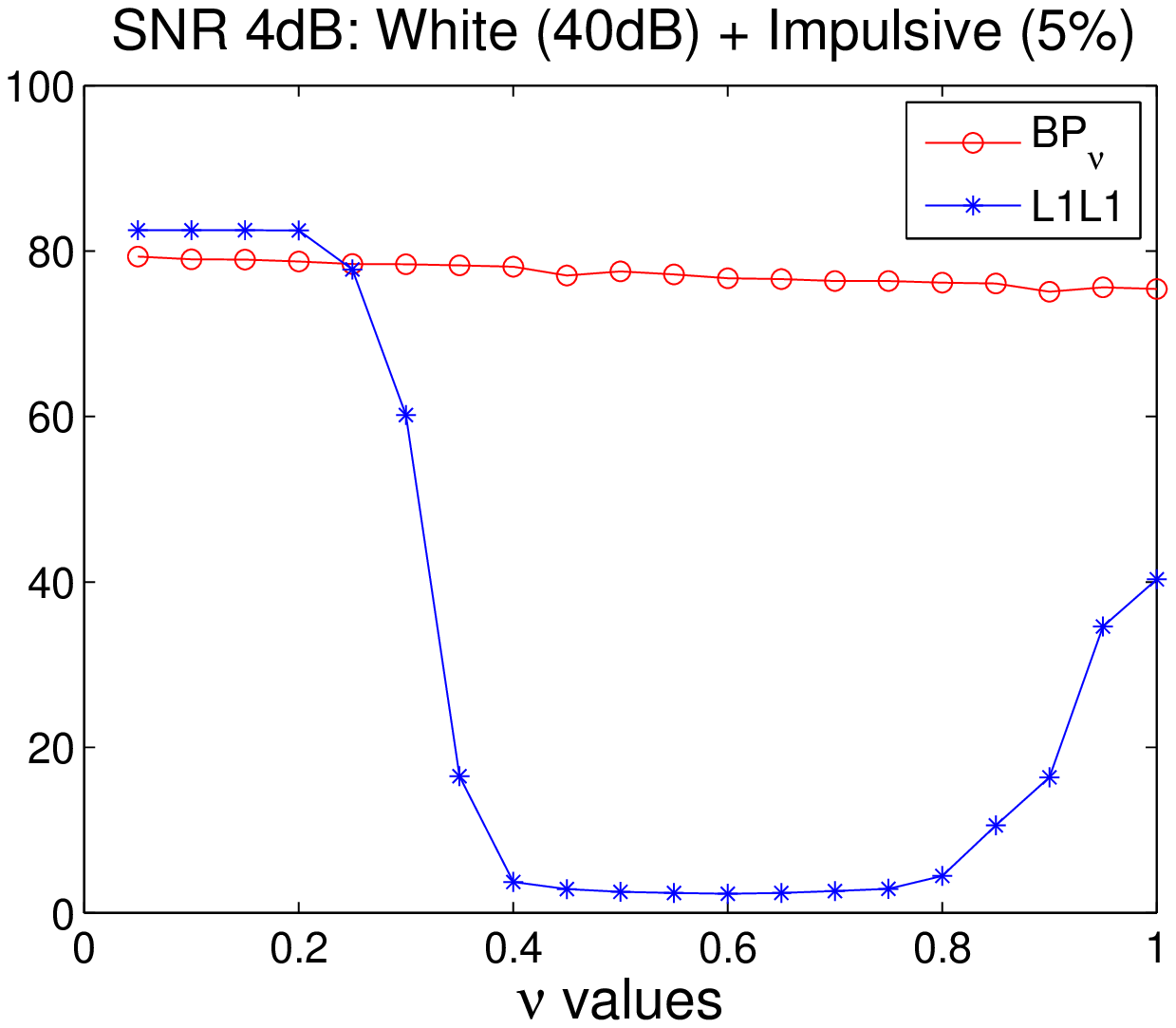}\hspace{-.5cm}
\includegraphics[scale=.38]{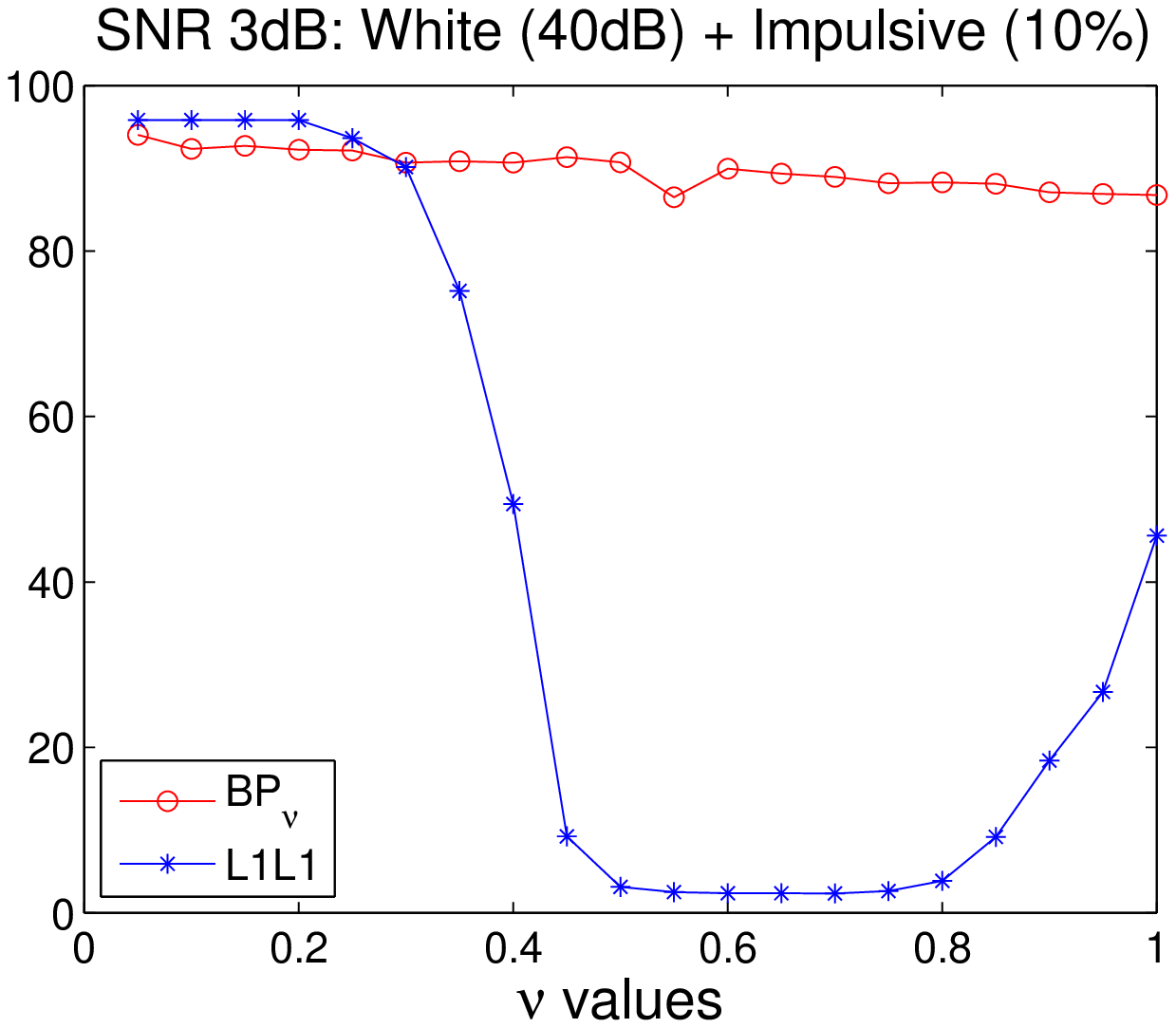}\\
\includegraphics[scale=.38]{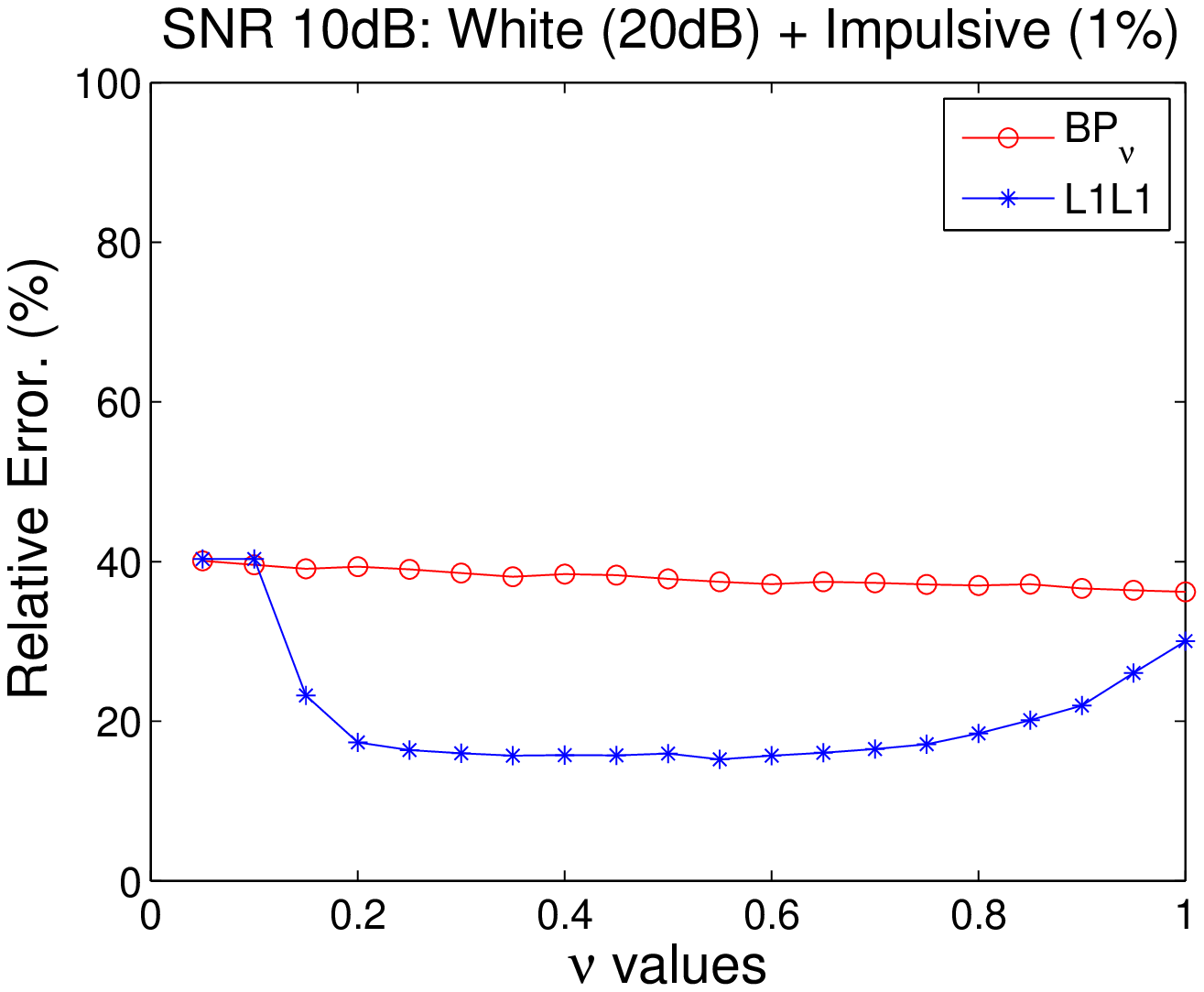}\hspace{-.5cm}
\includegraphics[scale=.38]{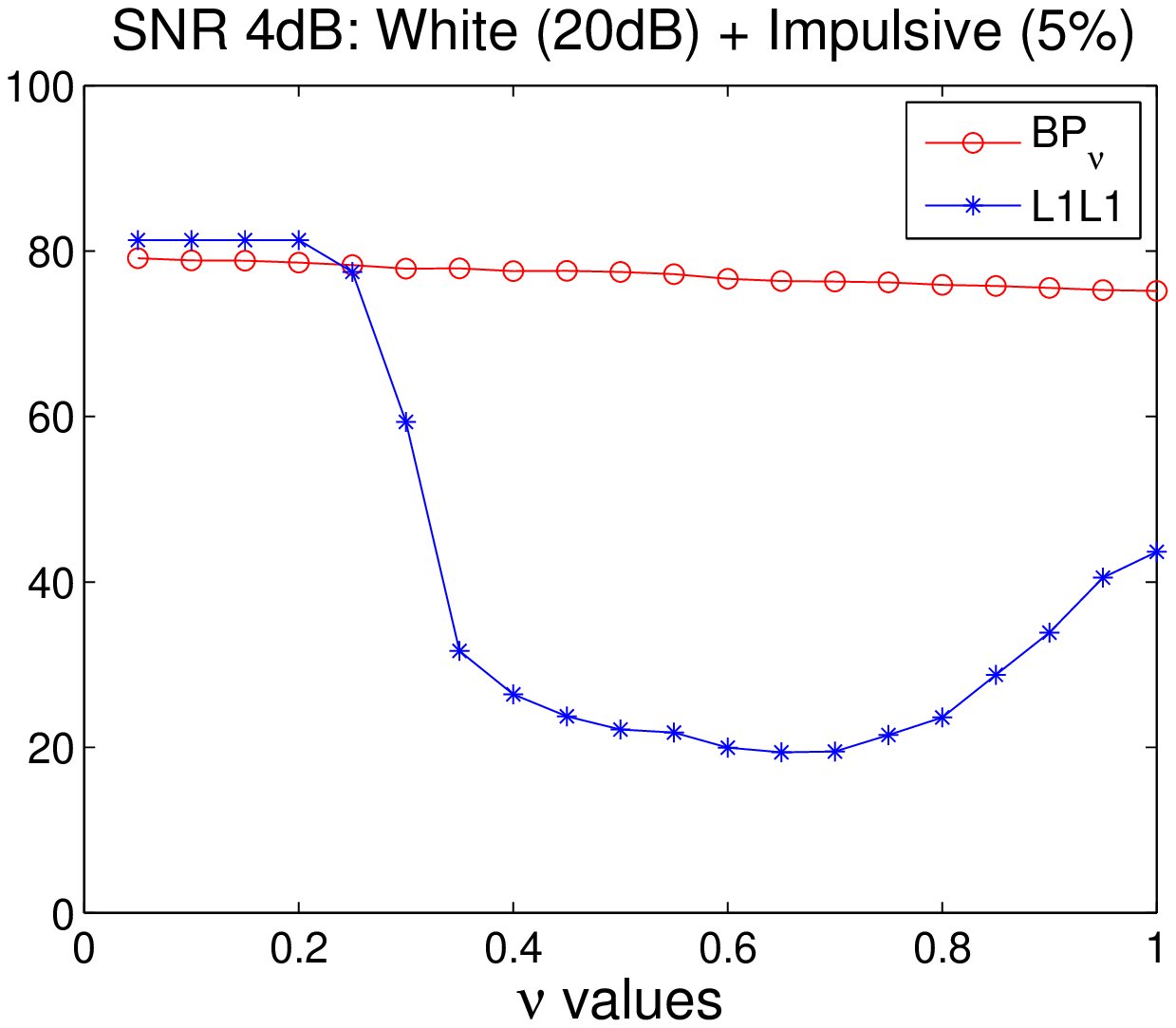}\hspace{-.5cm}
\includegraphics[scale=.38]{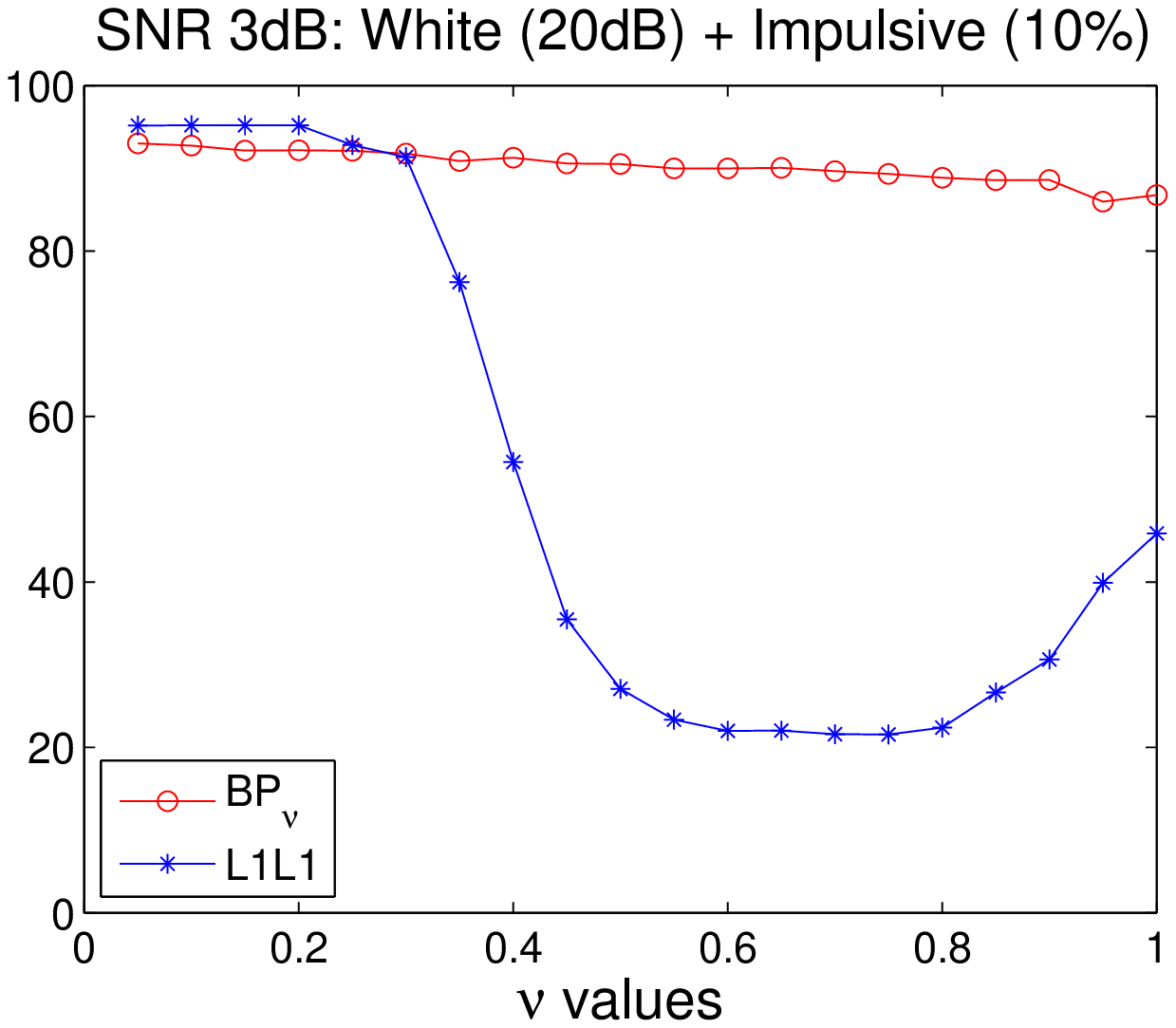}
}\caption{Recovered results from data corrupted by both white and
impulsive noise. SNR of $b_W$ is 40dB (top row) and 20dB (bottom
row); In each row, the ratios of impulsive noise corruption are 1\%,
5\% and 10\% from left to right; $x$-axes: model prameters in
BP$_\nu$ ($\delta\leftarrow\nu$ in \eqref{decoder-BPDN}) and
$\ell_1/\ell_1$ model \eqref{decoder-L1L1}; $y$-axes: relative
errors of recovered solutions to the true signal.
\label{fig-ip&white}}
\end{figure}

Which model should be used when data contain both white and
impulsive noise?  Let us examine results given in
Figure~\ref{fig-ip&white}, where \eqref{decoder-L1L1} is compared
with \eqref{decoder-BPDN} with data satisfying SNR$(b_W)$ = 40dB
(first row) and 20dB (second row). In both cases $p_I$ takes values
of 1\%, 5\% and 10\%. Similar to the noiseless case (free of white
noise), evidence strongly suggests that \eqref{decoder-L1L1} should
be the model of choice whenever there might be erroneous
measurements or impulsive noise in data even in the presence of
white noise. We did not present the results of \eqref{decoder-L1L2}
since they are similar to those of \eqref{decoder-BPDN}.

\begin{figure}[htbp]
\vspace{-0cm} \centering{
\includegraphics[scale=.38]{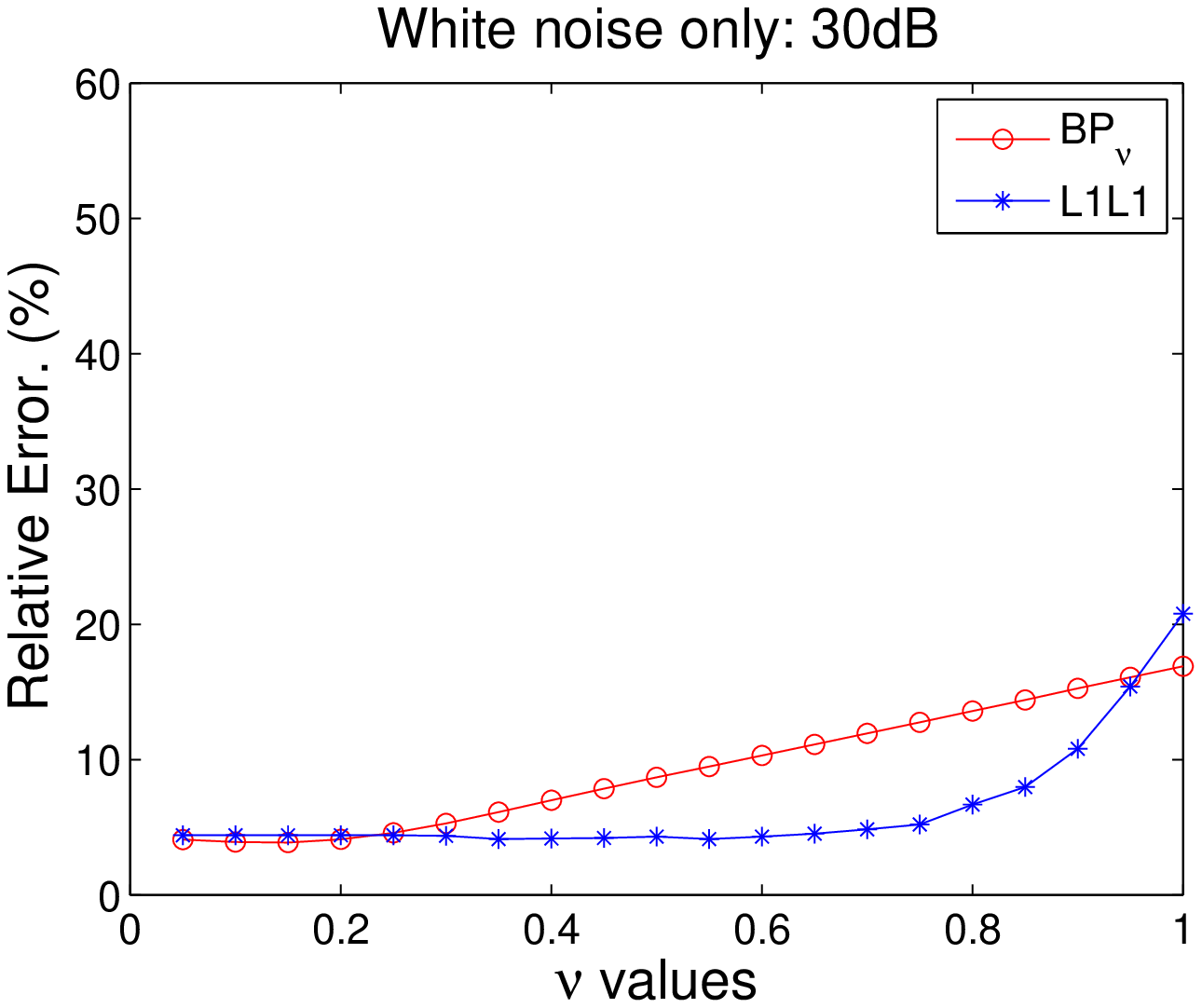}\hspace{-.4cm}
\includegraphics[scale=.38]{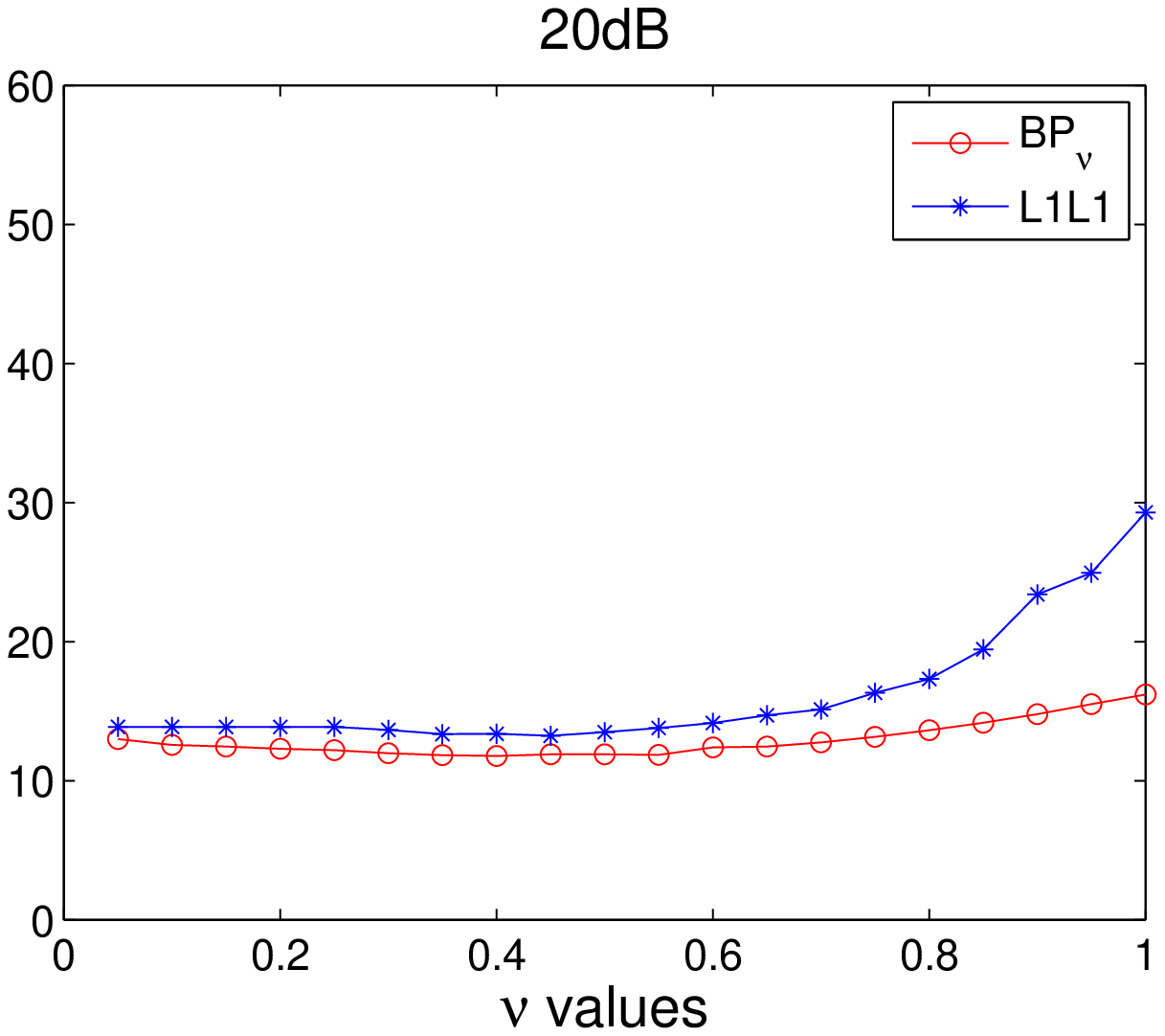}\hspace{-.4cm}
\includegraphics[scale=.38]{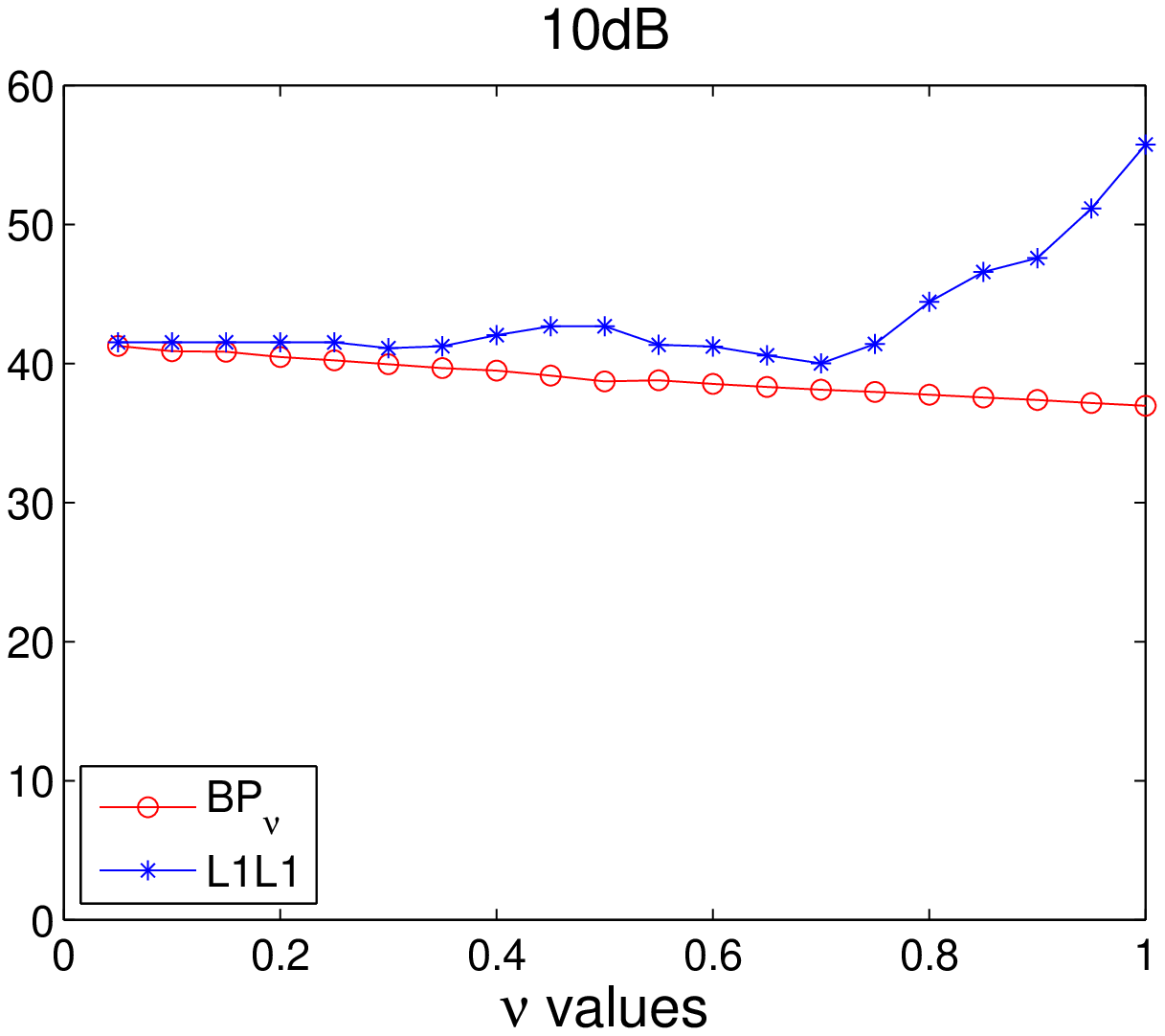}\\
\includegraphics[scale=.38]{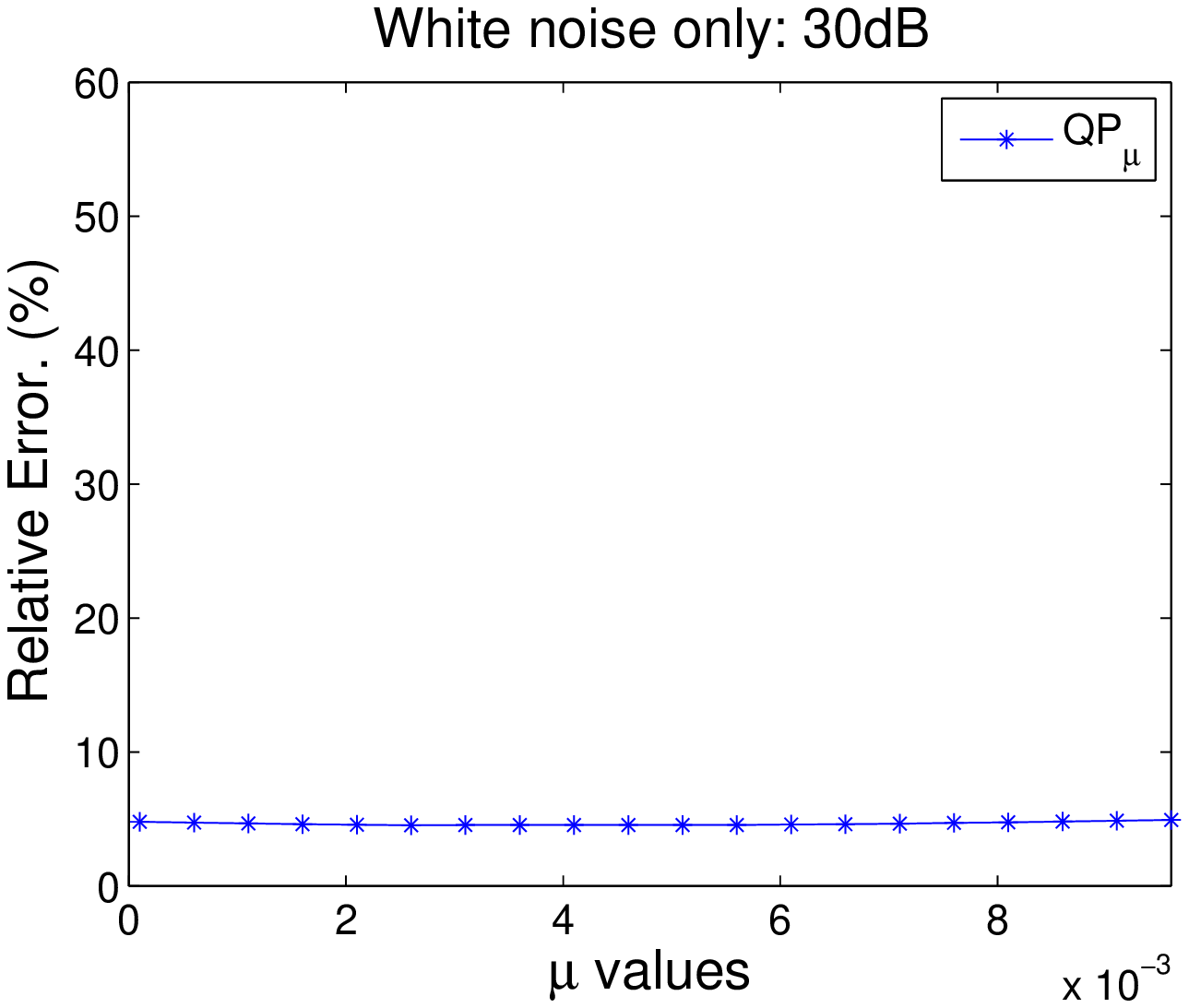}\hspace{-.4cm}
\includegraphics[scale=.38]{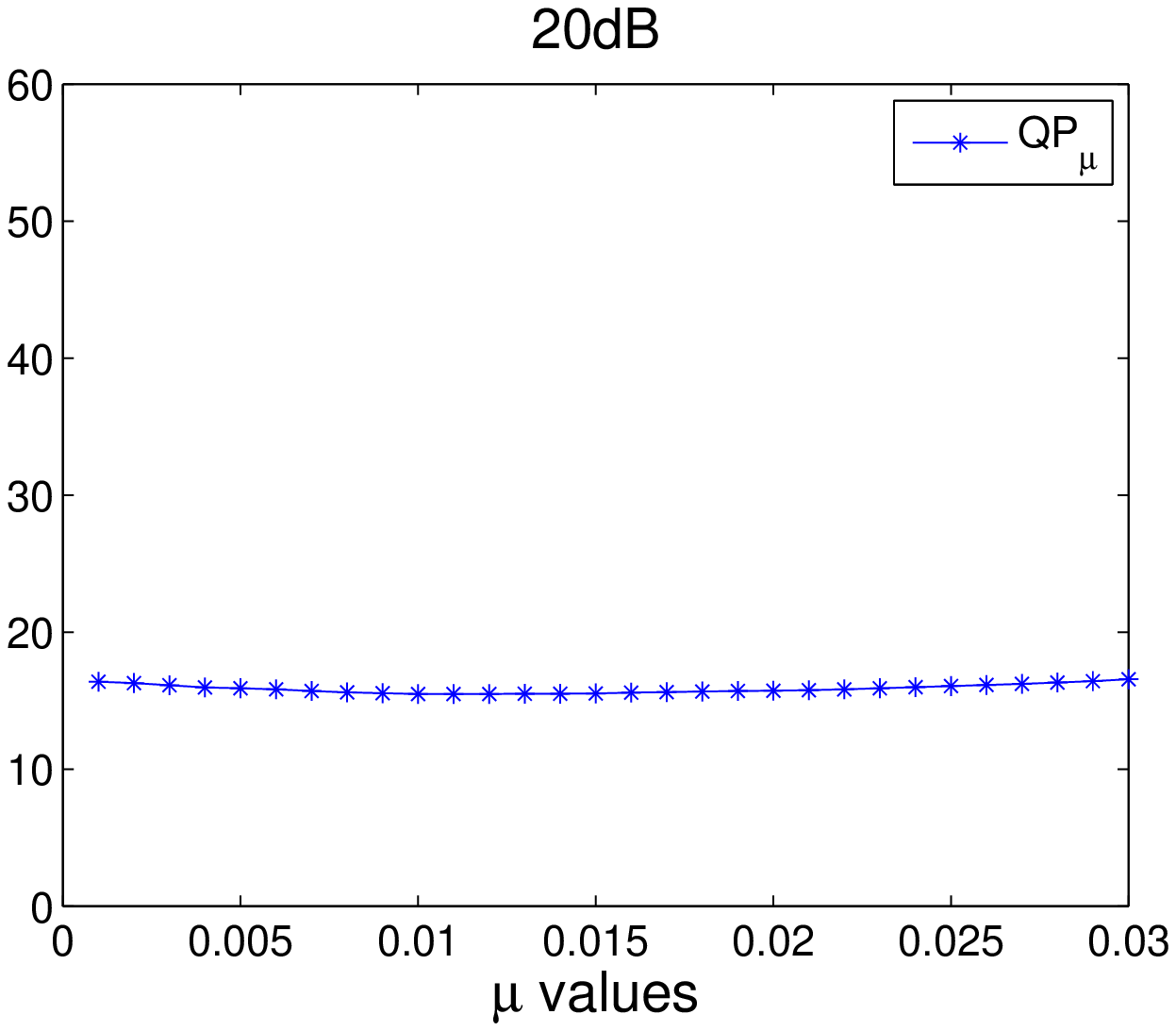}\hspace{-.4cm}
\includegraphics[scale=.38]{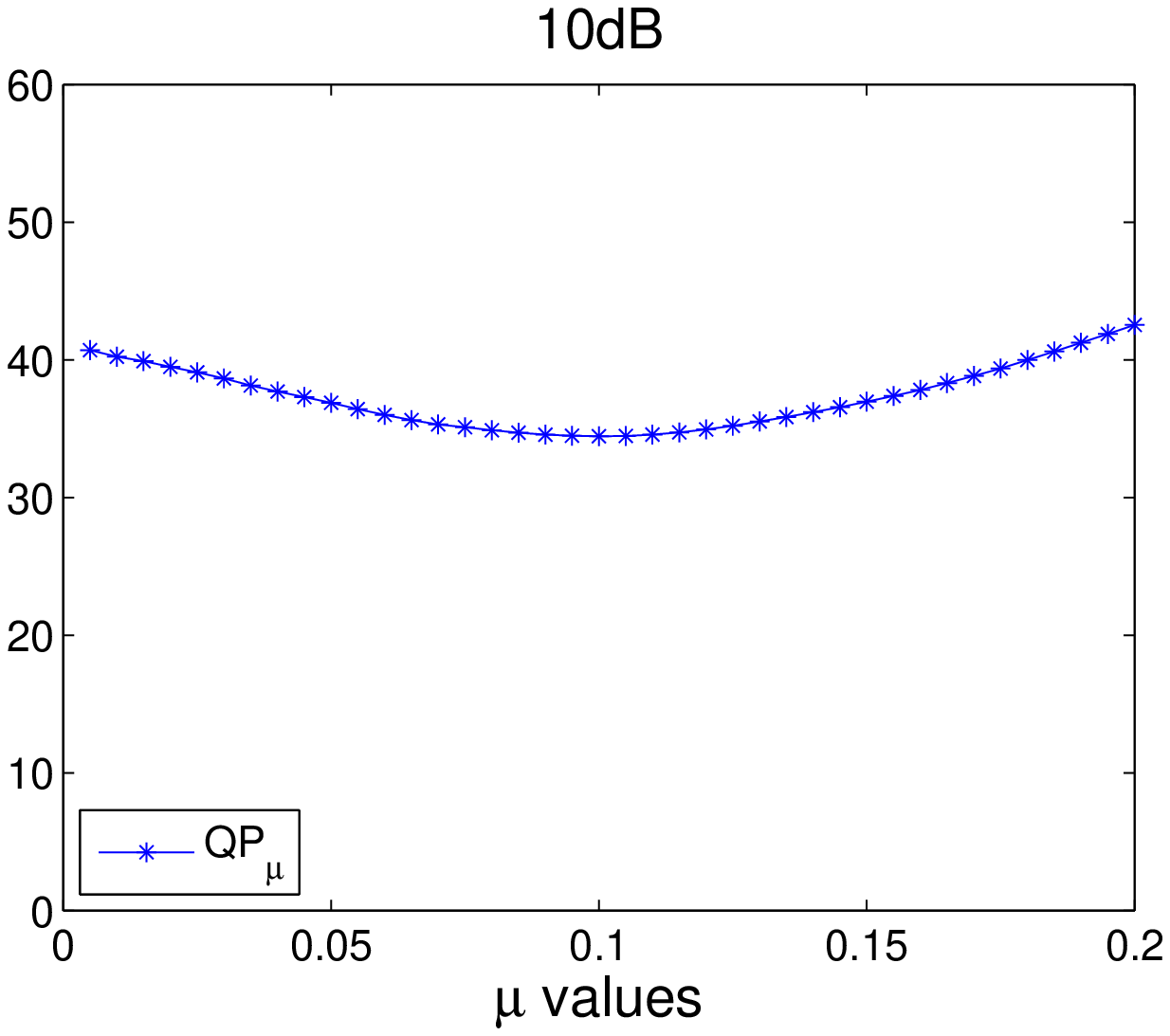}
} \caption{Recovered results from data contaminated by white noise
only. In both rows, SNR of $b_W$: 30dB, 20dB and 10dB from left to
right. Top row: results of BP$_\nu$ ($\delta\leftarrow\nu$ in
\eqref{decoder-BPDN})  and  $\ell_1/\ell_1$ model
\eqref{decoder-L1L1};  Bottom row: results of QP$_\mu$; In all
plots, the $x$-axes represent model parameters and the $y$-axes
represent relative errors of recovered solutions to the true signal.
\label{fig-white}}
\end{figure}

Is model \eqref{decoder-L1L1} still appropriate when there is no
impulsive noise?  Figure~\ref{fig-white} contains results obtained
from data with white noise only of SNR($b$) = SNR($b_W$) = 30dB,
20dB and 10dB, respectively.
Qualitatively and loosely speaking, these three types
of data can be characterized as good, fair and poor, respectively.
As can be seen from the top left plot, on good data
\eqref{decoder-BPDN} offers no improvement whatsoever to the basis
pursuit model ($\nu=0$) as $\nu$ decreases.  On the contrary, it
starts to degrade the quality of solution once $\nu > 0.25$.  On the
other hand, model \eqref{decoder-L1L1} essentially does no harm
until $\nu > 0.7$. From the top middle plot, it can be seen that on
fair data both models start to degrade the quality of solution after
$\nu > 0.7$, while the rate of degradation is faster for model
\eqref{decoder-L1L1}. Only in the case of poor data (the top right
plot), model \eqref{decoder-BPDN} always offers better solution
quality than that of model \eqref{decoder-L1L1}. However, for poor
data the solution quality is always poor for $\nu \in [0,1]$ (and
beyond as well).  At $\nu=1$, the relative error of model
\eqref{decoder-BPDN} is about 38\%, representing a less than 5\%
improvement over the relative error 42\% at $\nu=0.05$, while the
best error attained from model \eqref{decoder-L1L1} is about 40\%.
The results of \eqref{decoder-L1L2} are generally similar to those
of \eqref{decoder-BPDN} provided that model parameters are selected
properly, as is shown   on the bottom row of Figure \ref{fig-white}.

The sum of the computational evidence suggests the following three
guidelines, at least for random problems of the type tested here:
\begin{enumerate}

\item
Whenever data may contain large measurement errors, erroneous
observations or in general impulsive noise, the $\ell_1$-norm
fidelity used by model  \eqref{decoder-L1L1} should naturally be preferred
over the $\ell_2$-norm fidelity used by model \eqref{decoder-BPDN}
and its variants.
\item
Without impulsive noise, the $\ell_1$-norm  fidelity does no  harm
to solution quality, as long as data do not contain a large amount
of white noise and $\nu$ remains reasonably small.
\item
When data are contaminated by a large amount of white noise, the
$\ell_2$-norm  fidelity should be preferred. In this case, however,
high-quality recovery should not be expected regardless what model
is used.

\end{enumerate}

\section{Numerical results} \label{sc:numerical}

In this section, we compare the proposed ADMs for $\ell_1$ problems,
which will be referred as PADM and DADM, respectively, with several
state-of-the-art algorithms to demonstrate their practical
performance.  The rest of this section is organized as follows.  In
Section \ref{sc:Rerr-vs-opt}, we present experimental results to
emphasize a simple yet often overlooked point that algorithm speed
should be evaluated relative to solution accuracy. In Section
\ref{sc: alg-setting}, we describe experimental settings, including
parameter choices, stopping rules and generation of problem data
under {\tt MATLAB} environment. In Section \ref{sc:compQPmu}, we
compare PADM and DADM with FPC-BB \cite{FPC, FPC2}, a fixed-point
continuation method with non-monotone line search based the Barzilai
and Borwein (BB) steplength \cite{BB88}, SpaRSA \cite{SpaRSA} --- a
sparse reconstruction algorithm for more general regularizers, FISTA
\cite{FISTA} --- a fast IST algorithm that attains an optimal
convergence rate in function values, and CGD \cite{CGD} --- a block
coordinate gradient descent method for minimizing
$\ell_1$-regularized convex smooth function. In Section
\ref{sc:compBPdelta} we compare PADM and DADM with SPGL1
\cite{SPGL1} --- a spectral projected gradient algorithm for
\eqref{decoder-BPDN}, and NESTA \cite{NESTA}
--- a fast and accurate fast first-order algorithm based on
Nesterov's smoothing technique \cite{Nesterov-smoothing}. All
experiments were performed under Windows Vista Premium and {\tt MATLAB}
v7.8 (R2009a) running on a Lenovo laptop with an Intel Core 2 Duo
CPU at 1.8 GHz and 2 GB of memory.

\subsection{Relative error versus optimality}\label{sc:Rerr-vs-opt}

In algorithm assessment, the speed of an algorithm is often taken as
an important criterion, and rightfully so.   However, speed is a
relative concept and ought to be used within context, which
unfortunately has not always been the case. In track and field,
human speed is measured relative to the distance of a race.  The
fastest short distance runner may very well be a rather slow long
distance runner and vice versa.  Similarly, algorithm speed should
be measured relative to accuracy.  Using a single accuracy to
judge algorithm speed could be as misleading as comparing the
speeds of short and long distance runners.

\begin{figure}[htbp]
\centering{
\includegraphics[scale=.38]{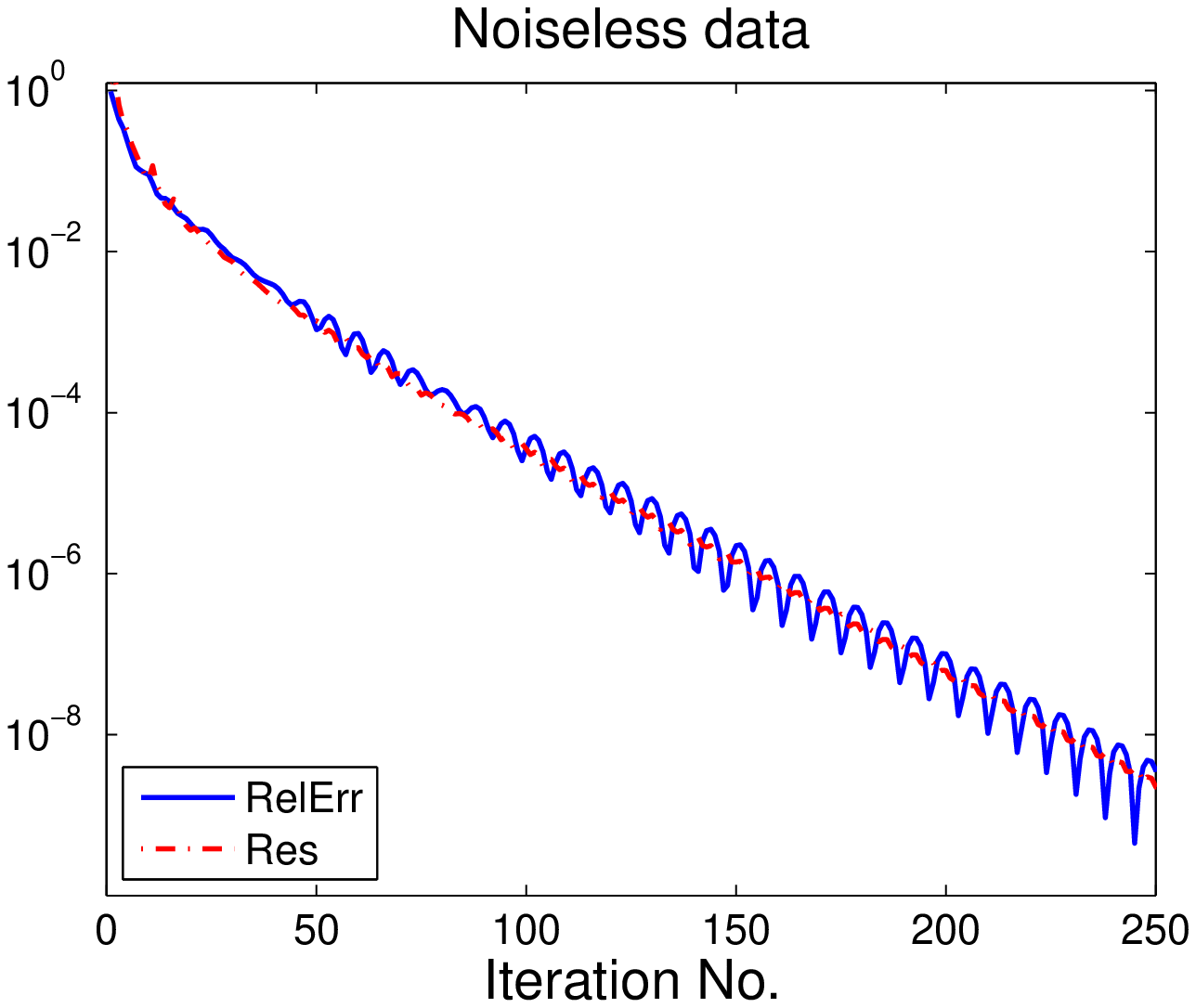}\hspace{-.5cm}
\includegraphics[scale=.38]{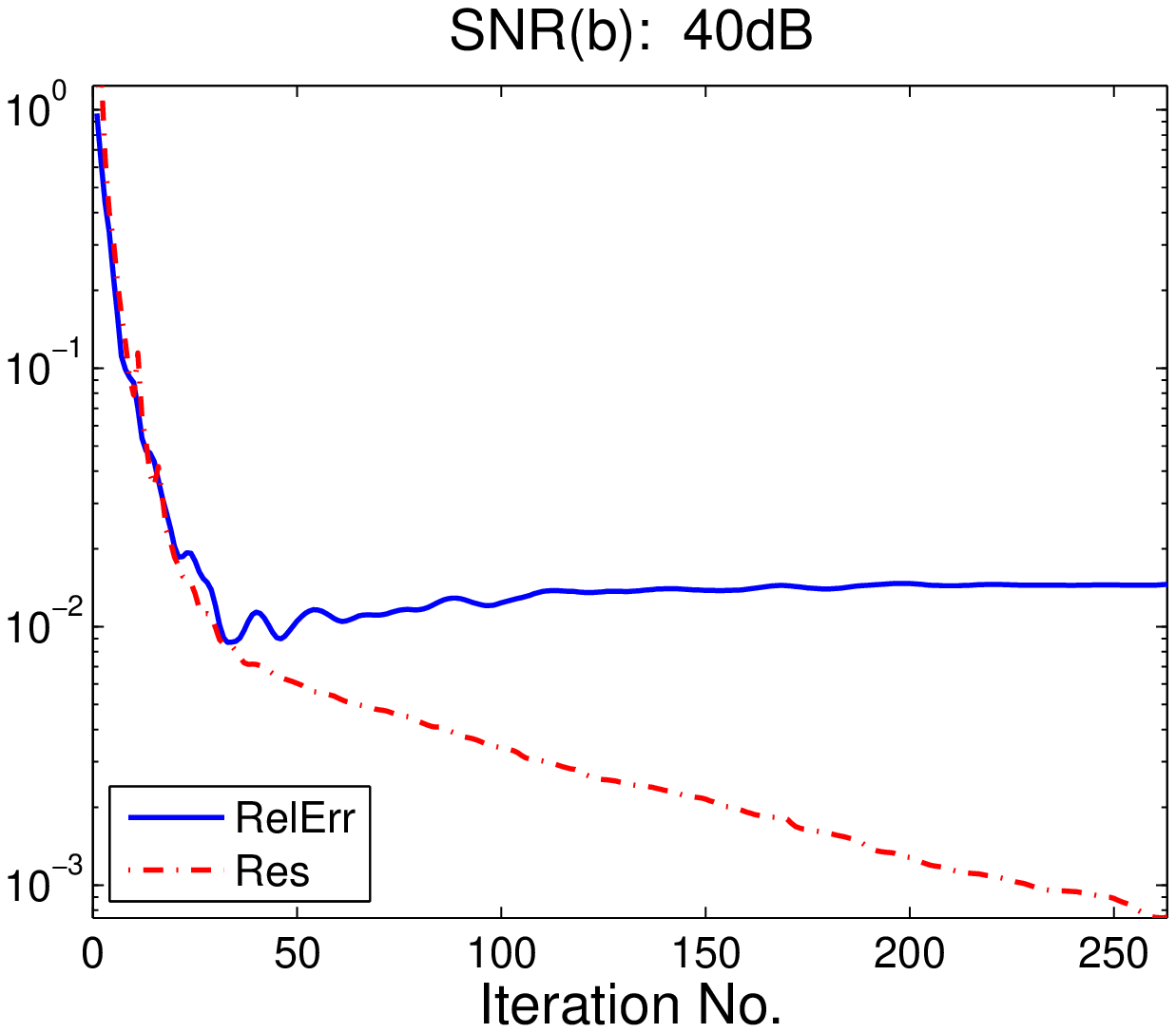}\hspace{-.5cm}
\includegraphics[scale=.38]{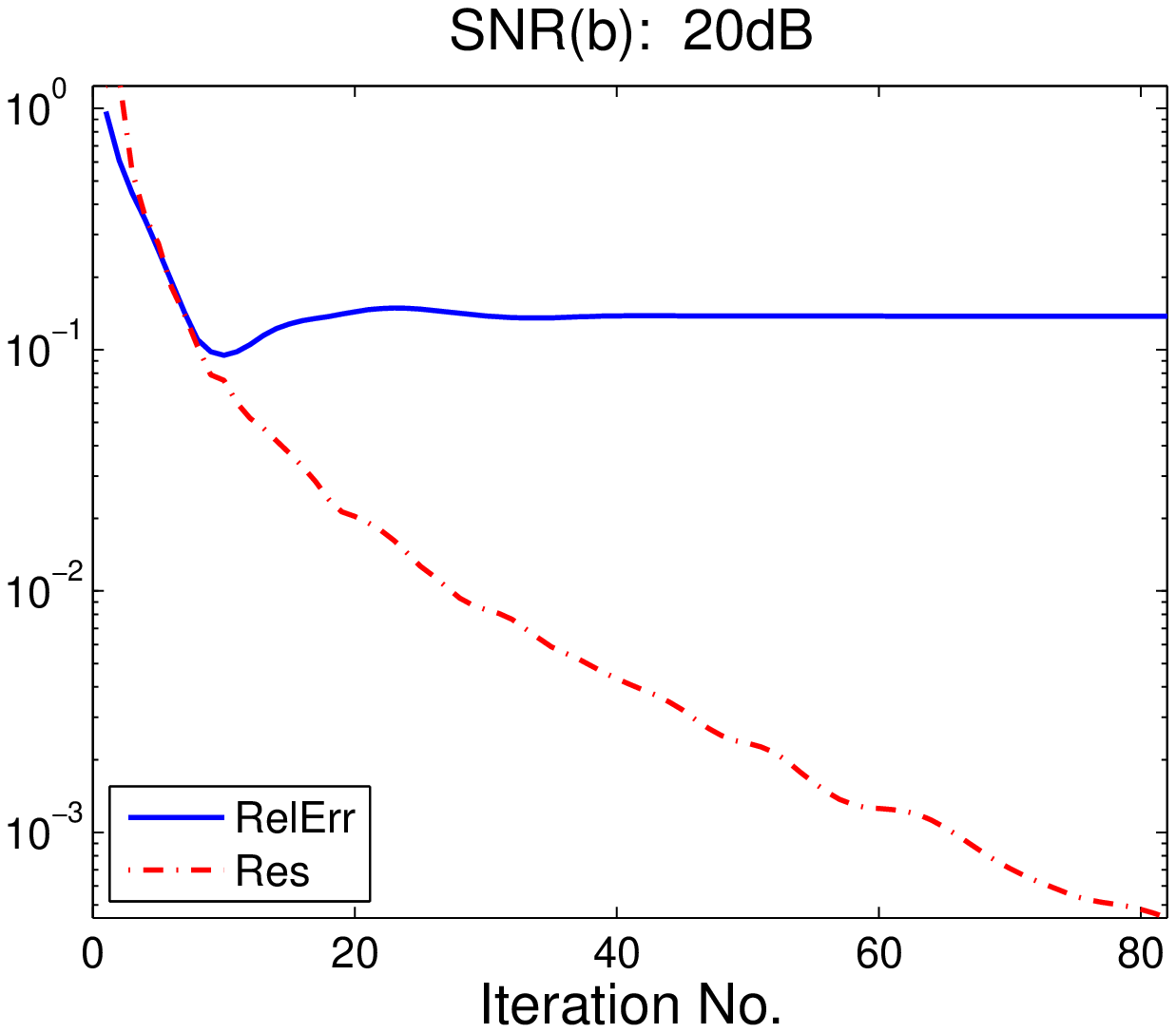}
} \caption{Relative error  versus optimality   for noiseless and
noisy data. The $x$-axes represent number of iterations, and $y$-axes
represent relative error (RelErr) defined in \eqref{def-RelErr} and
residue (Res) defined in \eqref{def-OptCond}.
\label{fig-Err-vs-optim-noiseless}}
\end{figure}

Clearly,  the definition for a good accuracy can vary with situation
and application.  Consequently, algorithm speed should be judged
within concrete context.  A relevant question here is what kind of
accuracy is reasonable for solving compressive sensing problems,
especially when data are noisy as is the case in most real applications.
To answer this question, we solved \eqref{decoder-L1} with noiseless
data and solved \eqref{decoder-L1L2} with data contaminated by white
noise of small to moderate levels.  In this experiment, the measurement
matrix was constructed by orthogonalizing and normalizing the rows of
a 330 by 1000 standard random Gaussian matrix.
The true signal $\bar x$ has 60 nonzeros whose positions are
determined at random, and the nonzero values are random Gaussian.
Both problems \eqref{decoder-L1} and \eqref{decoder-L1L2} were
solved by DADM to a relative high accuracy. The results on relative
error defined in \eqref{def-RelErr} and optimality measured in
residue defined in \eqref{def-OptCond} are given in Figure
\ref{fig-Err-vs-optim-noiseless}.

It is clear from Figure \ref{fig-Err-vs-optim-noiseless} that
solving $\ell_1$-problems to high accuracy improves solution quality
only when the observed data are free of noise. In the left plot of
Figure \ref{fig-Err-vs-optim-noiseless} where noiseless data were
used in \eqref{decoder-L1}, both relative error and optimality
measured by residue decrease as DADM proceeds. For noisy data, a
relatively low accuracy is sufficient to give the best relative
error that an $\ell_1$-denoising model can reach, e.g., in the
middle plot of Figure \ref{fig-Err-vs-optim-noiseless} where low
level noisy data were used in \eqref{decoder-L1L2}, relative error
does not decrease anymore after the residue is reduced to about
$10^{-2}$ in about 40 iterations. This phenomenon becomes more and
more obvious for higher level noisy data, as is shown in the right
plot of Figure \ref{fig-Err-vs-optim-noiseless}. Based on these
observations, we conclude that solving $\ell_1$-problems to high
accuracy is not necessary whenever observed data are contaminated by
noise. We choose to emphasize this rather mundane point because too
often such a common sense has been ignored in algorithmic studies.
In our numerical comparison, whenever noisy data are used we will
not compare how fast algorithms achieve a {\em high} accuracy, but
how fast they achieve an {\em appropriate} accuracy that is
consistent with the noise level of data.

\subsection{Experimental settings}\label{sc: alg-setting}
Now we describe experimental settings, including selection of
parameters, stopping rules  and  generation of data in {\tt MATLAB}. In
our experiments, we mainly tested randomized partial transform
matrices as encoding matrices for reasons given below. First,
partial transform matrices do not need explicit storage and thus are
suitable for large scale tests. Second, fast matrix-vector
multiplications are available, which are the main computation load
of all algorithms to be compared with. Finally, for this class of
matrices there hold $AA^*=I$, which is not only useful in parameter
choices in PADM but also required in DADM for exact minimization of
subproblems.  We note that condition $AA^*=I$ implies
$\lambda_{\max}(A^* A)=1$. Therefore, for convergence of PADM, in
our experiments we set $\tau=0.8$, $\gamma=1.199$ and
$\beta=2m/\|b\|_1$ in \eqref{alg-primal} and \eqref{alg-primal-Con},
which work quite well in practice, although suitably larger $\tau$
and $\gamma$ seem work better most of the time. For DADM, we used
the default settings in YALL1, i.e., $\gamma=1.618$ and
$\beta=\|b\|_1/m$. As described in subsection \ref{sc:Rerr-vs-opt},
high accuracy is not always necessary in CS problems especially for
noisy data. Thus, when comparing with other algorithms, we simply
terminated PADM and DADM when relative change of two consecutive
iterates becomes small, i.e.,
\begin{eqnarray}\label{def:RelChg}
\text{RelChg}\triangleq\frac{\|x^{k+1}-x^{k}\|}{\|x^{k}\|} <
\epsilon,
\end{eqnarray}
where $\epsilon>0$ is  a tolerance, although more complicated
stopping rules, such as the one based on optimality conditions
defined in \eqref{def-OptCond},  are possible. Parametric  settings
of FPC-BB, SpaRSA, FISTA, CGD, SPGL1 and NESTA will be specified
when we make comparison with them. In all experiments, we generated
data $b$ by {\tt MATLAB} scripts {\tt b = A*xbar +
sigma*randn(m,1)}, where {\tt A} is a randomized partial
Walsh-Hadamard transform matrix, {\tt xbar} represents a sparse
signal that we wish to recover, and {\tt sigma} is the standard
deviation of additive  Gaussian  noise. The partial Walsh-Hadamard
transform is implemented in the C language  with a {\tt MATLAB}
mex-interface available to all codes compared. In all tests, we set
$n = 8192$ and tested various combinations of $m$ and $k$ (the
number of nonzero components in {\tt xbar}). We initialized $x$ to
the zero vector in all algorithms unless otherwise specified.

\subsection{Comparison with FPC-BB, SpaRSA, FISTA and CGD}\label{sc:compQPmu}
In this subsection, we present comparison results of PADM and DADM
with FPC-BB~\cite{FPC,FPC2}, SpaRSA~\cite{SpaRSA},
FISTA~\cite{FISTA} and CGD~\cite{CGD}, all of which were developed
in the last two years for solving  \eqref{decoder-L1L2}. In the
comparison, we used random Gaussian spikes as  {\tt xbar}, i.e., the
location of nonzeros are selected uniformly at random while the
values of nonzero components are $iid$   Gaussian entries. The
standard deviation  of additive noise is set to be $10^{-3}$, while
model parameter $\mu$ in \eqref{decoder-L1L2} is set to be either
$10^{-3}$ or $10^{-4}$, both of which are suitable values for small
noise.  We aimed to compare both the efficiency and robustness of
the algorithms.

Since different algorithms use different stopping criteria, it is
rather difficult to compare their relative performance completely
fairly. Therefore, we present two classes of comparison results. In
the first class of results, we run PADM and DADM for a prescribed
number of iterations, and  choose proper stopping tolerance values
for other algorithms so that their iteration numbers are roughly
euqal to the prescribed iteration number. Then we examine how
relative errors and function values decrease as each algorithm
proceeds. In the second class of results, we terminate the ADM
algorithms by \eqref{def:RelChg}, while the stopping rules used for
other algorithms in comparison will be specified below.

\begin{figure}[htbp]
\centering{
\includegraphics[scale=.5]{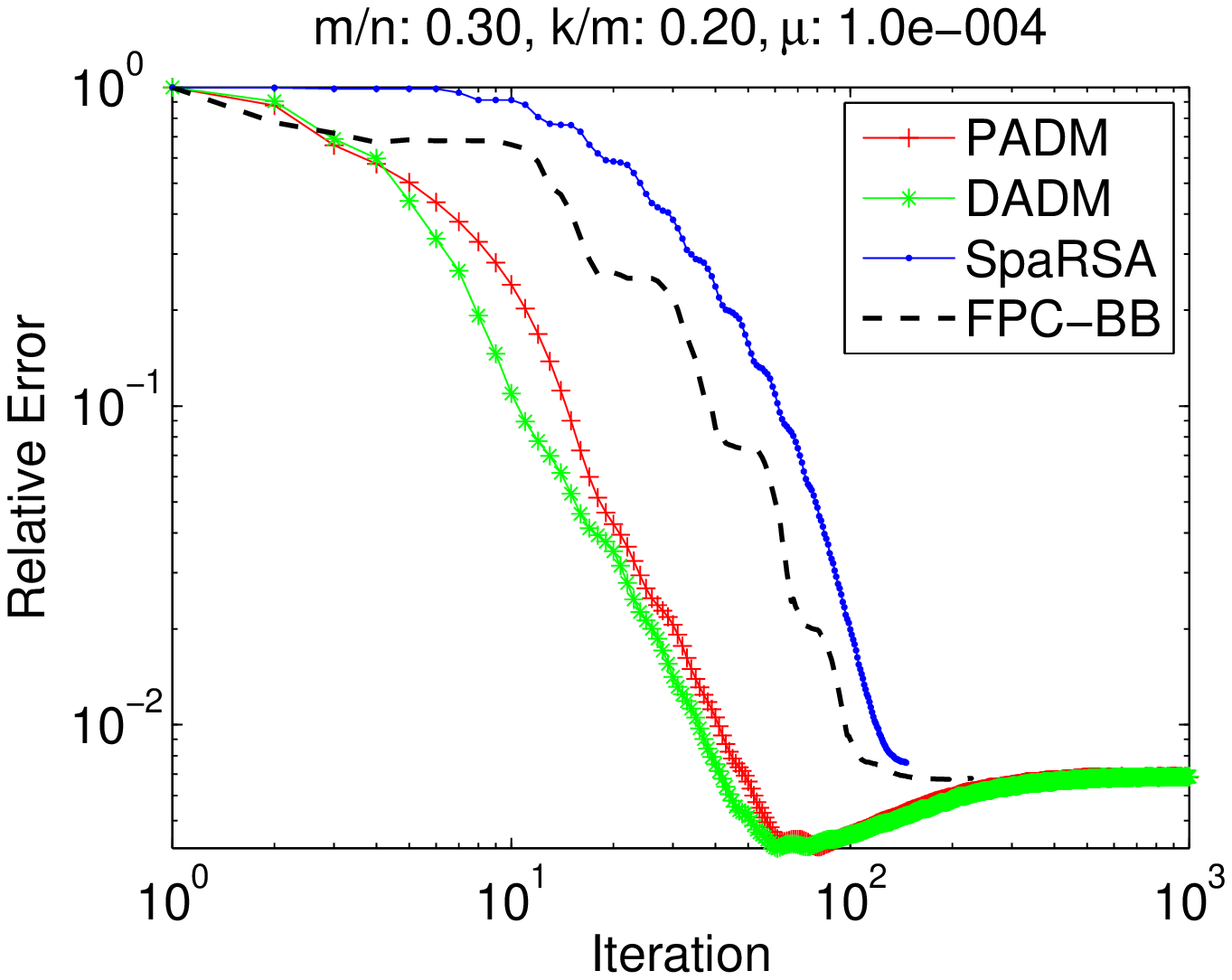}
\includegraphics[scale=.5]{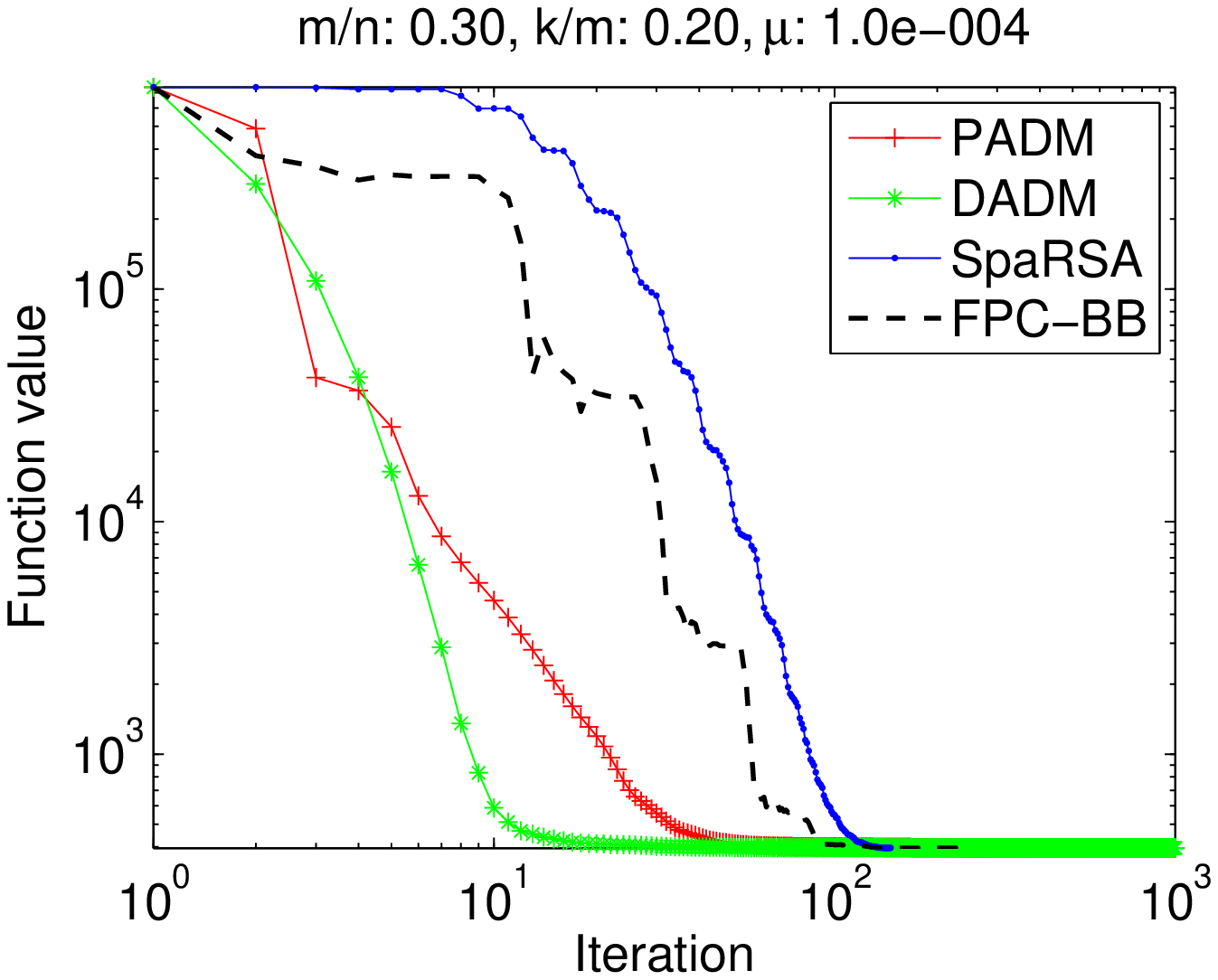}
\includegraphics[scale=.5]{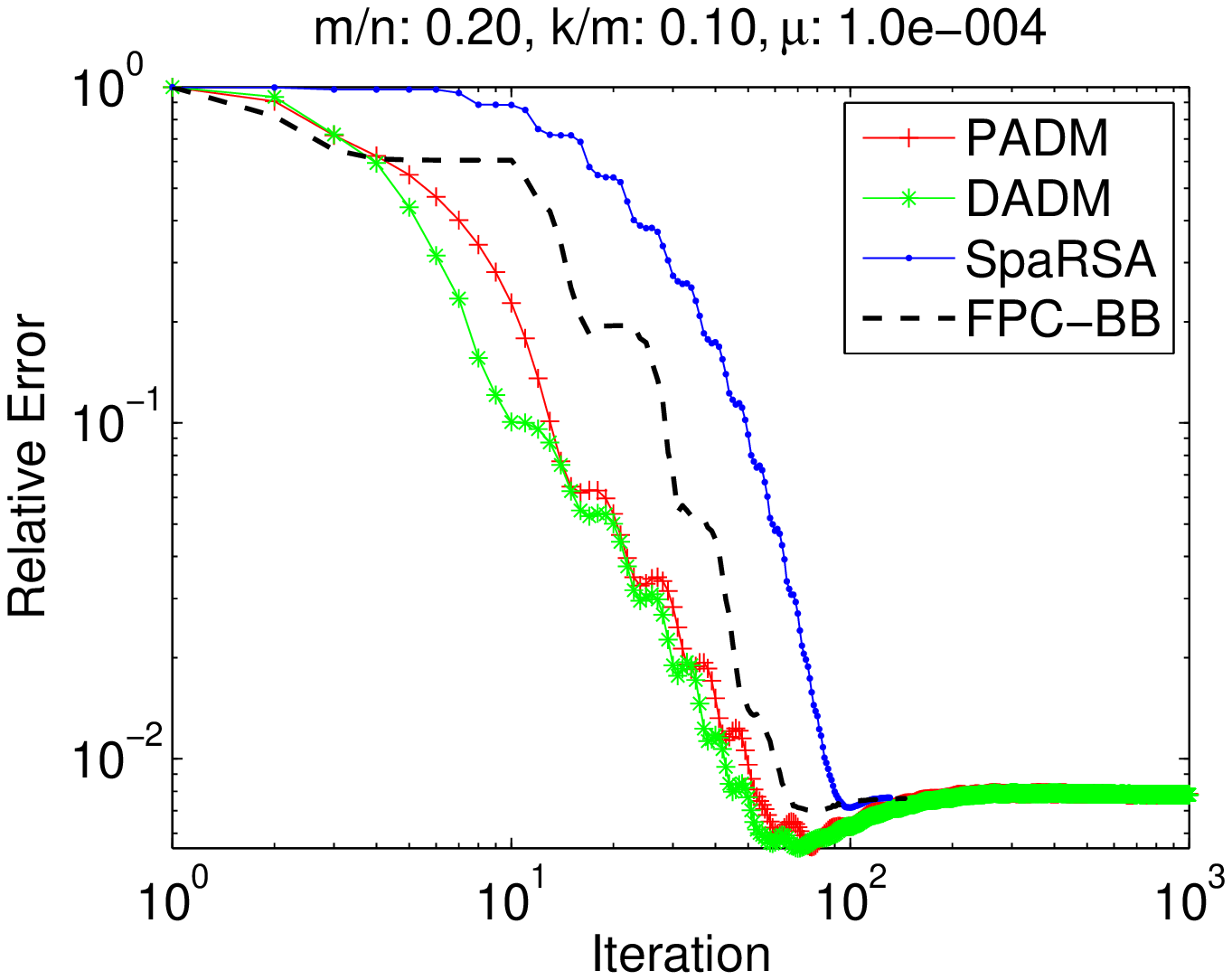}
\includegraphics[scale=.5]{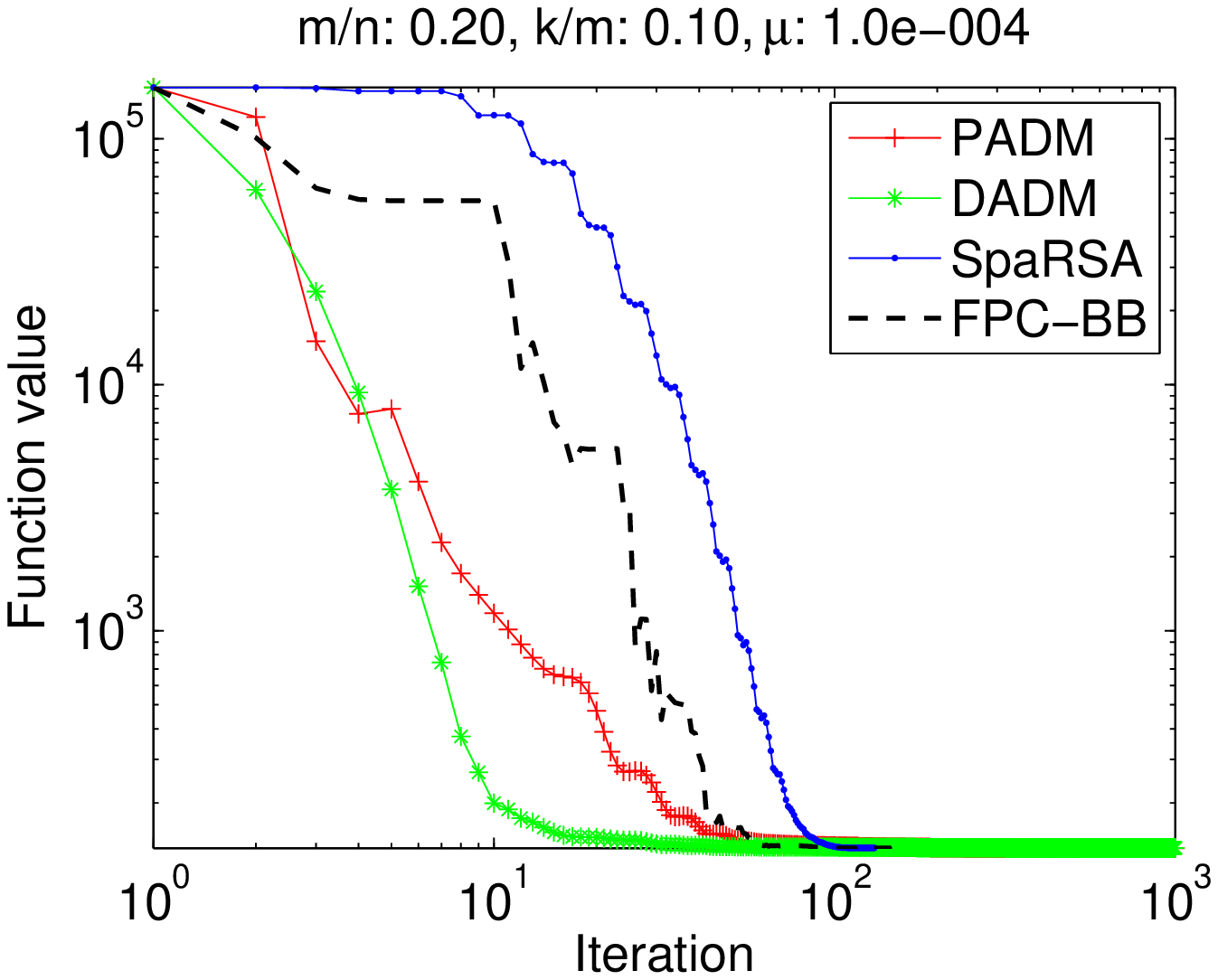}
} \caption{Comparison results of PADM, DADM, FPC-BB and SpaRSA on
\eqref{decoder-L1L2}. The $x$-axes represent number of iterations, and
$y$-axes represent relative errors (plots on the left) and function
values (plots on the right). The standard deviation of Gaussian
noise is {\tt sigma $=10^{-3}$}. The results are average of 50 runs.
\label{fig-QPmu}}
\end{figure}

For the purpose of clarity, we first compare PADM and DADM with
FPC-BB and SpaRSA, and then with FISTA and CGD. Since FPC-BB
implements continuation on regularization parameter but not on
stopping tolerance, we set all parameters as default except in the
last step of continuation we let {\tt xtol = $10^{-5}$} and {\tt
gtol = $0.02$}, which is more stringent than the default setting
{\tt xtol = $10^{-4}$} and {\tt gtol = $0.2$} because the latter
usually produces solutions of lower quality than that of other
algorithms in comparison.
For SpaRSA, we used its monotonic variant,
set continuation steps to be 20 and terminated it when relative
change in function values falls below $10^{-7}$. Since the
per-iteration cost is roughly two matrix-vector multiplications for
all compared algorithms, it is helpful to examine the
decreasing behavior of relative errors and function values as
functions of iteration numbers. Figure \ref{fig-QPmu} presents the
results of two different combinations of $m$ and $k$. Each result is
the average of 50 runs on randomly generated data.

As can be seen from Figure \ref{fig-QPmu} that, compared with FPC-BB
and SpaRSA, PADM and DADM usually decrease relative errors and
function values faster. Specifically, the  relative error and
function value curves of both PADM and DADM fall below those of
FPC-BB and SpaRSA almost throughout the entire iteration process.
With no more than 100 iterations, PADM and DADM reached lowest
relative errors and then started to increase, which is a problem of
the model \eqref{decoder-L1L2} rather than the algorithm since
function values keep decreasing. It is also clear that all
algorithms attain nearly equal relative errors and function values.

\begin{figure}[htbp]
\centering{
\includegraphics[scale=.5]{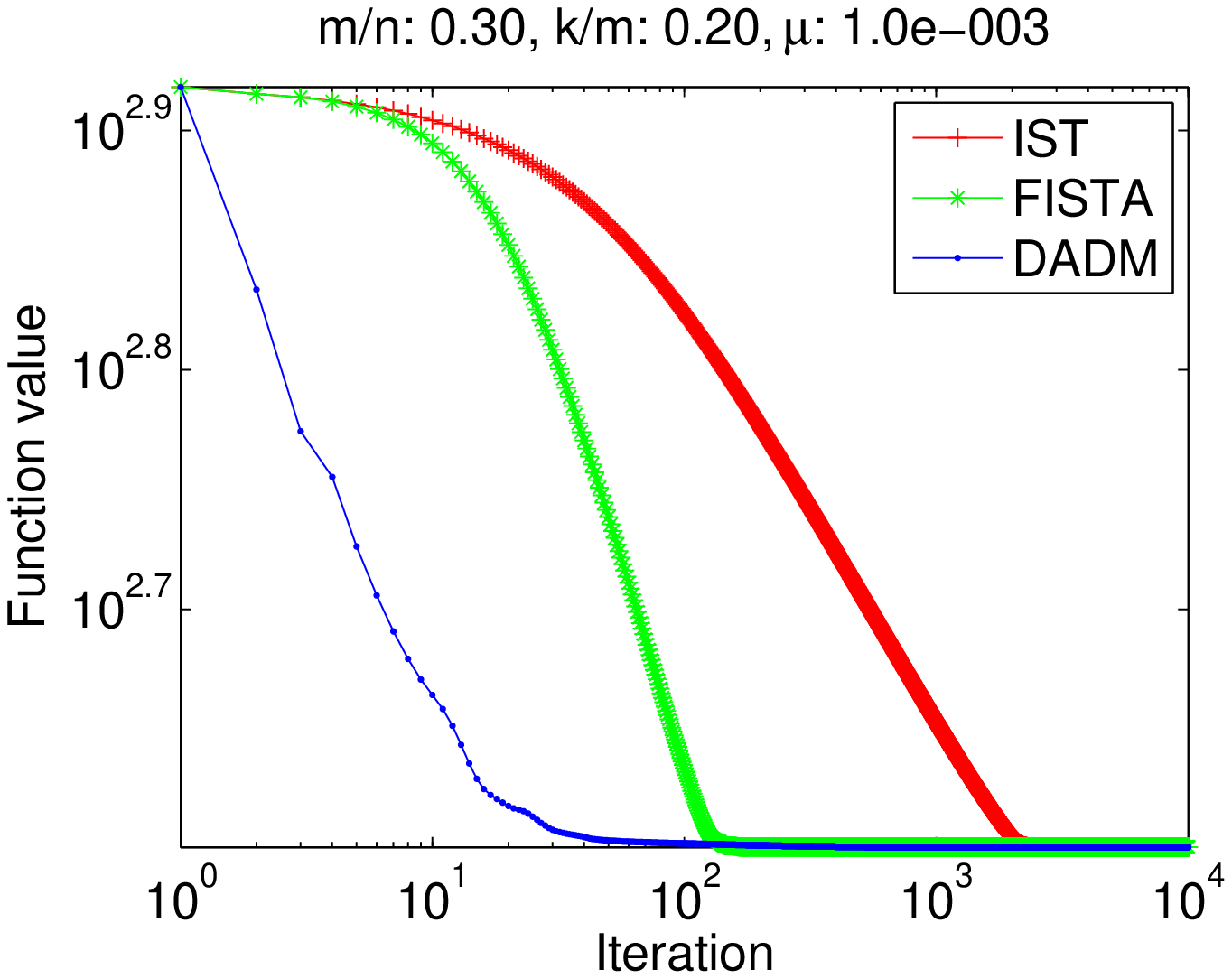}
\includegraphics[scale=.5]{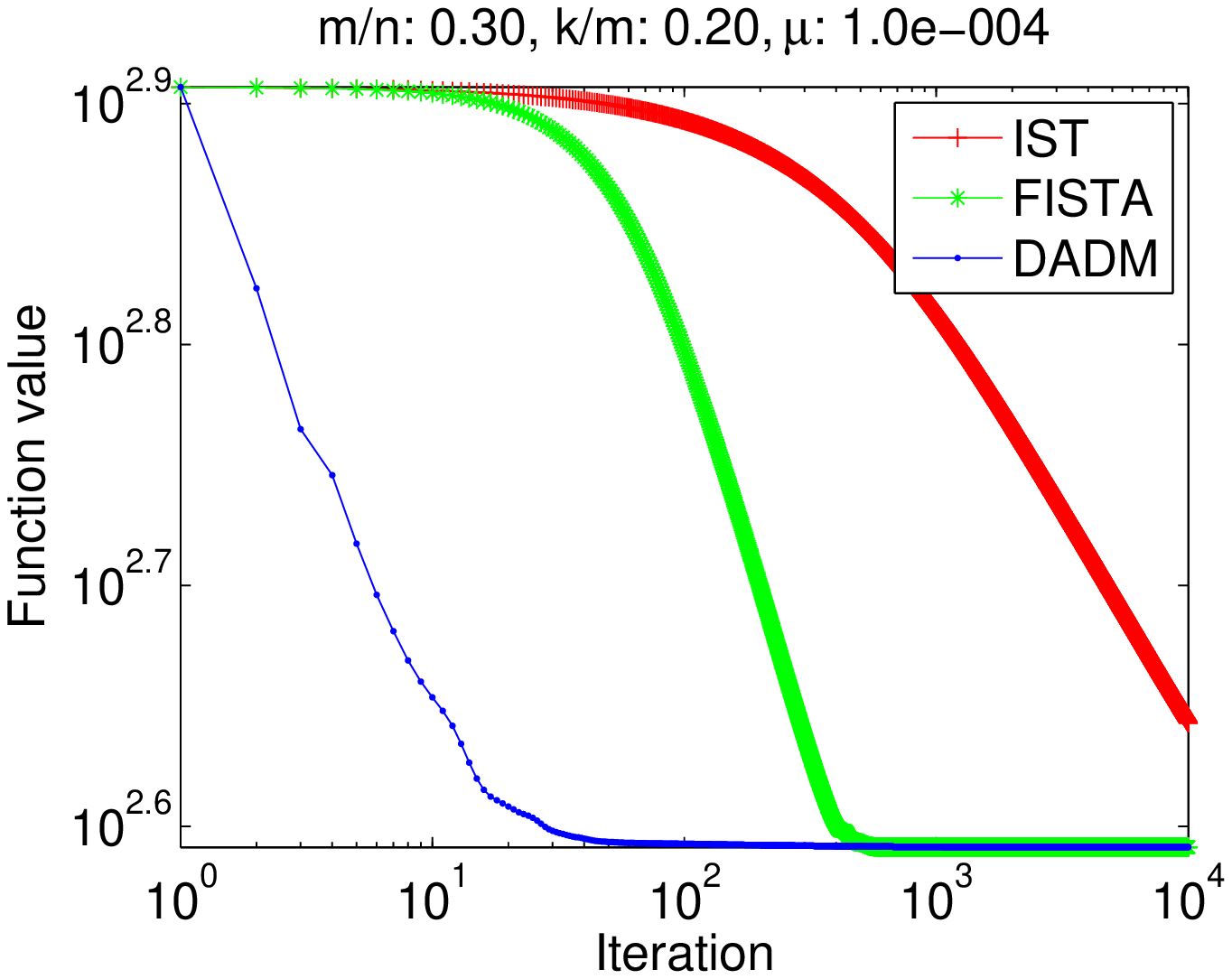}
} \caption{Comparison results of IST, FISTA and DADM on
\eqref{decoder-L1L2}. The $x$-axes represent number of iterations, and
$y$-axes represent function values, both in logarithmic scale.
The standard deviation of
Gaussian noise is $10^{-3}$. Left plot: $\mu=10^{-3}$; Right plot:
$\mu=10^{-4}$. The results are average of 50 runs.
\label{fig-C2Fista}}
\end{figure}

Next we compare DADM with FISTA  and CGD. Started at $x^0$, FISTA
iterates as follows
\begin{eqnarray*}\label{fista}
x^{k+1} = \text{Shrink}\left(y^k - \tau A^*(Ay^k-b),\tau/\mu\right),
\end{eqnarray*}
where $\tau>0$ is a parameter, and
\[
y^k=\left\{
      \begin{array}{ll}
        x^0, & \hbox{if $k=0$;} \\
        x^k + \frac{t_{k-1}-1}{t_{k}}(x^k-x^{k-1}), & \hbox{otherwise,}
      \end{array}
    \right.
\text{ where }\quad
t_k =\left\{
      \begin{array}{ll}
        1, & \hbox{if $k=0$;} \\
        \frac{1+\sqrt{1+4t_{k-1}^2}}{2}, & \hbox{otherwise.}
      \end{array}
    \right.
\]
It is shown in \cite{FISTA} that FISTA attains an optimal
convergence rate $O(1/k^2)$ in decreasing function values, where $k$
is the iteration counter. By letting $t_k\equiv1$, FISTA reduces to
the well-known primary iterative-soft-thresholding (IST) algorithm.
Figure \ref{fig-C2Fista} shows the
comparison results of DADM and FISTA, along with those of the
primary IST algorithm in order to illustrate the advantage of FISTA
over IST (after a small modification). Since the relative
performance of PADM and DADM on problem \eqref{decoder-L1L2} has
been illustrated in Figure \ref{fig-QPmu} and the dual approach
seems to be slightly more efficient than the primal one, here we merely
presented the results of DADM for simplicity. In this experiment, we
set $\tau=1$ for both FISTA and IST and initialized all algorithms
at $x=A^*b$ rather than at the origin in order to make function values
fall into a relatively small range.

\begin{figure}[htbp]
\centering{
\includegraphics[scale=.5]{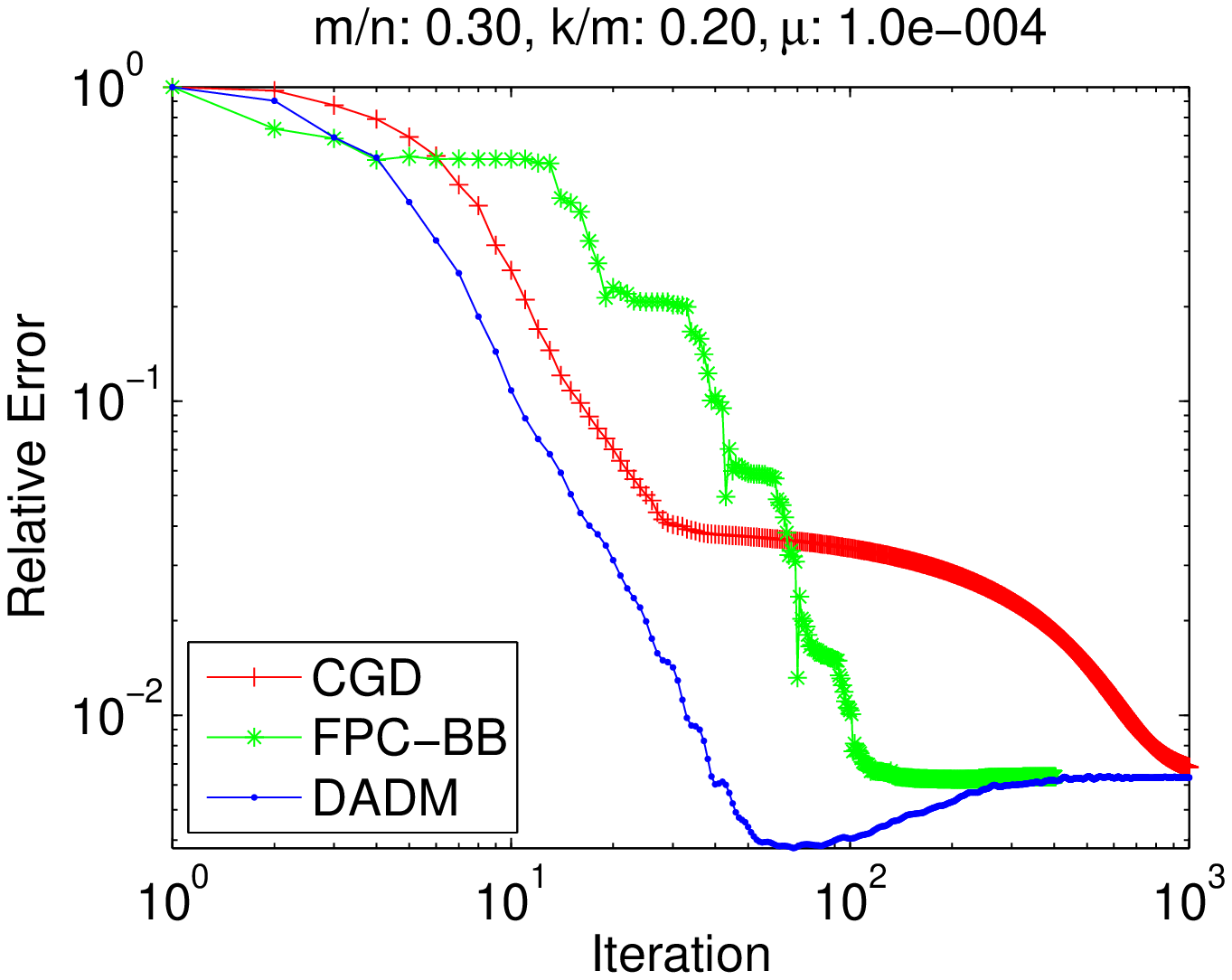}
\includegraphics[scale=.5]{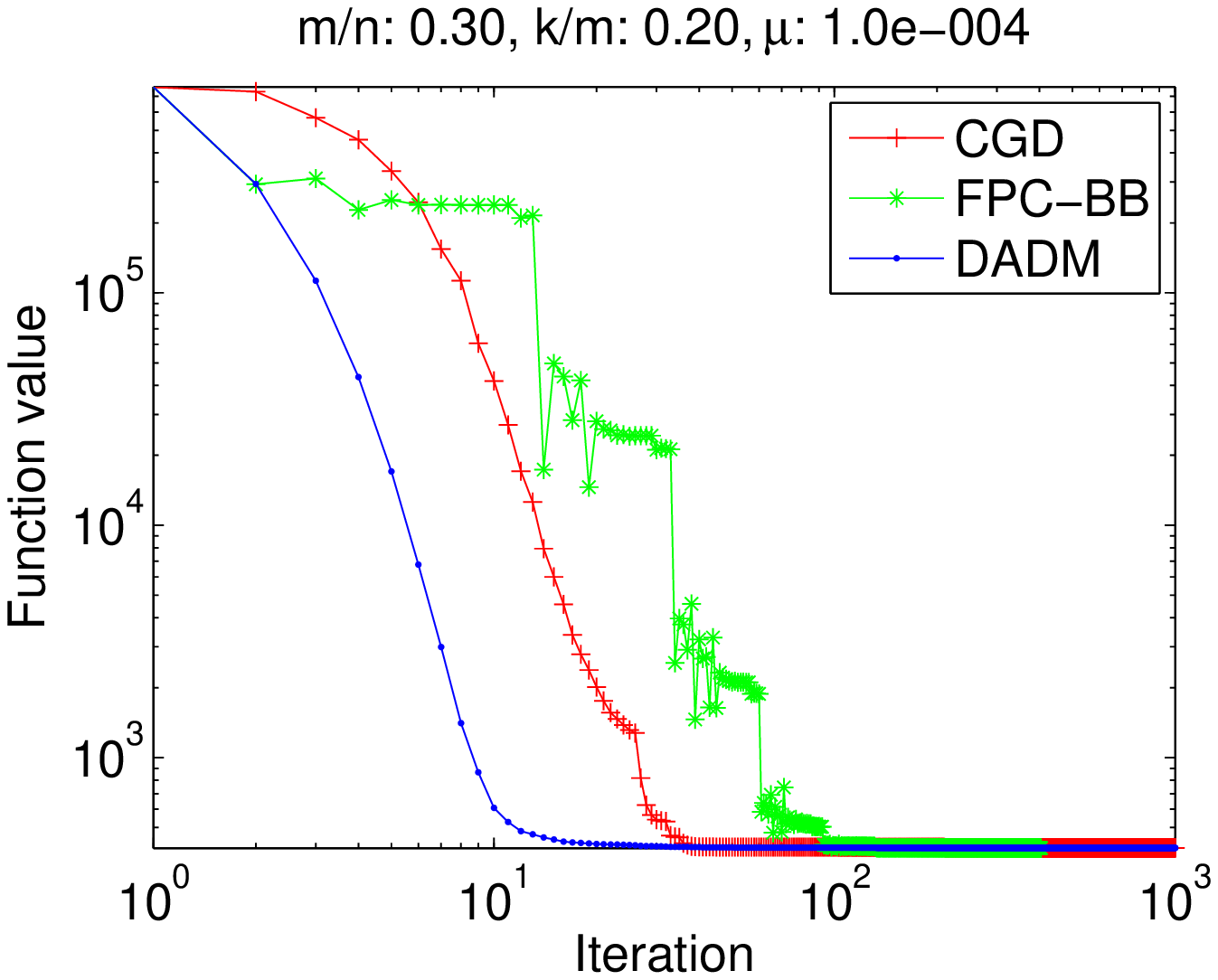}
} \caption{Comparison results of FPC-BB, CGD and DADM on
\eqref{decoder-L1L2}. The $x$-axes represent number of iterations, and
$y$-axes represent relative errors (plots on the left) and function
values (plots on the right), both in logarithmic scale.
The standard deviation of Gaussian
noise is $10^{-3}$. The results are average of 50 runs.
\label{fig-C2CGD}}
\end{figure}

It can be seen from Figure \ref{fig-C2Fista} that FISTA converges
much faster than IST. However, without heuristic techniques such as
continuation and line search as incorporated in many algorithms,
FISTA is generally slower than DADM. The slow convergence of FISTA
and IST as compared with DADM becomes more pronounced when $\mu$
becomes smaller in which case \eqref{decoder-L1L2} becomes more
difficult. From Figure \ref{fig-C2Fista}, when $\mu=10^{-3}$, FISTA
converges to a nearly optimal function value within no more than 200
iterations, while IST consumes more than 1000 iterations. When $\mu$
is decreased from $10^{-3}$ to $10^{-4}$, both FISTA and IST become
slower, while the convergence of DADM is not affected.

Similarly, the comparison results of DADM with CGD are given in
Figure \ref{fig-C2CGD}, along with those of FPC-BB in order to
evaluate the relative performance of CGD.  In this experiment, we
used the continuation variant of CGD (the code {\tt CGD$\_$cont} in
the {\tt MATLAB} package of CGD) and set all parameters as default
except setting the initial $\mu$ value to be  $\max( 0.01\|A^\top
x^0\|_\infty,2\mu)$ which works better than the default setting in
our tests when $\mu$ is small.    It can be seen from the plot on
the right of Figure \ref{fig-C2CGD} that CGD decreases the function
value faster than FPC-BB. However, from the plot on the left, CGD
does not necessarily decrease the relative error to the true signal
faster
than FPC-BB. 
In comparison, DADM converges faster in both the function
value and relative error.

We experimented on various combinations of $(m, k)$ with noisy data and
observed similar phenomenon.  As is the case in Figures \ref{fig-QPmu}
and \ref{fig-C2CGD}, the relative error produced by the ADM algorithms
tends to eventually increase, after the initial decrease when
problem \eqref{decoder-L1L2} is solved to a high accuracy.
This implies, as suggested in Section \ref{sc:Rerr-vs-opt}, that
it is unnecessary to solve $\ell_1$-problems to a higher accuracy
than what is warranted by the accuracy of the underlying data.
Therefore, in the next set of tests we terminate PADM and DADM
whenever \eqref{def:RelChg} is satisfied with
$\epsilon=5\times10^{-4}$ instead of using much smaller tolerances.

Next we compare PADM and DADM with FPC-BB and SpaRSA for various
combinations of $(m,k)$. Here we do not present results for FISTA
and CGD because they have been found to be less competitive in this
set of tests.  As is already mentioned earlier, without heuristic
continuation and line search techniques, FISTA is slower than ADM
algorithms.  On the other hand,  CGD is slower in terms of
decreasing relative error.
We set all parameters as default in FPC-BB and use the same setting
as before for SpaRSA except it is terminated when relative change in
function values falls below $10^{-4}$.  For each fixed pair $(m,k)$,
we take the average of 50 runs on random instances. Detailed results
including iteration number (Iter) and relative error to the true sparse
signal (RelErr) are given in Table \ref{Table-QPmu}.

\begin{table}[ht]
\centering \caption{Comparison results on \eqref{decoder-L1L2} ({\tt
sigma $= 10^{-3}$}, $\mu=10^{-4}$, average of 50 runs).  }
\vspace{-.2cm}
\begin{tabular}{c|c||c|c|c|c|c|c|c|c}\hline

\multicolumn{2}{c||}{n =
8192}&\multicolumn{2}{c|}{PADM}&\multicolumn{2}{c|}{DADM}&\multicolumn{2}{c|}{SpaRSA}&\multicolumn{2}{c}{FPC-BB}\\\hline

m/n&$k/m$&Iter&RelErr&Iter&RelErr&Iter&RelErr&Iter&RelErr\\
\hline\hline

0.3&0.1& 67.4&3.93E-3 & 55.6&3.96E-3 & 103.3&5.70E-3 & 55.1&5.88E-3
\\\hline

0.3&0.2 &  76.4& 4.03E-3  & 68.3 & 4.09E-3 & 139.4 & 7.12E-3 & 93.1
& 7.37E-3
\\\hline

0.2&0.1 & 102.2&6.34E-3 & 100.4 & 6.17E-3 & 114.2&7.47E-3 & 69.8&
7.56E-3
\\\hline

0.2&0.2&  103.3 & 7.90E-3  & 96.2& 7.72E-3 & 178.4 & 1.59E-2 & 121.3
& 2.32E-2
\\\hline

0.1&0.1 & 154.2&1.23E-2 & 167.5&1.22E-2 & 136.8&1.30E-2 & 83.5&
1.37E-2
\\\hline

0.1&0.2&  252.4& 1.09E-1  & 239.4 & 1.10E-1 & 210.2 & 1.65E-1 &
126.2& 2.00E-1
\\\hline
\end{tabular}
\label{Table-QPmu}
\end{table}

As can be seen from Table \ref{Table-QPmu}, in all cases PADM and
DADM obtained smaller relative errors in comparable numbers of
iterations as FPC-BB and SpaRSA. This is particularly true for more
difficult problems, e.g., for $m/n = 0.1$  and $k/m = 0.2$ where the
resulting relative errors of PADM and DADM are considerably smaller
than those of FPC-BB and SpaRSA.   We also tried to terminated
SpaRSA and FPC-BB using more stringent stopping rules. Specifically,
we set {\tt xtol}=$10^{-5}$ and {\tt gtol}=0.02 in FPC-BB and
terminated SpaRSA when relative change in function values falls
below $10^{-7}$. The relative error results  either remain roughly
the same as those presented in Table \ref{Table-QPmu} or were just
slightly better, while the iteration numbers required by both FPC-BB
and SpaRSA were increased from around 50\% to more than 100\%.  We
also tested partial DCT and sparse signals of different dynamic
ranges and observed similar phenomenon.

\subsection{Comparison  with SPGL1 and NESTA}\label{sc:compBPdelta}

\begin{figure}[htbp]
\centering{
\includegraphics[scale=.5]{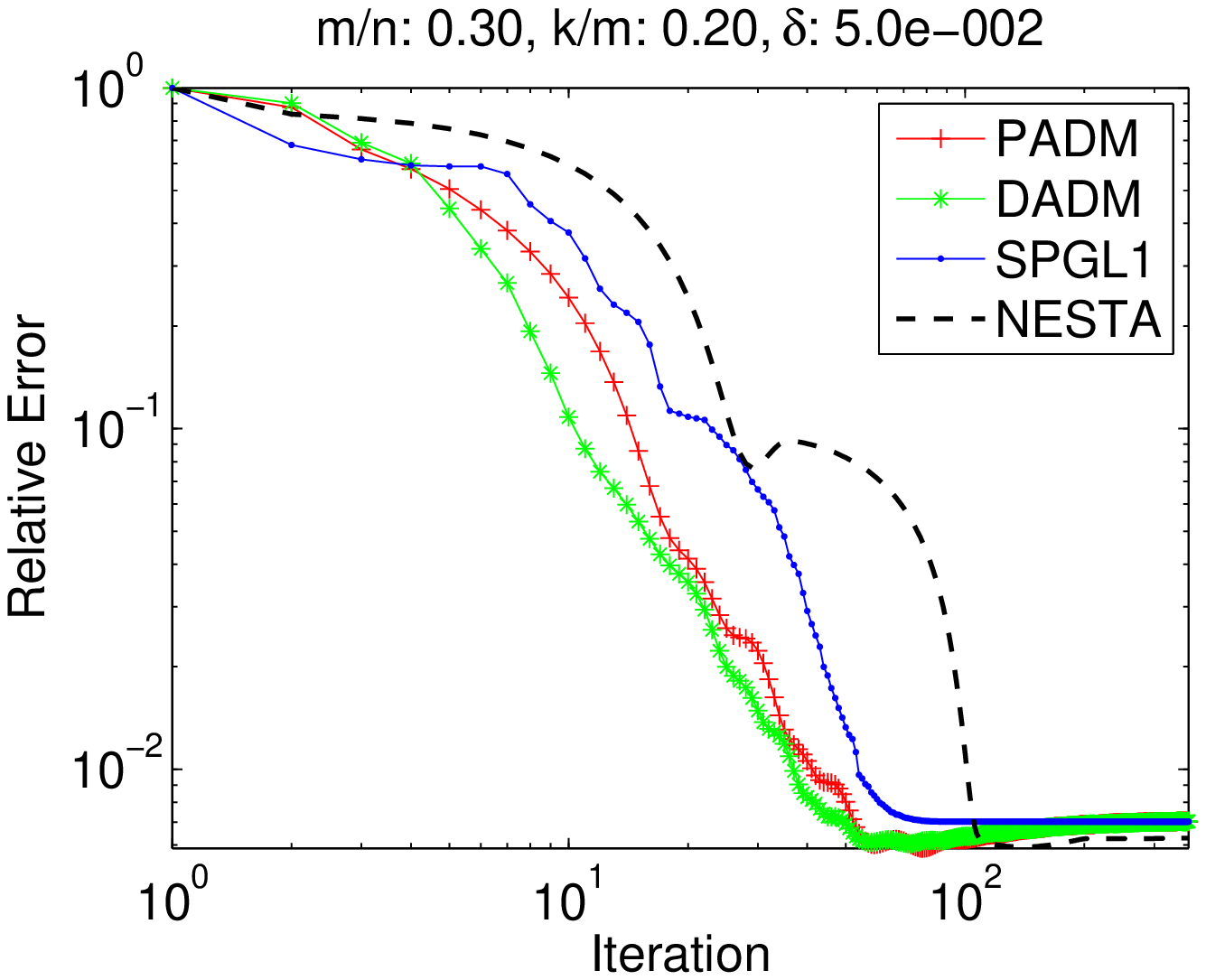}
\includegraphics[scale=.5]{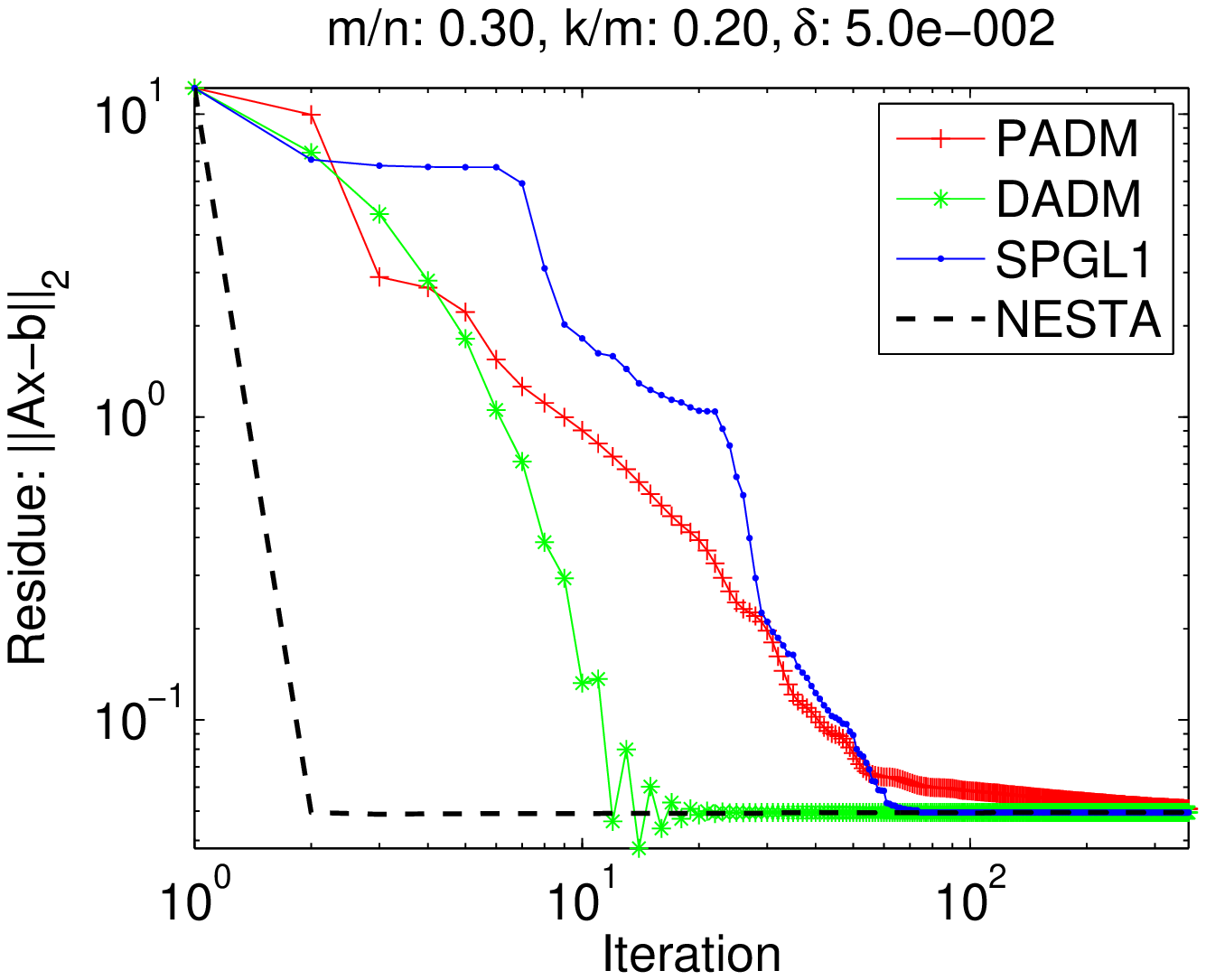}
\includegraphics[scale=.5]{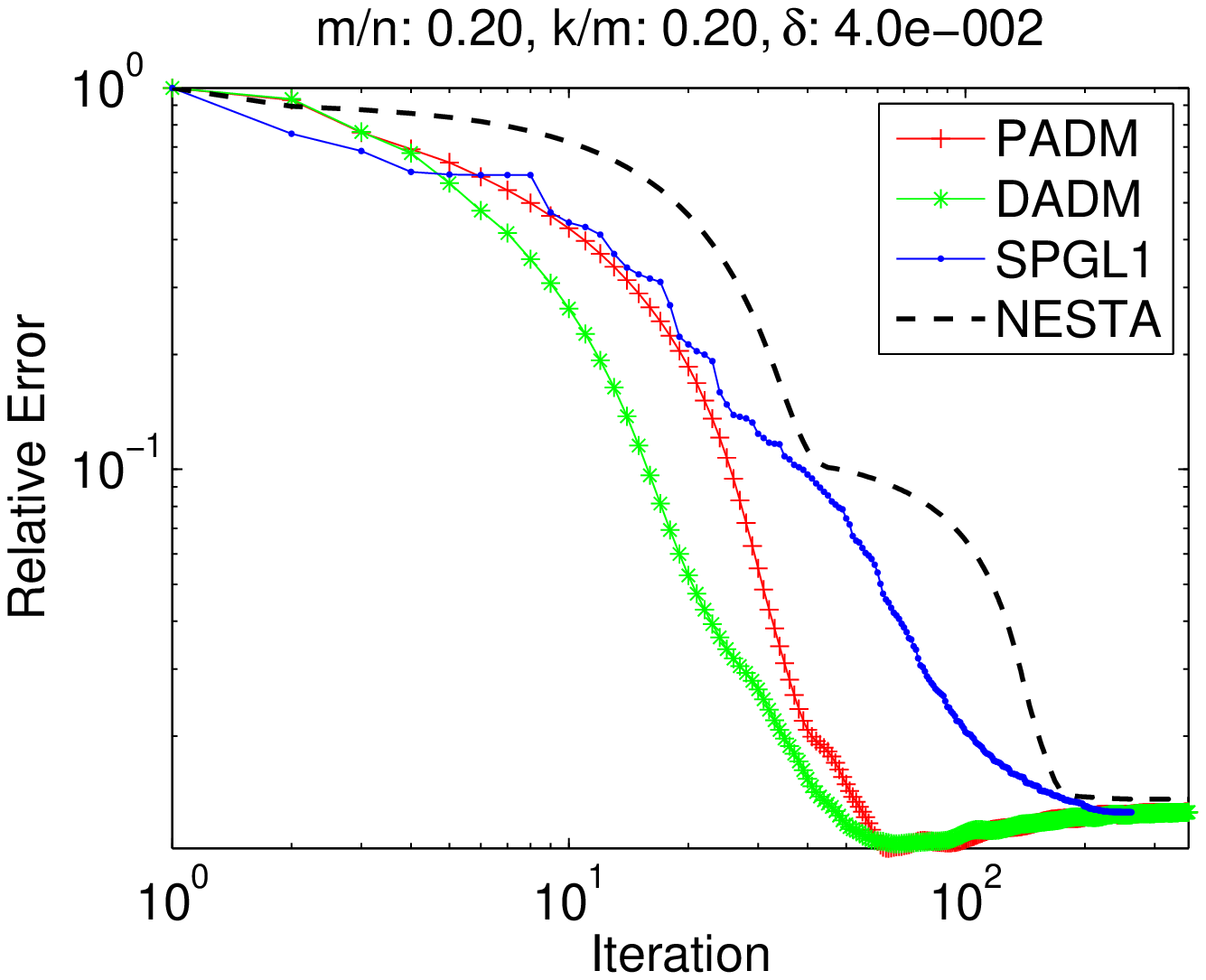}
\includegraphics[scale=.5]{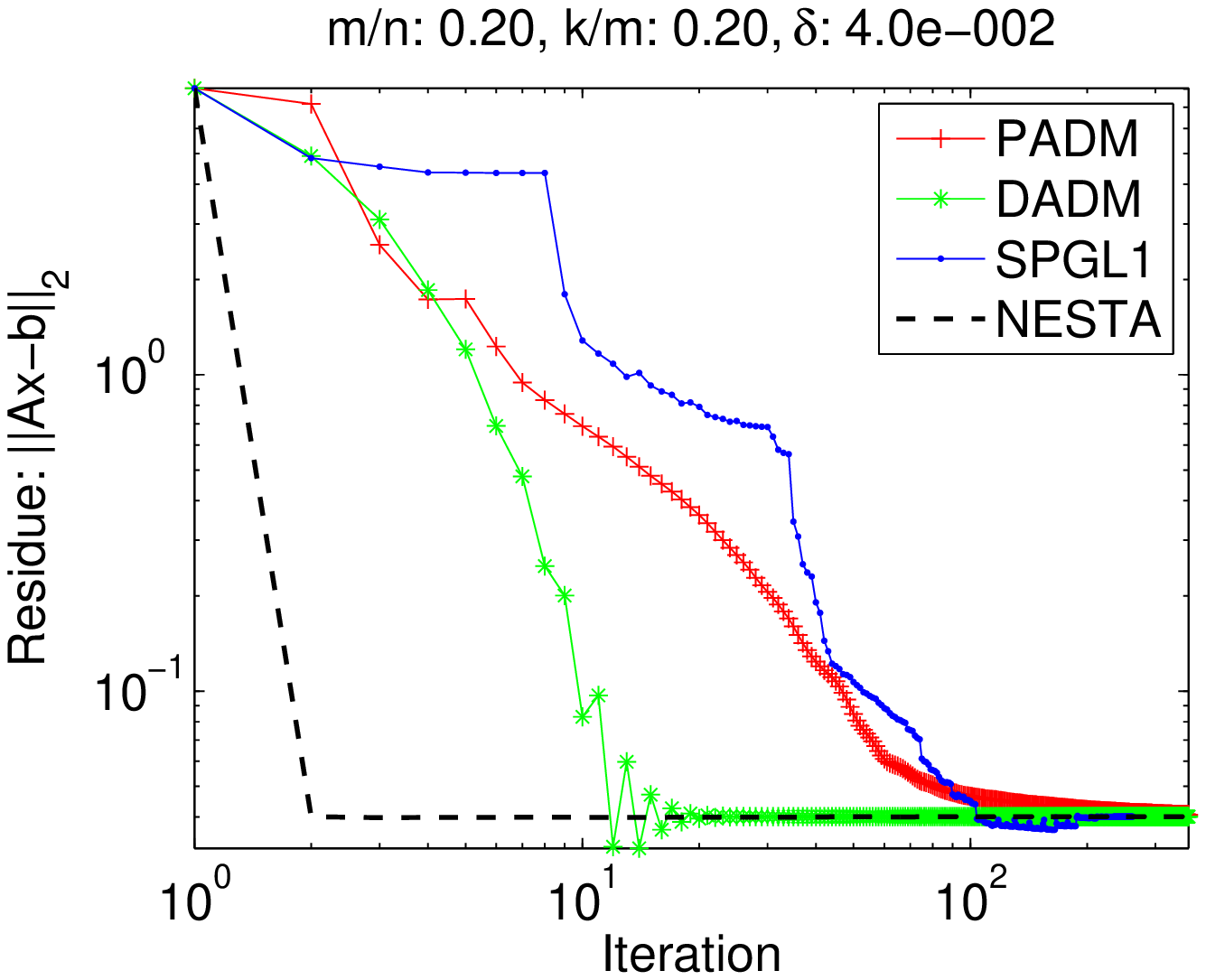}
} \caption{Comparison results of PADM, DADM, SPGL1 and NESTA on
\eqref{decoder-BPDN}. The $x$-axes represent number of iterations;
$y$-axes represent relative error (plots on the left) and residue
(plots on the right). The standard deviation of Gaussian noise is
$10^{-3}$. The results are average of 50 runs. \label{fig-BPdelta}}
\end{figure}

In this subsection, we compare PADM and DADM with SPGL1 and NESTA
for solving \eqref{decoder-BPDN}. The same as before, {\tt xbar} is
random Gaussian spikes, and  the standard deviation {\tt sigma} of
additive noise  is $10^{-3}$. Parameter $\delta$ in
\eqref{decoder-BPDN} was set to be the 2-norm of additive noise
(the ideal case).  As in the previous experiment, we performed two sets of tests.
In the first set, we ran PADM and DADM to a prescribed number of iterations
and terminated SPGL1 and NESTA with suitably chosen tolerances so that the
iteration numbers consumed are roughly identical.  Specifically, we set {\tt
optTol} $= 10^{-8}$ in SPGL1, which is a tolerance for optimality,
and {\tt TolVar} $= 10^{-8}$ in NESTA, which is a tolerance in
relative change in objective function. All other parameters
are set to their default values.  Figure \ref{fig-BPdelta} presents average
results of 50 random problems, where two combinations of $m$ and $k$
are used. The resulting relative error and fidelity residue (i.e., $\|Ax-b\|$)
are plotted as functions of iterations.

As can be seen from the left column of Figure \ref{fig-BPdelta}, compared with
SPGL1 and NESTA, both PADM and DADM attained smaller relative errors
throughout most of the iteration process (with the exception at the very beginning).
With no more than 100 iterations, both PADM and DADM reached lowest relative
errors and then started to increase slightly.
It is interesting to note that, as can be seen from the right column
of Figure \ref{fig-BPdelta}, NESTA is the fastest in terms of decreasing
fidelity residue but slowest in terms of decreasing the relative error.
In the applications of compressive sensing, we are interested in
recovering the true signal, which means that the smaller the
relative error is, the better.

After extensive experiments under various scenarios, including different
sensing matrices and sparse signals of different dynamic ranges, we have
found that when data are noisy the proposed ADM algorithms decrease
relative error faster than all other tested algorithms.

In the second set of tests, we terminated PADM and DADM  by
\eqref{def:RelChg} with $\epsilon = 5\times10^{-4}$. For SPGL1 and
NESTA, we set all parameters as default except {\tt TolVar} is set
to be $10^{-6}$ in NESTA (where the default value is $10^{-5}$) to obtain
solutions of comparable quality. The average results on 50 random
problems are given in Table \ref{Table-BPdelta}. As mentioned
before, matrix-vector multiplications are the main computational load
for all compared algorithms.  The number of matrix-vector multiplications
consumed by PADM and DADM is two times the number of iterations.
However, the number used by SPGL1 is not always proportional to the
number of iterations.   As for NESTA, its per-iteration cost is also two
matrix-vector multiplications for partial orthonormal matrices,
although this number increases to six for general matrices.
Therefore, instead of iteration numbers, we present in Table
\ref{Table-BPdelta} the number of matrix-vector multiplications,
denoted by \#AAt that includes both {\tt A*x} and {\tt A'*y}.

\begin{table}[ht]
\centering \caption{Comparison results on \eqref{decoder-BPDN} ({\tt
sigma $= 10^{-3}$}, $\delta=$ {\tt norm(noise)}, average of 50
runs).} \vspace{-.2cm}
\begin{tabular}{c|c||c|c|c|c|c|c|c|c}\hline

\multicolumn{2}{c||}{n =
8192}&\multicolumn{2}{c|}{PADM}&\multicolumn{2}{c|}{DADM}&\multicolumn{2}{c|}{SPGL1}&\multicolumn{2}{c}{NESTA}\\\hline

m/n&$k/m$&\#AAt&RelErr&\#AAt&RelErr&\#AAt&RelErr&\#AAt&RelErr\\
\hline\hline

0.3&0.1& 107.9& 6.40E-3  & 97.8& 6.35E-3  & 120.8& 5.39E-3& 300.6 &
5.83E-3
\\\hline

0.3&0.2& 118.1 &  6.05E-3  & 108.2& 6.26E-3 & 160.1& 7.26E-3 &
303.6& 8.33E-3
\\\hline

0.2&0.1& 124.9& 8.01E-3  & 112.7& 7.92E-3 & 149.7& 7.35E-3  & 316.9&
6.85E-3
\\\hline

0.2&0.2& 122.6 & 1.05E-2  & 110.0& 1.05E-2  & 168.2& 1.62E-2 &
311.5& 6.41E-2
\\\hline

0.1&0.1& 167.3& 1.29E-2  & 152.2& 1.26E-2  &  171.9&  1.29E-2  &
340.6& 1.73E-2
\\\hline

0.1&0.2& 161.2& 9.86E-2& 145.0& 1.06E-1  & 184.0& 1.49E-1 & 326.4&
2.96E-1
\\\hline

\end{tabular}
\label{Table-BPdelta}
\end{table}

As can be seen from Table \ref{Table-BPdelta}, compared with SGPL1
and NESTA, both PADM and DADM obtained solutions of comparable
quality within smaller numbers of matrix-vector multiplications.

\begin{table}[ht]
\centering \caption{Comparison results on \eqref{decoder-L1} ($b$ is
noiseless; stopping rule: $\epsilon=10^{-6}$ in \eqref{def:RelChg};
average of 50 runs).} \vspace{-.2cm}
\begin{tabular}{c|c||c|c|c|c|c|c|c|c}\hline

\multicolumn{2}{c||}{n =
8192}&\multicolumn{4}{c|}{DADM}&\multicolumn{4}{c}{SPGL1}\\\hline

m/n&$k/m$&RelErr&RelRes&CPU&\#AAt&RelErr&RelRes&CPU&\#AAt\\
\hline\hline
0.3&0.1& 7.29E-5& 4.41E-16  & 0.44& 258.8  & 1.55E-5& 9.19E-6& 0.39
& 114.9
\\\hline
0.3&0.2& 7.70E-5 & 4.65E-16  & 0.78& 431.4  & 2.50E-5& 6.77E-6& 1.11
& 333.4
\\\hline
0.2&0.1& 4.26E-5& 4.54E-16  & 0.66& 388.2  & 3.39E-5& 1.51E-5& 0.45
& 146.7
\\\hline
0.2&0.2& 7.04E-5& 4.85E-16  & 1.15& 681.8  & 1.40E-4& 1.03E-5& 2.50
& 791.0
\\\hline
0.1&0.1& 4.17E-5& 4.86E-16  & 1.11& 698.2  & 1.25E-4& 3.26E-5& 0.64
& 207.9
\\\hline

\end{tabular}
\label{Table-BP}
\end{table}

\subsection{Comparison with SPGL1 on the basis pursuit problem}
\label{sc:compBP}
In this subsection, we compare DADM with SPGL1 on the basis pursuit
problem \eqref{decoder-L1}, in which case $b$ is noiseless and
solving problems to a higher accuracy should result in higher-quality
solutions. Thus, we terminated DADM with a stringent stopping
tolerance of $\epsilon=10^{-6}$ in \eqref{def:RelChg}.  All
parameters in SPGL1 are set to be default values. Detailed
comparison results are given in Table \ref{Table-BP}, where, besides
relative error ({\tt RelErr}) and the number of matrix-vector
multiplications (\#AAt), the relative residue ${\tt
RelRes}=\|Ax-b\|/\|b\|$ and CPU time (in seconds) are also given.

As can be observed from previous experimental results, the proposed
ADM algorithms may slow down after a fast convergence stage at the
beginning. When measurements are free of noise and a highly accurate
solution is demanded, the ADM algorithms can sometimes be slower
than SPGL1. Table \ref{Table-BP} shows that DADM is slower than
SPGL1 in three of the five test cases in terms of CPU seconds while
getting slightly lower accuracy (the results for $m/n=0.1$ and
$k/m=0.2$ are omitted since both algorithms failed to recover an
accurate solution). We note that since SPGL1 requires some
non-trivial calculations other than matrix-vector multiplications, a
larger \#AAt number by DADM does not necessarily lead to a longer
CPU time.   We also comment that the relative residue results of
DADM are always numerically zero because when $AA^*=I$ the sequence
$\{x^k\}$ generated by DADM, applied to \eqref{decoder-L1},
satisfies $Ax^{k+1}-b=(1-\gamma)(Ax^k-b)$ and thus $\|Ax-b\|$
decreases geometrically.

\subsection{Summary}
We  provided supporting evidence to emphasize the important point that
algorithm speed should be evaluated relative to solution accuracy.
Since more often than not measurements are noisy applications,
solving $\ell_1$-problems to a very high accuracy is generally
not warranted and unnecessary.   It is more practically relevant to
evaluate the speed of an algorithm based on how fast it achieves
an appropriate accuracy that is consistent with noise levels in data.

We presented extensive experimental results on $\ell_1$-problems and
compared the proposed first-order, primal-dual algorithms,
derived from the ADM approach, with state-of-the-art
algorithms FPC-BB, SpaRSA, FISTA, CGD, SPGL1 and NESTA.
Our numerical results show that on the tested cases the proposed
ADM algorithms are efficient and stable.
In particular, under practically relevant conditions where
data contain noise and stopping tolerances are set to appropriate values
(where a more stringent tolerance would not lead to a more
accurate solution), the proposed ADM algorithms generally achieve
lower relative errors within fewer or comparable number of iterations
in comparison to other tested algorithms.

The experimental results also indicate that the dual-based ADMs are
generally more efficient than the primal-based ones.  One plausible
explanation is that when $A$ is orthonormal, the dual-based algorithms are
exact ADMs, while the primal ones are inexact which solve some subproblem
approximately.  The dual-based ADMs have been implemented in a
{\tt MATLAB} package called
YALL1 \cite{YALL1} (short for Your ALgorithm for L1), which is
applicable to eight different $\ell_1$-models including
\eqref{decoder-L1}, \eqref{decoder-BPDN}, \eqref{decoder-L1L2},
\eqref{decoder-L1L1} and the corresponding nonnegative counterparts.

\section{Concluding remarks}
We proposed to solve $\ell_1$-problems arising from compressive sensing
by first-order, primal-dual algorithms derived from the Alternating Direction
Method (ADM) framework which is based on the classic augmented Lagrangian
function and alternating minimization idea for structured optimization.
This ADM approach is applicable to numerious $\ell_1$-problems including
\eqref{decoder-L1}, \eqref{decoder-BPDN}, \eqref{decoder-L1L2},
\eqref{decoder-L1L1} and their nonnegative
counterparts, as well as others.   When applied to the $\ell_1$-problems,
the per-iteration cost of these algorithms is dominated by two
matrix-vector multiplications. Extensive experimental results show
that the proposed ADM algorithms, especially the dual-based ones,
perform competitively with several state-of-the-art algorithms.
On various classes of test problems with noisy data, the
proposed ADM algorithms have unmistakably exhibited the following
advantages over competing algorithms in the comparison:
(i) they converge relatively fast without the help of a continuation or a
line search technique;
(ii) their performance is relatively insensitive to changes in model and
algorithmic parameters; and
(iii) they demonstrate a notable ability to quickly decrease relative error to
true solutions.
Although the ADM algorithms are not necessarily the fastest
in reaching an extremely high accuracy when observed data
are noiseless, they are arguably the most effective in obtaining
the best achievable level of accuracy whenever data contain a nontrivial
level of noise, which is the case most relevant to practical applications.

The most influential feature of the ADM approach is perhaps its great
versatility and its seemingly universal effectiveness for a wide
range of optimization problems in signal, image and data analysis,
particular those involving $\ell_1$-like regularizations such as
nuclear-norm (sum of singular values) regularization in matrix rank
minimization like the matrix completion problem
\cite{RFP-MC07,CR-MC08, CT-MC09}, or the total variation (TV)
regularization \cite{ROF92} widely used in image processing. While
the nuclear-norm is just an extension of $\ell_1$-norm to the matrix
case, the TV regularization can be converted to
$\ell_1$-regularization after introducing a splitting variable
\cite{FTVd08, RecPF}. Therefore, the ADM approach is applicable to
both nuclear-norm and TV regularized problems (in either primal or
dual form) in a rather straightforward manner so that the
derivations and discussions are largely analogous to those for
$\ell_1$-problems as presented in this paper. Recently, the ADM has
also been applied to total variation based image reconstruction in
\cite{ADM09, RecPF, TVAL3} and to semi-definite programming in
\cite{Wen-Yin-Goldfarb09}.
A more recent application of the ADM approach is to the problem
of decomposing a given matrix into a sum of a low-rank matrix and
a sparse matrix simultaneously using $\ell_1$-norm and nuclear-norm
regularizations (see \cite{LRSD09}).  An ADM scheme has been proposed
and studied for this problem in \cite{YuanLRSD09}.

Although the ADM approach is classic and its convergence properties
have been well studied, its remarkable effectiveness in signal and image
reconstruction problems involving $\ell_1$-like regularizations
has just been recognized very recently.   These fruitful new
applications bring new research issues, such as convergence of
certain inexact ADM schemes, that should be interesting
for further investigations.


\section*{Acknowledgments} The first author would like to thank
Prof. Bingsheng He of Nanjing University and Dr. Wotao Yin of Rice University
for helpful discussions.

\bigskip
\appendix
\section{Proof of Theorem~\ref{conv-primal}}
\begin{proof}
Let $(\tilde{r},\tilde{x})$ be any solution of
\eqref{decoder-L1L2-2}. From optimization theory, there exists
$\tilde{y}\in\C^m$ such that the following conditions are satisfied:
\begin{eqnarray}\label{opt-primal}
    \tilde{r}/\mu - \tilde{y} = 0, \:
    A^* \tilde{y}\in \partial\|\tilde{x}\|_1 \text{\: and \:}
    A\tilde{x} + \tilde{r} = b.
\end{eqnarray}
For convenience, we let $\hat r \triangleq r^{k+1}$, $\hat x
\triangleq x^{k+1}$ and $\hat{y} \triangleq y^k - \beta(A\hat{x} +
\hat{r} - b)$. As such, $y^{k+1} = y^k - \gamma (y^k - \hat y)$.
From the definition of $\hat r$, $\hat x$ and $\hat y$, equation
\eqref{iter-s} can be reformulated as $\hat{r}/\mu - \hat{y} + \beta
A(x^k - \hat{x}) = 0$. Further considering $\tilde{r}/\mu -
\tilde{y} = 0$, we have $\hat{y} - \tilde{y} - \beta A(x^k -
\hat{x}) = (\hat{r}-\tilde{r})/\mu$, and thus
\begin{eqnarray}\label{opt&iter-s}
(\hat{r} -\tilde{r})^*\left(\hat{y} - \tilde{y} - \beta A(x^k -
\hat{x})\right) = \|\hat{r}-\tilde{r}\|^2/\mu\geq 0.
\end{eqnarray}
Similarly, equation \eqref{iter-x} is equivalent to $A^* \hat{y} -
\beta A^* A (x^k - \hat{x}) + \frac{\beta}{\tau} (x^k - \hat{x}) \in
\partial \|\hat{x}\|_1$. Further considering $A^* \tilde{y}\in
\partial\|\tilde{x}\|_1$ and the convexity of $\|\cdot\|_1$, it can be
shown that
\begin{eqnarray}\label{opt&iter-x}
(\hat{x} - \tilde{x})^* \left(A^* (\hat{y} - \tilde{y}) - \beta A^*
A(x^k-\hat{x}) + \frac{\beta}{\tau}(x^k - \hat{x})\right)\geq 0.
\end{eqnarray}
 From
$A\tilde{x}+\tilde{r}=b$ and $\beta (A\hat{x} + \hat{r} - b) = y^k -
\hat{y}$, the addition of \eqref{opt&iter-s} and \eqref{opt&iter-x}
gives
\begin{eqnarray}\label{key-ineqaulity-0}
\frac{1}{\beta}(\hat{y} - \tilde{y})^* (y^k - \hat{y}) +
\frac{\beta}{\tau}(\hat{x}-\tilde{x})^*(x^k - \hat{x})\geq (y^k -
\hat{y})^* A(x^k - \hat{x}).
\end{eqnarray}
Let $I_n$ be the identity matrix of order $n$. For convenience, we
define
\begin{eqnarray}\label{def:G}
G_0=\left(
  \begin{array}{cc}
 I_n & \\
& \gamma I_m
\end{array}
\right), \; G_1=\left(
  \begin{array}{cc}
\frac{\beta}{\tau}I_n & \\
& \frac{1}{\beta}I_m
\end{array}
\right), \; G=\left(
  \begin{array}{cc}
\frac{\beta}{\tau}I_n & \\
& \frac{1}{\beta\gamma}I_m
\end{array}
\right)  \text{ and } u=\left(
  \begin{array}{c}
x\\
y
\end{array}
\right).
\end{eqnarray}
By using this notation and considering equality $\hat{u}-\tilde{u} =
(\hat{u}-u^k) + (u^k-\tilde{u})$,  \eqref{key-ineqaulity-0} implies
\begin{eqnarray}\label{key-ineqaulity}
(u^k - \tilde{u})^*G_1(u^k-\hat{u}) \geq \|u^k - \hat{u}\|_{G_1}^2 +
(y^k-\hat{y})^*A(x^k-\hat{x}).
\end{eqnarray}
The iteration of $u$ in \eqref{alg-primal} can be written as $
u^{k+1} = u^k - G_0(u^k - \hat u)$, from which it can be shown
\begin{eqnarray}\label{convergence}
\nonumber\|u^k-\tilde{u}\|_G^2 - \|u^{k+1}-\tilde{u}\|_G^2 &=&
2(u^k-\tilde{u})^*G_1(u^k-\hat u)
- \|G_0(u^k-\hat u)\|_G^2\\
\nonumber \text{  (from \eqref{key-ineqaulity})  } &\geq& 2\|u^k -
\hat{u}\|_{G_1}^2 + 2(y^k-\hat{y})^*A(x^k-\hat{x}) -
\|u^k-\hat u\|_{G_0GG_0}^2 \\
\text{  (from \eqref{def:G})  } &=&\frac{\beta}{\tau}\|x^k-\hat
x\|^2 + \frac{2-\gamma}{\beta}\|y^k -\hat y \|^2 +
2(y^k-\hat{y})^*A(x^k-\hat{x})
\end{eqnarray}
From condition $\tau\lambda_{\max} + \gamma <2$, it holds that
$\delta\triangleq1-\tau\lambda_{\max}/(2-\gamma)>0$. Let
$\rho\triangleq(2-\gamma)/(\beta+\beta\delta)>0$, inequality
\eqref{convergence} implies
\begin{eqnarray}\label{convergence-2}
\nonumber\|u^k-\tilde{u}\|_G^2 - \|u^{k+1}-\tilde{u}\|_G^2
&\geq&\frac{\beta}{\tau}\|x^k-\hat x\|^2 +
\frac{2-\gamma}{\beta}\|y^k
-\hat y \|^2 - \rho \|y^k-\hat{y}\|^2 - \frac{1}{\rho}\|A(x^k-\hat{x})\|^2\\
\nonumber&\geq&\left(\frac{\beta}{\tau}-
\frac{\lambda_{\max}}{\rho}\right)\|x^k-\hat x\|^2 +
\left(\frac{2-\gamma}{\beta}-\rho\right)\|y^k -\hat y \|^2\\
\nonumber \text{  (from definitions of $\delta$ and $\rho$)  } &=&
\frac{\beta\delta^2}{\tau}\|x^k-\hat x\|^2 +
\frac{2-\gamma}{\beta}\frac{\delta}{1+\delta}\|y^k -\hat y \|^2
\\
\text{  (from definitions of $\hat x$, $\hat y$ and $G$) } &=&
\eta\|u^k- u^{k+1}\|_G^2,
\end{eqnarray}
where $\eta\triangleq\min\left(\delta^2,
\frac{\delta(2-\gamma)}{\gamma(1+\delta)}\right)>0$. It follows from
\eqref{convergence-2} that \vspace{.1cm}
\begin{enumerate}
  \item [(a)] $\|u^k-u^{k+1}\|_G$ converges to 0, and thus
$\lim_{k\rightarrow\infty}Ax^k + r^k = b$;
  \item [(b)] $\{u^k\}$ lies in a compact region;
  \item [(c)] $\|u^k-\tilde{u}\|_G^2$ is monotonically non-increasing and thus
converges.
\end{enumerate}
From (b), $\{u^k\}$ has a subsequence $\{u^{k_j}\}$ converges to
$u^\star=(x^\star;y^\star)$. From $y^{k_j}\rightarrow y^{\star}$,
$\lim_{k\rightarrow\infty}Ax^k + r^k = b$ and the iterative formula
for $y$, we have $y^k\rightarrow y^\star$. From (c),
$\|u^k-\tilde{u}\|_G^2\rightarrow\|u^\star-\tilde{u}\|_G^2$, which,
by further considering $y^k\rightarrow y^\star$, implies that any
limit point of $\{x^k\}$, if more than one, must have an equal
distance to $\tilde{x}$. For convenience, we let $z(x,y)\triangleq
x-\tau A^*(Ax-b-y/\beta)/(1+\mu\beta)$. By eliminating $r^{k+1}$,
the second equation in \eqref{alg-primal} becomes $x^{k+1} =
\text{Shrink}(z(x^k,y^k), \tau/\beta)$, which implies
\begin{eqnarray}
\nonumber x^{k_j+1} &=&
\text{Shrink}\left(z(x^{k_j},y^{k_j}),\tau/\beta\right) \rightarrow
\text{Shrink} \left(z(x^\star,y^\star),\tau/\beta\right)\triangleq
x^{\star\star}.
\end{eqnarray}
Thus, $x^{\star\star}$ is also a limit point of $\{x^k\}$ and must
have an equal distance to $\tilde{x}$ as $x^\star$ does, i.e.,
\begin{eqnarray}\label{key-conv}
\|x^\star-\tilde{x}\| = \|x^{\star\star} - \tilde{x}\| =
\|\text{Shrink} \left(z(x^\star,y^\star),\tau/\beta\right) -
\text{Shrink} \left(z(\tilde{x},\tilde{y}),\tau/\beta\right)\|,
\end{eqnarray}
where the second equality is because  $\tilde{x}=\text{Shrink}
\left(z(\tilde{x},\tilde{y}),\tau/\beta\right)$.  From the property
of Shrink$(\cdot,\tau/\beta)$, equality \eqref{key-conv} implies
(see e.g., \cite{FPC})
\[
x^\star - \tilde{x}  =\text{Shrink}
\left(z(x^\star,y^\star),\tau/\beta\right) - \text{Shrink}
\left(z(\tilde{x},\tilde{y}),\tau/\beta\right) = \text{Shrink}
\left(z(x^\star,y^\star),\tau/\beta\right) - \tilde{x}.
\]
 Thus, $x^\star = \text{Shrink}
\left(z(x^\star,y^\star),\tau/\beta\right)$. By letting $r^\star = b
- Ax^\star$,  it is easy to show from the above discussions  that
$(r^\star,x^\star,y^\star)$ is a solution to \eqref{decoder-L1L2-2}.
By letting $\tilde{u}=(\tilde{x},\tilde{y}) =
(x^\star,y^\star)=u^\star$ and $\tilde{r} = r^\star$ at the
beginning, we get the convergence of $\{u^k\}$ and thus that of
$\{r^k,x^k,y^k\}$.
\end{proof}

\end{document}